\documentstyle[12pt, leqno, amstex]{amsart}

\makeatletter

\begin{document}
\bibliographystyle{amsalpha}
\newcommand\lra{\longrightarrow}
\newcommand\ra{\rightarrow}
\newtheorem{Assumption}{Assumption}[section]
\newtheorem{Theorem}{Theorem}[section]
\newtheorem{Lemma}{Lemma}[section]
\newtheorem{Remark}{Remark}[section]
\newtheorem{Corollary}{Corollary}[section]
\newtheorem{Conjecture}{Conjecture}[section]
\newtheorem{Proposition}{Proposition}[section]
\newtheorem{Example}{Example}[section]
\newtheorem{Definition}{Definition}[section]
\newtheorem{Proof}{Proof}[section]
\newtheorem{Problem}{Problem}[section]
\newtheorem{Question}{Question}[section]
\renewcommand{\thesubsection}{\it}

\baselineskip=14pt
\addtocounter{section}{-1}

\title[Symmetry in the Painlev\'e systems]{Symmetry in the Painlev\'e systems and their extensions to four-dimensional systems \\}

\author{Yusuke Sasano }
\maketitle

\begin{abstract}
We give a new approach to the symmetries of the Painlev\'e equations $P_{V},P_{IV},P_{III}$ and $P_{II}$, respectively. Moreover, we make natural extensions to fourth-order analogues for each of the Painlev\'e equations $P_{V}$ and $P_{III}$, respectively, which are natural in the sense that they preserve the symmetries.
\end{abstract}

\section{Introduction}
This is the third paper in a series of four papers (see \cite{Sasa5,Sasa6}), aimed at giving a complete study of the following problem:

\begin{Problem}
For each affine root system $A$ with affine Weyl group $W(A)$, find a system of differential equations for which $W(A)$ acts as its B{\"a}cklund transformations.
\end{Problem}

At first, let us summarize the results obtained up to now in the following list.
\begin{center}
\begin{tabular}{|c|c|c|c|} \hline
Type & System  & Dimension & References  \\ \hline
$A_4^{(1)}$ & 2-coupled Painlev\'e IV system & 4 & \cite{N1,N2,Sasa5}  \\ \hline
$A_5^{(1)}$ & 2-coupled Painlev\'e V system & 4 & \cite{N1,N2,Sasa5}  \\ \hline
$D_5^{(1)}$ & 2-coupled Painlev\'e V system & 4 & \cite{Sasa6}  \\ \hline
$D_6^{(1)}$ & 2-coupled Painlev\'e VI system & 4 &  \cite{Sasa5} \\ \hline
$B_4^{(1)}$ & 2-coupled Painlev\'e III system & 4 & \cite{Sasa6} \\ \hline
\end{tabular}
\end{center}
Each of them is a family of coupled Painlev\'e systems with the following form for its Hamiltonian:
\begin{align}
\begin{split}
&\frac{dx}{dt}=\frac{\partial H}{\partial y}, \quad \frac{dy}{dt}=-\frac{\partial H}{\partial x}, \quad \frac{dz}{dt}=\frac{\partial H}{\partial w}, \quad \frac{dw}{dt}=-\frac{\partial H}{\partial z},\\
&H(x,y,z,w,t;\alpha_0,\alpha_1,\dots,\beta_0,\beta_1, \dots)\\
&=H_{*}(x,y,t;\alpha_0,\alpha_1,\dots)+H_{*}(z,w,t;\beta_0,\beta_1,\dots)+R\\
&(*=VI,V,IV,III).
\end{split}
\end{align}
Here the symbol $R$ denotes the interaction term for each system.

Our idea is to find a system in the following way:
\begin{enumerate}
\item We make a set of invariant divisors given by connecting two copies of them given in the case of the Painlev\'e systems by adding the term with invariant divisor $x-z$.

\item We make the symmetry associated with a set of invariant divisors given by 1.

\item We make the holomorphy conditions $r_i$ associated with the symmetry in 2.

\item We look for a polynomial Hamiltonian system with the holomorphy conditions given by 3.
\end{enumerate}
The crucial idea of this work is to use the holomorphy characterization of each system, which can be considered as a generalization of Takano's theory \cite{M,T1}.

In the next stage, following the above results, we try to seek a system with $W(B_3^{(1)})$-symmetry. At first, one might try to seek a system in dimension four with its symmetry. However, such a system can not be obtained. In this paper, we will change our viewpoint, by seeking a system not in dimension four but in dimension two. In dimension two, it is well-known that the Painlev\'e systems $P_J, \ (J=VI,V,IV,III,II,I)$ have the affine Weyl group symmetries explicitly given in the following table.
\begin{center}
\begin{tabular}{|c|c|c|c|c|c|c|} \hline
$P_J$ & $P_I$ & $P_{II}$  & $P_{IIII}$ & $P_{IV}$ & $P_{V}$ & $P_{VI}$  \\ \hline
Symmetry & none & $W(A_1^{(1)})$  & $W(C_2^{(1)})$ & $W(A_2^{(1)})$ & $W(A_3^{(1)})$ & $W(D_4^{(1)})$  \\ \hline
\end{tabular}
\end{center}

This paper is the stage in this project where we find a new viewpoint for the symmetries of the Painlev\'e equations $P_{V},P_{IV},P_{III}$ and $P_{II}$, that is, we will show that each of the Painlev\'e equations $P_{V},P_{IV},P_{III}$ and $P_{II}$ has hidden affine Weyl group symmetry of types $B_3^{(1)},G_2^{(1)},D_3^{(2)}$ and $A_2^{(2)}$, respectively. We seek these symmetries for a Hamiltonian system in charts other than the original chart in each space of initial conditions constructed by K. Okamoto. In other charts, we can find hidden symmetries different from the ones in the original charts.
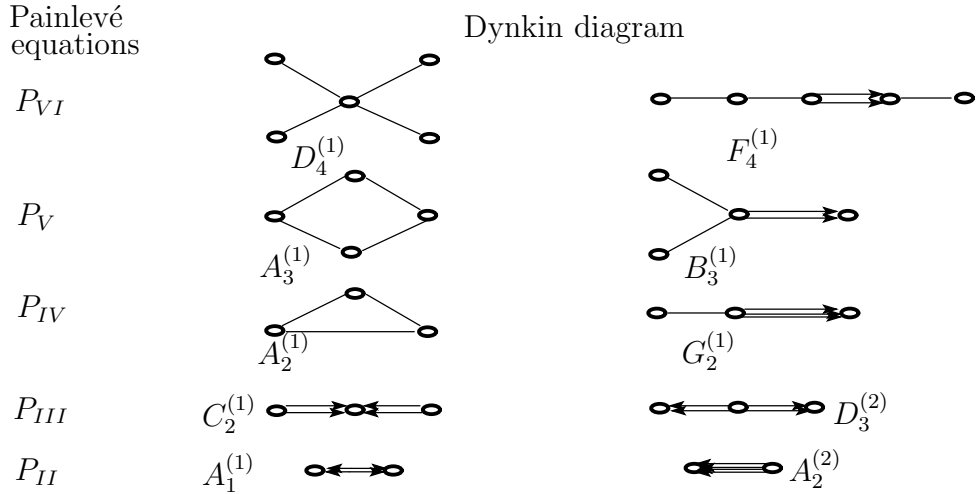
\begin{figure}[h]
\unitlength 0.1in
\begin{picture}(50.45,25.15)(16.00,-28.35)
%
\special{pn 8}%
\special{pa 5570 2777}%
\special{pa 5200 2777}%
\special{fp}%
\special{sh 1}%
\special{pa 5200 2777}%
\special{pa 5267 2797}%
\special{pa 5253 2777}%
\special{pa 5267 2757}%
\special{pa 5200 2777}%
\special{fp}%
%
\special{pn 20}%
\special{ar 3380 882 45 24  0.0000000 6.2831853}%
%
\special{pn 20}%
\special{ar 2990 657 45 23  0.0000000 6.2831853}%
%
\special{pn 20}%
\special{ar 3800 662 45 23  0.0000000 6.2831853}%
%
\special{pn 20}%
\special{ar 3000 1066 45 24  0.0000000 6.2831853}%
%
\special{pn 20}%
\special{ar 3800 1072 45 22  0.0000000 6.2831853}%
%
\special{pn 20}%
\special{ar 5010 867 45 23  0.0000000 6.2831853}%
%
\special{pn 20}%
\special{ar 5410 867 45 23  0.0000000 6.2831853}%
%
\special{pn 20}%
\special{ar 5800 867 45 23  0.0000000 6.2831853}%
%
\special{pn 20}%
\special{ar 6210 867 45 23  0.0000000 6.2831853}%
%
\special{pn 20}%
\special{ar 6600 862 45 23  0.0000000 6.2831853}%
%
\special{pn 20}%
\special{ar 3410 1271 45 23  0.0000000 6.2831853}%
%
\special{pn 20}%
\special{ar 2990 1482 45 22  0.0000000 6.2831853}%
%
\special{pn 20}%
\special{ar 3390 1670 45 24  0.0000000 6.2831853}%
%
\special{pn 20}%
\special{ar 3790 1476 45 23  0.0000000 6.2831853}%
%
\special{pn 20}%
\special{ar 5000 1266 45 24  0.0000000 6.2831853}%
%
\special{pn 20}%
\special{ar 5000 1681 45 23  0.0000000 6.2831853}%
%
\special{pn 20}%
\special{ar 5420 1471 45 23  0.0000000 6.2831853}%
%
\special{pn 20}%
\special{ar 5990 1476 45 23  0.0000000 6.2831853}%
%
\special{pn 20}%
\special{ar 3410 1886 45 23  0.0000000 6.2831853}%
%
\special{pn 20}%
\special{ar 2990 2080 45 23  0.0000000 6.2831853}%
%
\special{pn 20}%
\special{ar 3790 2086 45 23  0.0000000 6.2831853}%
%
\special{pn 20}%
\special{ar 4990 1988 45 23  0.0000000 6.2831853}%
%
\special{pn 20}%
\special{ar 5400 1988 45 23  0.0000000 6.2831853}%
%
\special{pn 20}%
\special{ar 6000 1988 45 23  0.0000000 6.2831853}%
%
\special{pn 20}%
\special{ar 3000 2500 45 23  0.0000000 6.2831853}%
%
\special{pn 20}%
\special{ar 3410 2495 45 23  0.0000000 6.2831853}%
%
\special{pn 20}%
\special{ar 3810 2495 45 23  0.0000000 6.2831853}%
%
\special{pn 20}%
\special{ar 5005 2487 45 23  0.0000000 6.2831853}%
%
\special{pn 20}%
\special{ar 5415 2482 45 24  0.0000000 6.2831853}%
%
\special{pn 20}%
\special{ar 5815 2482 45 24  0.0000000 6.2831853}%
%
\special{pn 20}%
\special{ar 3200 2813 45 22  0.0000000 6.2831853}%
%
\special{pn 20}%
\special{ar 3610 2813 45 22  0.0000000 6.2831853}%
%
\special{pn 20}%
\special{ar 5185 2800 45 22  0.0000000 6.2831853}%
%
\special{pn 20}%
\special{ar 5595 2800 45 22  0.0000000 6.2831853}%
\put(16.3000,-9.2800){\makebox(0,0)[lb]{$P_{VI}$}}%
\put(16.4000,-15.2700){\makebox(0,0)[lb]{$P_{V}$}}%
\put(16.4000,-20.0900){\makebox(0,0)[lb]{$P_{IV}$}}%
\put(16.2000,-25.4100){\makebox(0,0)[lb]{$P_{III}$}}%
\put(16.3000,-28.6400){\makebox(0,0)[lb]{$P_{II}$}}%
\put(39.8000,-5.4400){\makebox(0,0)[lb]{Dynkin diagram}}%
%
\special{pn 8}%
\special{pa 3030 682}%
\special{pa 3330 857}%
\special{fp}%
%
\special{pn 8}%
\special{pa 3760 688}%
\special{pa 3420 862}%
\special{fp}%
%
\special{pn 8}%
\special{pa 3030 1046}%
\special{pa 3340 898}%
\special{fp}%
%
\special{pn 8}%
\special{pa 3420 902}%
\special{pa 3750 1051}%
\special{fp}%
%
\special{pn 8}%
\special{pa 5060 862}%
\special{pa 5360 862}%
\special{fp}%
%
\special{pn 8}%
\special{pa 5470 862}%
\special{pa 5750 862}%
\special{fp}%
%
\special{pn 8}%
\special{pa 6270 862}%
\special{pa 6530 862}%
\special{fp}%
%
\special{pn 8}%
\special{pa 3360 1271}%
\special{pa 3020 1461}%
\special{fp}%
%
\special{pn 8}%
\special{pa 3010 1502}%
\special{pa 3330 1650}%
\special{fp}%
%
\special{pn 8}%
\special{pa 3470 1282}%
\special{pa 3750 1446}%
\special{fp}%
%
\special{pn 8}%
\special{pa 3450 1661}%
\special{pa 3790 1502}%
\special{fp}%
%
\special{pn 8}%
\special{pa 5050 1282}%
\special{pa 5360 1450}%
\special{fp}%
%
\special{pn 8}%
\special{pa 5050 1666}%
\special{pa 5380 1486}%
\special{fp}%
%
\special{pn 8}%
\special{pa 3360 1891}%
\special{pa 3030 2054}%
\special{fp}%
%
\special{pn 8}%
\special{pa 3050 2086}%
\special{pa 3720 2086}%
\special{fp}%
%
\special{pn 8}%
\special{pa 3470 1891}%
\special{pa 3750 2054}%
\special{fp}%
%
\special{pn 8}%
\special{pa 5050 1988}%
\special{pa 5350 1988}%
\special{fp}%
%
\special{pn 8}%
\special{pa 5830 842}%
\special{pa 6170 842}%
\special{fp}%
\special{sh 1}%
\special{pa 6170 842}%
\special{pa 6103 822}%
\special{pa 6117 842}%
\special{pa 6103 862}%
\special{pa 6170 842}%
\special{fp}%
%
\special{pn 8}%
\special{pa 5830 887}%
\special{pa 6170 887}%
\special{fp}%
\special{sh 1}%
\special{pa 6170 887}%
\special{pa 6103 867}%
\special{pa 6117 887}%
\special{pa 6103 907}%
\special{pa 6170 887}%
\special{fp}%
%
\special{pn 8}%
\special{pa 5460 1456}%
\special{pa 5930 1456}%
\special{fp}%
\special{sh 1}%
\special{pa 5930 1456}%
\special{pa 5863 1436}%
\special{pa 5877 1456}%
\special{pa 5863 1476}%
\special{pa 5930 1456}%
\special{fp}%
%
\special{pn 8}%
\special{pa 5450 1491}%
\special{pa 5930 1491}%
\special{fp}%
\special{sh 1}%
\special{pa 5930 1491}%
\special{pa 5863 1471}%
\special{pa 5877 1491}%
\special{pa 5863 1511}%
\special{pa 5930 1491}%
\special{fp}%
%
\special{pn 8}%
\special{pa 5450 1968}%
\special{pa 5930 1968}%
\special{fp}%
\special{sh 1}%
\special{pa 5930 1968}%
\special{pa 5863 1948}%
\special{pa 5877 1968}%
\special{pa 5863 1988}%
\special{pa 5930 1968}%
\special{fp}%
%
\special{pn 8}%
\special{pa 5450 1994}%
\special{pa 5950 1994}%
\special{fp}%
\special{sh 1}%
\special{pa 5950 1994}%
\special{pa 5883 1974}%
\special{pa 5897 1994}%
\special{pa 5883 2014}%
\special{pa 5950 1994}%
\special{fp}%
%
\special{pn 8}%
\special{pa 5420 2009}%
\special{pa 5950 2009}%
\special{fp}%
\special{sh 1}%
\special{pa 5950 2009}%
\special{pa 5883 1989}%
\special{pa 5897 2009}%
\special{pa 5883 2029}%
\special{pa 5950 2009}%
\special{fp}%
%
\special{pn 8}%
\special{pa 3750 2474}%
\special{pa 3460 2474}%
\special{fp}%
\special{sh 1}%
\special{pa 3460 2474}%
\special{pa 3527 2494}%
\special{pa 3513 2474}%
\special{pa 3527 2454}%
\special{pa 3460 2474}%
\special{fp}%
%
\special{pn 8}%
\special{pa 3760 2510}%
\special{pa 3460 2510}%
\special{fp}%
\special{sh 1}%
\special{pa 3460 2510}%
\special{pa 3527 2530}%
\special{pa 3513 2510}%
\special{pa 3527 2490}%
\special{pa 3460 2510}%
\special{fp}%
%
\special{pn 8}%
\special{pa 3050 2474}%
\special{pa 3360 2474}%
\special{fp}%
\special{sh 1}%
\special{pa 3360 2474}%
\special{pa 3293 2454}%
\special{pa 3307 2474}%
\special{pa 3293 2494}%
\special{pa 3360 2474}%
\special{fp}%
%
\special{pn 8}%
\special{pa 3040 2510}%
\special{pa 3350 2510}%
\special{fp}%
\special{sh 1}%
\special{pa 3350 2510}%
\special{pa 3283 2490}%
\special{pa 3297 2510}%
\special{pa 3283 2530}%
\special{pa 3350 2510}%
\special{fp}%
%
\special{pn 8}%
\special{pa 5370 2474}%
\special{pa 5060 2474}%
\special{fp}%
\special{sh 1}%
\special{pa 5060 2474}%
\special{pa 5127 2494}%
\special{pa 5113 2474}%
\special{pa 5127 2454}%
\special{pa 5060 2474}%
\special{fp}%
%
\special{pn 8}%
\special{pa 5380 2500}%
\special{pa 5070 2500}%
\special{fp}%
\special{sh 1}%
\special{pa 5070 2500}%
\special{pa 5137 2520}%
\special{pa 5123 2500}%
\special{pa 5137 2480}%
\special{pa 5070 2500}%
\special{fp}%
%
\special{pn 8}%
\special{pa 5470 2474}%
\special{pa 5770 2474}%
\special{fp}%
\special{sh 1}%
\special{pa 5770 2474}%
\special{pa 5703 2454}%
\special{pa 5717 2474}%
\special{pa 5703 2494}%
\special{pa 5770 2474}%
\special{fp}%
%
\special{pn 8}%
\special{pa 5450 2500}%
\special{pa 5790 2500}%
\special{fp}%
\special{sh 1}%
\special{pa 5790 2500}%
\special{pa 5723 2480}%
\special{pa 5737 2500}%
\special{pa 5723 2520}%
\special{pa 5790 2500}%
\special{fp}%
%
\special{pn 8}%
\special{pa 3390 2802}%
\special{pa 3270 2802}%
\special{fp}%
\special{sh 1}%
\special{pa 3270 2802}%
\special{pa 3337 2822}%
\special{pa 3323 2802}%
\special{pa 3337 2782}%
\special{pa 3270 2802}%
\special{fp}%
%
\special{pn 8}%
\special{pa 3390 2802}%
\special{pa 3570 2802}%
\special{fp}%
\special{sh 1}%
\special{pa 3570 2802}%
\special{pa 3503 2782}%
\special{pa 3517 2802}%
\special{pa 3503 2822}%
\special{pa 3570 2802}%
\special{fp}%
%
\special{pn 8}%
\special{pa 3380 2820}%
\special{pa 3260 2820}%
\special{fp}%
\special{sh 1}%
\special{pa 3260 2820}%
\special{pa 3327 2840}%
\special{pa 3313 2820}%
\special{pa 3327 2800}%
\special{pa 3260 2820}%
\special{fp}%
%
\special{pn 8}%
\special{pa 3380 2820}%
\special{pa 3560 2820}%
\special{fp}%
\special{sh 1}%
\special{pa 3560 2820}%
\special{pa 3493 2800}%
\special{pa 3507 2820}%
\special{pa 3493 2840}%
\special{pa 3560 2820}%
\special{fp}%
%
\special{pn 8}%
\special{pa 5570 2822}%
\special{pa 5200 2822}%
\special{fp}%
\special{sh 1}%
\special{pa 5200 2822}%
\special{pa 5267 2842}%
\special{pa 5253 2822}%
\special{pa 5267 2802}%
\special{pa 5200 2822}%
\special{fp}%
%
\special{pn 8}%
\special{pa 5560 2797}%
\special{pa 5220 2797}%
\special{fp}%
\special{sh 1}%
\special{pa 5220 2797}%
\special{pa 5287 2817}%
\special{pa 5273 2797}%
\special{pa 5287 2777}%
\special{pa 5220 2797}%
\special{fp}%
%
\special{pn 8}%
\special{pa 5550 2807}%
\special{pa 5210 2807}%
\special{fp}%
\special{sh 1}%
\special{pa 5210 2807}%
\special{pa 5277 2827}%
\special{pa 5263 2807}%
\special{pa 5277 2787}%
\special{pa 5210 2807}%
\special{fp}%
\put(30.7000,-12.3000){\makebox(0,0)[lb]{$D_4^{(1)}$}}%
\put(53.5000,-11.9000){\makebox(0,0)[lb]{$F_4^{(1)}$}}%
\put(29.1000,-17.9800){\makebox(0,0)[lb]{$A_3^{(1)}$}}%
\put(29.0000,-22.5400){\makebox(0,0)[lb]{$A_2^{(1)}$}}%
\put(26.1000,-28.9000){\makebox(0,0)[lb]{$A_1^{(1)}$}}%
\put(56.8000,-28.7000){\makebox(0,0)[lb]{$A_2^{(2)}$}}%
\put(51.3000,-18.0400){\makebox(0,0)[lb]{$B_3^{(1)}$}}%
\put(51.2000,-22.5000){\makebox(0,0)[lb]{$G_2^{(1)}$}}%
\put(59.1000,-25.7000){\makebox(0,0)[lb]{$D_3^{(2)}$}}%
\put(26.1000,-25.8000){\makebox(0,0)[lb]{$C_2^{(1)}$}}%
\put(16.0000,-4.9000){\makebox(0,0)[lb]{Painlev\'e}}%
\put(16.1000,-6.3000){\makebox(0,0)[lb]{equations}}%
\end{picture}%
\label{fig:LIST.bak}
\caption{Symmetries of the Painlev\'e equations}
\end{figure}
Furthermore, in the case of dimension four we make natural extensions for each of the Painlev\'e equations $P_{V}$ and $P_{III}$, natural in the sense that they preserve the symmetries.

This paper is organized as follows. In Sections 1 through 4, we present two-dimensional polynomial Hamiltonian systems with $W(B_3^{(1)}),W(G_2^{(1)}),W(D_3^{(2)})$ and $W(A_2^{(2)})$-symmetry, respectively. We will show that each system coincides with the Painlev\'e V (resp. IV,III,II) system. We also give an explicit confluence process from the system of type $D_4^{(1)}$ (resp. $B_3^{(1)},B_3^{(1)},G_2^{(1)}$) to the system of type $B_3^{(1)}$ (resp. $G_2^{(1)},D_3^{(2)},A_2^{(2)}$). In Sections 5 and 6, we present a family of coupled Painlev\'e V (resp. III) systems in dimension four with $W(B_5^{(1)}$) (resp. $W(D_5^{(2)})$)-symmetry. We also show that this system coincides with a family of coupled Painlev\'e V (resp. III) systems in dimension four with $W(D_5^{(1)})$ (resp. $W(B_4^{(1)}$)-symmetry (see \cite{Sasa6}). In the final section, we propose further problems on Problem 0.1.

\section{The system of $B_3^{(1)}$}

In this section, we present a 3-parameter family of two-dimensional polynomial Hamiltonian systems given by
\begin{equation}\label{B3}
  \left\{
  \begin{aligned}
   \frac{dx}{dt} &=\frac{-2x^3y+2x^2y-(\alpha_1+2\alpha_2)x^2-(t-1+2\alpha_3)x+t}{t},\\
   \frac{dy}{dt} &=\frac{3x^2y^2-2xy^2+2(\alpha_1+2\alpha_2)xy+(t-1+2\alpha_3)y+(\alpha_1+\alpha_2)\alpha_2}{t}\\
   \end{aligned}
  \right. 
\end{equation}
with the polynomial Hamiltonian (cf. \cite{Oka,N3})
\begin{align}\label{HB3}
\begin{split}
&H_{B_3^{(1)}}(x,y,t;\alpha_1,\alpha_2,\alpha_3)\\
&=-\frac{-ty+x^3y^2-x^2y^2+(\alpha_1+2\alpha_2)x^2y+(t-1+2\alpha_3)xy+(\alpha_1+\alpha_2)\alpha_2x}{t}.
\end{split}
\end{align}
Here $x$ and $y$ denote unknown complex variables and $\alpha_0,\alpha_1,\alpha_2$ and $\alpha_3$ are complex parameters satisfying the relation:
\begin{equation}
\alpha_0+\alpha_1+2\alpha_2+2\alpha_3=1.
\end{equation}
We note that the Hamiltonian system \eqref{B3} is a polynomial in the canonical variables $x,y$. In this sense, we call the system \eqref{B3} as a polynomial Hamiltonian system.

\begin{Theorem}\label{th:B3}
The system \eqref{B3} admits extended affine Weyl group symmetry of type $B_3^{(1)}$ as the group of its B{\"a}cklund transformations {\rm (cf. \cite{N3}), \rm} whose generators are explicitly given as follows{\rm : \rm}
\begin{figure}[h]
\unitlength 0.1in
\begin{picture}(35.03,16.31)(9.10,-22.27)
%
\special{pn 20}%
\special{ar 1837 863 293 267  0.0000000 6.2831853}%
%
\special{pn 20}%
\special{ar 1850 1960 293 267  0.0000000 6.2831853}%
%
\special{pn 20}%
\special{ar 2800 1380 293 267  0.0000000 6.2831853}%
%
\special{pn 20}%
\special{ar 4120 1410 293 267  0.0000000 6.2831853}%
%
\special{pn 20}%
\special{pa 2100 960}%
\special{pa 2540 1220}%
\special{fp}%
%
\special{pn 20}%
\special{pa 2150 1930}%
\special{pa 2580 1570}%
\special{fp}%
%
\special{pn 20}%
\special{pa 3090 1310}%
\special{pa 3830 1310}%
\special{fp}%
\special{sh 1}%
\special{pa 3830 1310}%
\special{pa 3763 1290}%
\special{pa 3777 1310}%
\special{pa 3763 1330}%
\special{pa 3830 1310}%
\special{fp}%
%
\special{pn 20}%
\special{pa 3080 1480}%
\special{pa 3820 1480}%
\special{fp}%
\special{sh 1}%
\special{pa 3820 1480}%
\special{pa 3753 1460}%
\special{pa 3767 1480}%
\special{pa 3753 1500}%
\special{pa 3820 1480}%
\special{fp}%
\put(15.7000,-9.4000){\makebox(0,0)[lb]{$x-\infty$}}%
\put(16.3000,-20.3000){\makebox(0,0)[lb]{$x-1$}}%
\put(26.8000,-14.9000){\makebox(0,0)[lb]{$y$}}%
%
\special{pn 8}%
\special{pa 910 1400}%
\special{pa 2500 1400}%
\special{dt 0.045}%
\special{pa 2500 1400}%
\special{pa 2499 1400}%
\special{dt 0.045}%
%
\special{pn 20}%
\special{pa 1320 1320}%
\special{pa 1303 1285}%
\special{pa 1285 1252}%
\special{pa 1267 1222}%
\special{pa 1248 1199}%
\special{pa 1228 1183}%
\special{pa 1206 1177}%
\special{pa 1183 1181}%
\special{pa 1160 1194}%
\special{pa 1136 1215}%
\special{pa 1114 1242}%
\special{pa 1092 1275}%
\special{pa 1072 1313}%
\special{pa 1055 1355}%
\special{pa 1041 1399}%
\special{pa 1031 1445}%
\special{pa 1023 1491}%
\special{pa 1020 1535}%
\special{pa 1019 1576}%
\special{pa 1022 1612}%
\special{pa 1029 1642}%
\special{pa 1039 1665}%
\special{pa 1052 1678}%
\special{pa 1069 1680}%
\special{pa 1089 1671}%
\special{pa 1112 1652}%
\special{pa 1137 1627}%
\special{pa 1163 1598}%
\special{pa 1170 1590}%
\special{sp}%
%
\special{pn 20}%
\special{pa 1130 1620}%
\special{pa 1230 1520}%
\special{fp}%
\special{sh 1}%
\special{pa 1230 1520}%
\special{pa 1169 1553}%
\special{pa 1192 1558}%
\special{pa 1197 1581}%
\special{pa 1230 1520}%
\special{fp}%
\put(10.0000,-11.5000){\makebox(0,0)[lb]{$\pi$}}%
\end{picture}%
\label{fig:CPVB53}
\caption{Dynkin diagram of type $B_3^{(1)}$}
\end{figure}
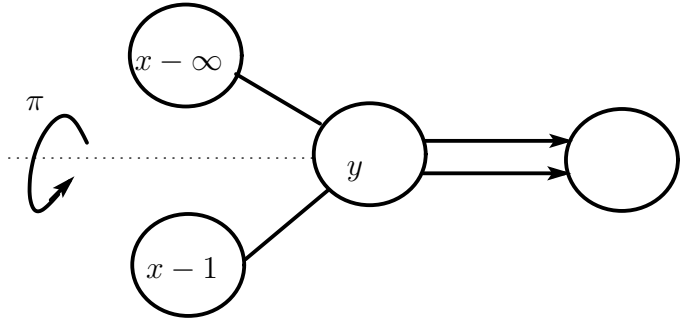
\begin{align}\label{symmetryB3}
\begin{split}
        s_{0}: (x,y,t;\alpha_0,\alpha_1,\alpha_2,\alpha_3) &\rightarrow (x,y-\frac{\alpha_0}{x-1},t;-\alpha_0,\alpha_1,\alpha_2+\alpha_0,\alpha_3), \\
        s_{1}: (x,y,t;\alpha_0,\alpha_1,\alpha_2,\alpha_3) &\rightarrow (x,y,t;\alpha_0,-\alpha_1,\alpha_2+\alpha_1,\alpha_3), \\
        s_{2}: (x,y,t;\alpha_0,\alpha_1,\alpha_2,\alpha_3) &\rightarrow (x+\frac{\alpha_2}{y},y,t;\alpha_0+\alpha_2,\alpha_1+\alpha_2,-\alpha_2,\alpha_3+\alpha_2), \\
        s_{3}: (x,y,t;\alpha_0,\alpha_1,\alpha_2,\alpha_3) &\rightarrow (x,y-\frac{2\alpha_3}{x}+\frac{t}{x^2},-t;\alpha_0,\alpha_1,\alpha_2+2\alpha_3,-\alpha_3), \\
        {\pi}: (x,y,t;\alpha_0,\alpha_1,\alpha_2,\alpha_3) &\rightarrow (\frac{x}{x-1},-(x-1)((x-1)y+\alpha_2),-t;\alpha_1,\alpha_0,\alpha_2,\alpha_3).
\end{split}
\end{align}
\end{Theorem}

\noindent
The list \eqref{symmetryB3} should be read as
\begin{align*}
&s_0(\alpha_0)=-\alpha_0, \quad s_0(\alpha_1)=\alpha_1, \quad s_0(\alpha_2)=\alpha_2+\alpha_0, \quad s_0(\alpha_3)=\alpha_3,\\
&s_0(x)=x, \quad s_0(y)=y-\frac{\alpha_0}{x-1}, \quad s_0(t)=t.
\end{align*}
The above figure denotes the Dynkin diagram of type $B_3^{(1)}$. Let us set
\begin{equation*}
f_0:=x-1, \quad f_1:=x-\infty, \quad f_2:=y.
\end{equation*}
Following \cite{N2}, we define the actions $w_i \ (i=0,1,2)$ as
\begin{equation}
w_i(g)=g+\frac{\alpha_i}{f_i}\{g,f_i\}, \quad g \in {\Bbb C}(t)[x,y],
\end{equation}
where $\{,\}$ is the Poisson bracket such that $\{x,x\}=\{y,y\}=0, \ \{x,y\}=1$. These actions of $w_i \ (i=0,1,2)$ are equivalent to the actions of $s_i \ (i=0,1,2)$ given in Theorem \ref{th:B3}. However, the actions of $s_3$ given in Theorem \ref{th:B3} are different from the actions defined by Noumi and Yamada in \cite{N2}. We also remark that $f:=x$ is not an invariant divisor of the system \eqref{B3}.

In order to prove Theorem \ref{th:B3}, we recall the definition of a symplectic transformation and its properties (see \cite{M,T1}). Let
$$
\varphi : x=x(X,Y,t), \  y=y(X,Y,t), \  t=t
$$
be a biholomorphic mapping from a domain $D$ in ${\Bbb C}^3 \ni (X,Y,t)$ into ${\Bbb C}^3 \ni (x,y,t)$. We say that the mapping is symplectic if
\begin{equation*}
dx \wedge dy= dX \wedge dY,
\end{equation*}
where $t$ is considered as a constant or a parameter, namely, if, for $t=t_0$, ${\varphi}_{t_0}=\varphi|_{t=t_0}$ is a symplectic mapping from the $t_0$-section $D_{t_0}$ of $D$ to $\varphi(D_{t_0})$. Suppose that the mapping is symplectic. Then any Hamiltonian system
$$
dx/dt=\partial H/\partial y, \ \  dy/dt=-\partial H/\partial x
$$
is transformed to
$$
dX/dt=\partial K/\partial Y, \ \  dY/dt=-\partial K/\partial X,
$$
where
\begin{equation*}\label{A}
(A) \ \ \ dx \wedge dy- dH \wedge dt =dX \wedge dY- dK \wedge dt.
\end{equation*}
Here $t$ is considered as a variable. By this equation, the function $K$ is determined by $H$ uniquely modulo functions of $t$, namely, modulo functions independent of $X$ and $Y$.

{\bf Proof of Theorem \ref{th:B3}.}
At first, we consider the case of the transformation $s_0$.  Set
\begin{align*}
\begin{split}
&X:=x, \quad Y:=y-\frac{\alpha_0}{x-1}, \quad T:=t,\\
&A_0:=-\alpha_0, \quad A_1=\alpha_1, \quad A_2=\alpha_2+\alpha_0, \quad A_3=\alpha_3.
\end{split}
\end{align*}
By resolving in $x,y,t,\alpha_0,\dots,\alpha_3$, we obtain
\begin{align}\tag*{$S_0:$}
\begin{split}
&x=X, \quad y=Y-\frac{A_0}{X-1}, \quad t=T,\\
&\alpha_0=-A_0, \quad \alpha_1=A_1, \quad \alpha_2=A_2+A_0, \quad \alpha_3=A_3.
\end{split}
\end{align}
By $S_0$, we obtain the polynomial Hamiltonian $S_0(H_1)$, and we see that
$$
H_1=(S_0(H_1)-A_0)|_{\{X \rightarrow x, \ Y \rightarrow y, \ T \rightarrow t, \ A_0 \rightarrow \alpha_0, \ A_1 \rightarrow \alpha_1, \ A_2 \rightarrow \alpha_2, \ A_3 \rightarrow \alpha_3\}}.
$$
Since $H_1$ is modulo functions of $t$, we can check in the case of $s_0$.

The cases of $s_1,s_2$ are similar. We note the relation between $H_1$ and the transformed Hamiltonian $K_i \ (i=1,2)$, respectively: with the notation $res:=\{X \rightarrow x, \ Y \rightarrow y, \ T \rightarrow t, \ A_0 \rightarrow \alpha_0, \ A_1 \rightarrow \alpha_1, \ A_2 \rightarrow \alpha_2, \ A_3 \rightarrow \alpha_3\}$
\begin{align*}
H_1=&K_1|_{res},\\
H_2=&\left(K_2+\frac{A_2(A_2+2A_3-1+T)}{T}\right)\vert_{res}.
\end{align*}
Next, we consider the case of $s_3$. Setting
\begin{align}\tag*{$S_3:$}
\begin{split}
&x=X, \quad y=Y-\frac{2A_3}{X}+\frac{T}{X^2}, \quad t=-T,\\
&\alpha_0=A_0, \quad \alpha_1=A_1, \quad \alpha_2=A_2+2A_3, \quad \alpha_3=-A_3.
\end{split}
\end{align}
Applying the transformation in $t$ and the transformation of the symplectic 2-form:
\begin{equation*}
dx \wedge dy=dX \wedge dY-d\left(\frac{1}{X}\right) \wedge dT,
\end{equation*}
 we obtain the rational Hamiltonian $S_3(H_1)$, and we see that
$$
H_1=-\{S_3(H_1)+\frac{1}{X}-(A_1+2A_2+\frac{2A_3(T+1)}{T})\}|_{res}.
$$
Then we can check in the case of $s_3$.

The case of $\pi$ is similar. We note the relation between $H_1$ and the transformed Hamiltonian $\varPi(H_1)$ is given as follows:
$$
H_1=-\{\varPi(H_1)-A_2\}|_{res}.
$$

This completes the proof. \qed

Consider the following birational and symplectic transformations $r_i$(cf. \cite{M,T1}):
\begin{align}\label{holomorphyB3}
\begin{split}
r_0&:x_0=-((x-1)y-\alpha_0)y, \ y_0=\frac{1}{y},\\
r_1&:x_1=\frac{1}{x}, \ y_1=-(yx+\alpha_1+\alpha_2)x, \\
r_2&:x_2=\frac{1}{x}, \ y_2=-(yx+\alpha_2)x, \\
r_3&:x_3=x, \ y_3=y-\frac{2\alpha_3}{x}+\frac{t}{x^2}.
\end{split}
\end{align}
These transformations are appeared as the patching data in the space of initial conditions of the system \eqref{B3}. The fact that the space of initial conditions of the system \eqref{B3} is covered by this data will be cleared in the following paper.

Since each transformation $r_i$ is symplectic, the system \eqref{B3} is transformed into a Hamiltonian system, whose Hamiltonian may have poles. It is remarkable that the transformed system becomes again a polynomial system for any $i=0,1,2,3$. Furthermore, this holomorphy property uniquely characterizes the system \eqref{B3}.

\begin{Theorem}\label{th:HoB3}
Let us consider a polynomial Hamiltonian system with Hamiltonian $H \in {\Bbb C}(t)[x,y]$. We assume that

$(A1)$ $deg(H)=5$ with respect to $x,y$.

$(A2)$ This system becomes again a polynomial Hamiltonian system in each coordinate $r_i \ (i=0,1,2,3)$.

\noindent
Then such a system coincides with the system \eqref{B3}.
\end{Theorem}
We remark that if we look for a polynomial Hamiltonian system which admits the symmetry \eqref{symmetryB3}, we must consider cumbersome polynomials in variables $x,y,t,\alpha_i$. On the other hand, in the holomorphy requirement \eqref{holomorphyB3}, we only need to consider polynomials in $x,y$. This reduces the number of unknown coefficients drastically. On relations between $s_i$ and $r_i$, see \cite{Sasa7}.

{\bf Proof of Theorem \ref{th:HoB3}.}
At first, resolving the coordinate $r_0$ in the variables $x,y$, we obtain
\begin{equation}\tag*{($R_0$)}
x=-x_0y_0^2+\alpha_0 y_0+1, \quad y=\frac{1}{y_0}.
\end{equation}
The polynomial $H$ satisfying $(A1)$ has 20 unknown coefficients in ${\Bbb C}(t)$. By $R_0$, we transform $H$ into $R_0(H)$, which has poles in only $y_0$. For $R_0(H)$, we only have to determine the unknown coefficients so that they cancel the poles of $R_0(H)$.

In this way, we can obtain the Hamiltonian $H_1$. \qed

Here we recall that the sixth Painlev\'e system is given by
\begin{equation}\label{PVI}
\frac{dx}{dt}=\frac{\partial H_{VI}}{\partial y}, \ \ \frac{dy}{dt}=-\frac{\partial H_{VI}}{\partial x}
\end{equation}
with the polynomial Hamiltonian
\begin{align}
\begin{split}
&H_{VI}(x,y,t;\alpha_0,\alpha_1,\alpha_2,\alpha_3,\alpha_4)\\
&=\frac{1}{t(t-1)}[y^2(x-t)(x-1)x-\{(\alpha_0-1)(x-1)x+\alpha_3(x-t)x\\
&+\alpha_4(x-t)(x-1)\}y+\alpha_2(\alpha_1+\alpha_2)x]  \quad (\alpha_0+\alpha_1+2\alpha_2+\alpha_3+\alpha_4=1). 
\end{split}
\end{align}
This system \eqref{PVI} admits affine Weyl group symmetry of type $D_4^{(1)}$ as the group of its B{\"a}cklund transformations {\rm (see \cite{N3}), \rm} whose generators are explicitly given as follows{\rm : \rm}
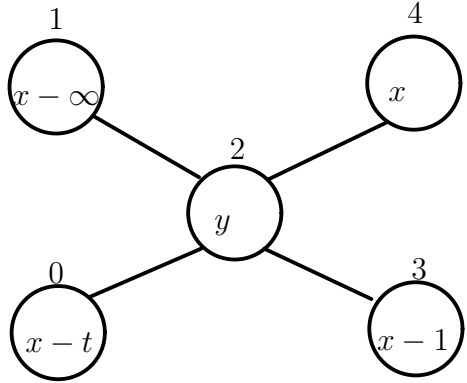
\begin{figure}[ht]
\unitlength 0.1in
\begin{picture}(23.71,20.35)(22.18,-23.65)
%
\special{pn 20}%
\special{ar 2462 834 244 244  0.0000000 6.2831853}%
%
\special{pn 20}%
\special{ar 2473 2121 244 244  0.0000000 6.2831853}%
%
\special{pn 20}%
\special{ar 3397 1494 244 244  0.0000000 6.2831853}%
%
\special{pn 20}%
\special{ar 4334 814 244 244  0.0000000 6.2831853}%
%
\special{pn 20}%
\special{ar 4345 2101 244 244  0.0000000 6.2831853}%
%
\special{pn 20}%
\special{pa 2660 988}%
\special{pa 3210 1307}%
\special{fp}%
%
\special{pn 20}%
\special{pa 2638 1934}%
\special{pa 3221 1681}%
\special{fp}%
%
\special{pn 20}%
\special{pa 3573 1318}%
\special{pa 4189 1021}%
\special{fp}%
%
\special{pn 20}%
\special{pa 3551 1681}%
\special{pa 4112 1945}%
\special{fp}%
\put(23.0000,-22.2000){\makebox(0,0)[lb]{$x-t$}}%
\put(22.3000,-9.3000){\makebox(0,0)[lb]{$x-\infty$}}%
\put(32.9000,-15.8000){\makebox(0,0)[lb]{$y$}}%
\put(42.0000,-9.1000){\makebox(0,0)[lb]{$x$}}%
\put(41.4000,-22.1000){\makebox(0,0)[lb]{$x-1$}}%
\put(24.2000,-18.6000){\makebox(0,0)[lb]{$0$}}%
\put(24.2000,-5.3000){\makebox(0,0)[lb]{$1$}}%
\put(33.7000,-12.1000){\makebox(0,0)[lb]{$2$}}%
\put(43.2000,-18.4000){\makebox(0,0)[lb]{$3$}}%
\put(43.0000,-5.0000){\makebox(0,0)[lb]{$4$}}%
\end{picture}%
\label{fig:PVID4}
\caption{Dynkin diagram of type $D_4^{(1)}$}
\end{figure}
\begin{align*}
        w_{0}: (x,y,t;\alpha_0,\alpha_1,\alpha_2,\alpha_3,\alpha_4) \rightarrow &(x,y-\frac{\alpha_0}{x-t},t;-\alpha_0,\alpha_1,\alpha_2+\alpha_0,\alpha_3,\alpha_4), \\
        w_{1}: (x,y,t;\alpha_0,\alpha_1,\alpha_2,\alpha_3,\alpha_4) \rightarrow &(x,y,t;\alpha_0,-\alpha_1,\alpha_2+\alpha_1,\alpha_3,\alpha_4), \\
        w_{2}: (x,y,t;\alpha_0,\alpha_1,\alpha_2,\alpha_3,\alpha_4) \rightarrow &(x+\frac{\alpha_2}{y},y,t;\\
        &\alpha_0+\alpha_2,\alpha_1+\alpha_2,-\alpha_2,\alpha_3+\alpha_2,\alpha_4+\alpha_2), \\
        w_{3}: (x,y,t;\alpha_0,\alpha_1,\alpha_2,\alpha_3,\alpha_4) \rightarrow &(x,y-\frac{\alpha_3}{x-1},t;\alpha_0,\alpha_1,\alpha_2+\alpha_3,-\alpha_3,\alpha_4), \\
        w_{4}: (x,y,t;\alpha_0,\alpha_1,\alpha_2,\alpha_3,\alpha_4) \rightarrow &(x,y-\frac{\alpha_4}{x},t;\alpha_0,\alpha_1,\alpha_2+\alpha_4,\alpha_3,-\alpha_4).
        \end{align*}

In addition to Theorems \ref{th:B3} and \ref{th:HoB3}, we give an explicit description of a confluence from the sixth Painlev\'e system \eqref{PVI} to the system \eqref{B3}. 
\begin{Theorem}\label{th:D4B3}
For the sixth Painlev\'e system \eqref{PVI}, we make the change of parameters and variables
\begin{gather}
\begin{gathered}\label{PD4B3}
\alpha_0={\varepsilon}^{-1}, \quad \alpha_1=A_0, \quad \alpha_2=A_2, \quad \alpha_3=\frac{2A_3\varepsilon -1}{\varepsilon}, \quad \alpha_4=A_1,
\end{gathered}\\
\begin{gathered}\label{VD4B3}
t=1-{\varepsilon}T, \quad x=\frac{1}{1-X}, \quad y=X^2Y-2XY+A_2X+Y-A_2
\end{gathered}
\end{gather}
from $\alpha_0,\alpha_1,\alpha_2,\alpha_3,\alpha_4,t,x,y$ to $A_0,A_1,A_2,A_3,\varepsilon,T,X,Y$. Then the system \eqref{PVI} can also be written in the new variables $T,X,Y$ and parameters $A_0,A_1,A_2,A_3,\varepsilon$ as a Hamiltonian system. This new system tends to the system \eqref{B3} as $\varepsilon \ra 0$.
\end{Theorem}

By proving the following theorem, we see how the degeneration process given in Theorem \ref{th:D4B3} works on the B{\"a}cklund transformation group $W(D_4^{(1)})$ (cf. \cite{T2}).
\begin{Theorem}\label{th:WD4WB3}
For the degeneration process in Theorem \ref{th:D4B3}, we can choose a subgroup $W_{D_4^{(1)} \rightarrow B_3^{(1)}}$ of the B{\"a}cklund transformation group $W(D_4^{(1)})$ so that $W_{D_4^{(1)} \rightarrow B_3^{(1)}}$ converges to the B{\"a}cklund transformation group $W(B_3^{(1)})$ of the system \eqref{B3} as $\varepsilon \rightarrow 0$.
\end{Theorem}

{\bf Proof of Theorem \ref{th:WD4WB3}.}
Notice that
$$
A_0+A_1+2A_2+2A_3=\alpha_0+\alpha_1+2\alpha_2+\alpha_3+\alpha_4=1
$$
and the change of variables from $(x,y)$ to $(X,Y)$ is symplectic.

Let us see the actions of the generators $w_i, \ i=0,1,2,3,4$ on the parameters $A_i, \ i=0,1,2,3$ and $\varepsilon$ where
\begin{equation*}
A_0=\alpha_1, \quad A_1=\alpha_4, \quad A_2=\alpha_2, \quad A_3=\frac{\alpha_0+\alpha_3}{2}, \quad \varepsilon=\frac{1}{\alpha_0}.
\end{equation*}
By a direct calculation, we have
\begin{align*}
w_0&(A_0,A_1,A_2,A_3,\varepsilon) \longrightarrow (A_0,A_1,A_2+\frac{1}{\varepsilon},A_3-\frac{1}{\varepsilon},-\varepsilon),\\
w_1&(A_0,A_1,A_2,A_3,\varepsilon) \longrightarrow (-A_0,A_1,A_2+A_0,A_3,\varepsilon),\\
w_2&(A_0,A_1,A_2,A_3,\varepsilon) \longrightarrow (A_0+A_2,A_1+A_2,-A_2,A_3+A_2,\frac{\varepsilon}{1+\varepsilon A_2}),\\
w_3&(A_0,A_1,A_2,A_3,\varepsilon) \longrightarrow (A_0,A_1,A_2+2A_3-\frac{1}{\varepsilon},-A_3+\frac{1}{\varepsilon},\varepsilon),\\
w_4&(A_0,A_1,A_2,A_3,\varepsilon) \longrightarrow (A_0,-A_1,A_2+A_1,A_3,\varepsilon).
\end{align*}
We remark that $w_0(A_2),w_0(A_3),w_3(A_2)$ and $w_3(A_3)$ diverge as $\varepsilon \rightarrow 0$.

Observing these relations, we take a subgroup $W_{D_4^{(1)} \rightarrow B_3^{(1)}}$ of $W(D_4^{(1)})$ generated by $S_0,S_1,S_2,S_3$ defined by
$$
S_0:=w_1, \quad S_1:=w_4, \quad S_2:=w_2, \quad S_3:=w_0w_3.
$$
We can easily check
\begin{align}\label{pr1:B3}
\begin{split}
S_0&(A_0,A_1,A_2,A_3,\varepsilon) \longrightarrow (-A_0,A_1,A_2+A_0,A_3,\varepsilon),\\
S_1&(A_0,A_1,A_2,A_3,\varepsilon) \longrightarrow (A_0,-A_1,A_2+A_1,A_3,\varepsilon),\\
S_2&(A_0,A_1,A_2,A_3,\varepsilon) \longrightarrow (A_0+A_2,A_1+A_2,-A_2,A_3+A_2,\frac{\varepsilon}{1+\varepsilon A_2}),\\
S_3&(A_0,A_1,A_2,A_3,\varepsilon) \longrightarrow (A_0,A_1,A_2+2A_3,-A_3,-\varepsilon),
\end{split}
\end{align}
and the generators satisfy the following relations:
$$
(S_i)^2=1, \ (S_0S_1)^2=(S_0S_3)^2=(S_1S_3)^2=1, \ (S_0S_2)^3=(S_1S_2)^3=1, \ (S_2S_3)^4=1.
$$
In short, the group $W_{D_4^{(1)} \rightarrow B_3^{(1)}}=<S_0,S_1,S_2,S_3>$ can be considered to be an affine Weyl group of the affine Lie algebra of type $B_3^{(1)}$ with simple roots $A_0,A_1,A_2,A_3$.

Now we investigate how the generators of $W_{D_4^{(1)} \rightarrow B_3^{(1)}}$ act on $T,X$ and $Y$. We can verify
\begin{align}\label{pr2:B3}
\begin{split}
S_0&(X,Y,T) \longrightarrow (X,Y-\frac{A_0}{X-1},T),\\
S_1&(X,Y,T) \longrightarrow (X,Y,T),\\
S_2&(X,Y,T) \longrightarrow (X+\frac{A_2}{Y},Y,T(1+\varepsilon A_2)),\\
S_3&(X,Y,T) \longrightarrow (X,\frac{(\varepsilon T-1)X^2Y-\varepsilon TXY-2A_3(\varepsilon T-1)X+(2\varepsilon A_3-1)T}{\{(\varepsilon T-1)X-\varepsilon T\}X},-T).
\end{split}
\end{align}
Here $S_3(T)=-T$ can be understood as follows: by using the relation $T=\alpha_0(1-t)$, the action of $S_3$ on $T$ is obtained as
\begin{align*}
S_3(T)=&w_0 \circ w_3(\alpha_0(1-t))\\
=&w_0(\alpha_0(1-t))\\
=&-\alpha_0(1-t)\\
=&-T.
\end{align*}

By comparing \eqref{pr1:B3},\eqref{pr2:B3} with $s_i \ (i=0,1,2,3)$ given in Theorem \ref{th:B3}, we see that our theorem holds. \qed

By the following theorem, we will show that the system \eqref{B3} coincides with the system of type $A_3^{(1)}$ (see \cite{N1,N2,N3}).

\begin{Theorem}\label{th:B3A3}
For the system \eqref{B3}, we make the change of parameters and variables
\begin{gather}
\begin{gathered}\label{PB3A3}
\beta_0=\alpha_2+2\alpha_3,  \quad \beta_1=\alpha_1, \quad \beta_2=\alpha_2, \quad \beta_3=\alpha_0,
\end{gathered}\\
\begin{gathered}\label{VB3A3}
X=\frac{1}{x}, \quad Y=-x(xy+\alpha_2), \quad T=-t
\end{gathered}
\end{gather}
from $\alpha_0,\alpha_1,\alpha_2,\alpha_3,t,x,y$ to $\beta_0,\beta_1,\beta_2,\beta_3,T,X,Y$. Then the system \eqref{B3} can also be written in the new variables $T,X,Y$ and parameters $\beta_0,\beta_1,\beta_2,\beta_3$ as a Hamiltonian system. This new system tends to the system $A_3^{(1)}${\rm : \rm}
\begin{equation}\label{A3}
  \left\{
  \begin{aligned}
   \frac{dX}{dT} &=\frac{2X^2Y}{T}+X^2-\frac{2XY}{T}-(1+\frac{\beta_1+\beta_3}{T})X+\frac{\beta_1}{T},\\
   \frac{dY}{dT} &=-\frac{2XY^2}{T}+\frac{Y^2}{T}-2XY+(1+\frac{\beta_1+\beta_3}{T})Y-\beta_2
   \end{aligned}
  \right. 
\end{equation}
with the Hamiltonian
\begin{align}
\begin{split}
&H_{A_3^{(1)}}(X,Y,T;\beta_1,\beta_2,\beta_3)\\
&=\frac{Y(Y+T)X(X-1)+\beta_2 XT-\beta_3 XY-\beta_1 Y(X-1)}{T}.
\end{split}
\end{align}
\end{Theorem}
By putting $q=1-\frac{1}{X}$, we have the Painlev\'e V equation (see \cite{Tsuda}):
\begin{equation*}
\frac{d^2q}{dT^2}=\left(\frac{1}{2q}+\frac{1}{q-1}\right)\left(\frac{dq}{dT}\right)^2-\frac{1}{T}\frac{dq}{dT}+\frac{(q-1)^2}{T^2}\left(aq+\frac{b}{q}\right)+c\frac{q}{T}+d\frac{q(q+1)}{q-1},
\end{equation*}
where
\begin{equation*}
a=\frac{\beta_1^2}{2}, \quad b=-\frac{\beta_3^2}{2}, \quad c=\beta_0-\beta_2, \quad d=-\frac{1}{2}.
\end{equation*}

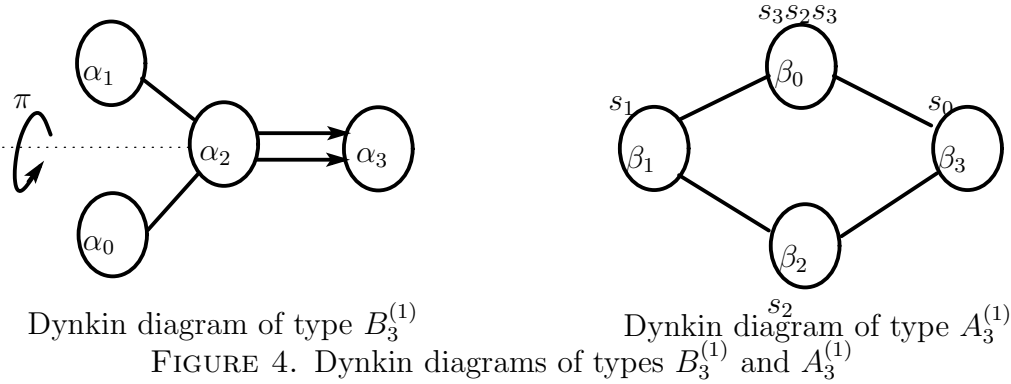
\begin{figure}[ht]
\unitlength 0.1in
\begin{picture}(52.30,16.60)(9.10,-20.70)
%
\special{pn 20}%
\special{ar 1478 814 179 219  0.0000000 6.2831853}%
%
\special{pn 20}%
\special{ar 1486 1712 179 218  0.0000000 6.2831853}%
%
\special{pn 20}%
\special{ar 2068 1237 179 219  0.0000000 6.2831853}%
%
\special{pn 20}%
\special{ar 2876 1262 180 218  0.0000000 6.2831853}%
%
\special{pn 20}%
\special{pa 1639 894}%
\special{pa 1908 1106}%
\special{fp}%
%
\special{pn 20}%
\special{pa 1670 1687}%
\special{pa 1933 1393}%
\special{fp}%
%
\special{pn 20}%
\special{pa 2245 1180}%
\special{pa 2699 1180}%
\special{fp}%
\special{sh 1}%
\special{pa 2699 1180}%
\special{pa 2632 1160}%
\special{pa 2646 1180}%
\special{pa 2632 1200}%
\special{pa 2699 1180}%
\special{fp}%
%
\special{pn 20}%
\special{pa 2239 1319}%
\special{pa 2693 1319}%
\special{fp}%
\special{sh 1}%
\special{pa 2693 1319}%
\special{pa 2626 1299}%
\special{pa 2640 1319}%
\special{pa 2626 1339}%
\special{pa 2693 1319}%
\special{fp}%
\put(19.4000,-13.2000){\makebox(0,0)[lb]{$\alpha_2$}}%
%
\special{pn 8}%
\special{pa 910 1254}%
\special{pa 1884 1254}%
\special{dt 0.045}%
\special{pa 1884 1254}%
\special{pa 1883 1254}%
\special{dt 0.045}%
%
\special{pn 20}%
\special{pa 1161 1188}%
\special{pa 1147 1150}%
\special{pa 1133 1115}%
\special{pa 1117 1089}%
\special{pa 1101 1074}%
\special{pa 1082 1074}%
\special{pa 1063 1088}%
\special{pa 1044 1114}%
\special{pa 1025 1148}%
\special{pa 1008 1189}%
\special{pa 994 1235}%
\special{pa 984 1283}%
\special{pa 977 1331}%
\special{pa 973 1376}%
\special{pa 973 1417}%
\special{pa 977 1450}%
\special{pa 984 1473}%
\special{pa 996 1484}%
\special{pa 1011 1481}%
\special{pa 1030 1463}%
\special{pa 1051 1435}%
\special{pa 1069 1409}%
\special{sp}%
%
\special{pn 20}%
\special{pa 1045 1434}%
\special{pa 1106 1352}%
\special{fp}%
\special{sh 1}%
\special{pa 1106 1352}%
\special{pa 1050 1394}%
\special{pa 1074 1395}%
\special{pa 1082 1417}%
\special{pa 1106 1352}%
\special{fp}%
\put(9.6500,-10.4900){\makebox(0,0)[lb]{$\pi$}}%
\put(13.3000,-17.9000){\makebox(0,0)[lb]{$\alpha_0$}}%
\put(13.3000,-9.0000){\makebox(0,0)[lb]{$\alpha_1$}}%
\put(27.6000,-13.3000){\makebox(0,0)[lb]{$\alpha_3$}}%
%
\special{pn 20}%
\special{ar 5090 830 180 218  0.0000000 6.2831853}%
%
\special{pn 20}%
\special{ar 4320 1260 180 218  0.0000000 6.2831853}%
%
\special{pn 20}%
\special{ar 5110 1770 180 218  0.0000000 6.2831853}%
%
\special{pn 20}%
\special{ar 5960 1260 180 218  0.0000000 6.2831853}%
%
\special{pn 20}%
\special{pa 4920 880}%
\special{pa 4450 1110}%
\special{fp}%
%
\special{pn 20}%
\special{pa 4450 1400}%
\special{pa 4920 1690}%
\special{fp}%
%
\special{pn 20}%
\special{pa 5260 880}%
\special{pa 5770 1140}%
\special{fp}%
%
\special{pn 20}%
\special{pa 5300 1720}%
\special{pa 5800 1390}%
\special{fp}%
\put(48.6000,-5.8000){\makebox(0,0)[lb]{$s_3s_2s_3$}}%
\put(40.8000,-10.6000){\makebox(0,0)[lb]{$s_1$}}%
\put(49.1000,-21.2000){\makebox(0,0)[lb]{$s_2$}}%
\put(57.5000,-10.6000){\makebox(0,0)[lb]{$s_0$}}%
\put(49.5000,-9.1000){\makebox(0,0)[lb]{$\beta_0$}}%
\put(41.7000,-13.6000){\makebox(0,0)[lb]{$\beta_1$}}%
\put(49.7000,-18.8000){\makebox(0,0)[lb]{$\beta_2$}}%
\put(58.0000,-13.6000){\makebox(0,0)[lb]{$\beta_3$}}%
\put(10.4000,-22.3000){\makebox(0,0)[lb]{Dynkin diagram of type $B_3^{(1)}$}}%
\put(41.6000,-22.4000){\makebox(0,0)[lb]{Dynkin diagram of type $A_3^{(1)}$}}%
\end{picture}%
\label{fig:CPVB510}
\caption{Dynkin diagrams of types $B_3^{(1)}$ and $A_3^{(1)}$}
\end{figure}

{\bf Proof of Theorem \ref{th:B3A3}.}
Notice that
\begin{equation*}
\alpha_0+\alpha_1+2\alpha_2+2\alpha_3=\beta_0+\beta_1+\beta_2+\beta_3=1
\end{equation*}
and the change of variables from $(x,y,t)$ to $(X,Y,T)$ in Theorem \ref{th:B3A3} is symplectic. Choose $S_i$ as
\begin{equation*}
S_0:=s_3s_2s_3, \quad S_1:=s_1, \quad S_2:=s_2, \quad S_3:=s_0, \quad S_4:=s_3, \quad S_5:=\pi.
\end{equation*}
The transformations $S_0,S_1,S_2,S_3$ are reflections of
\begin{equation*}
\beta_0=\alpha_2+2\alpha_3,  \quad \beta_1=\alpha_1, \quad \beta_2=\alpha_2, \quad \beta_3=\alpha_0, \ {\rm respectively \rm}.
\end{equation*}
We can verify
\begin{align*}
S_{0}: (X,Y,T;\beta_0,\beta_1,\beta_2,\beta_3) &\rightarrow (X+\frac{\beta_0}{Y+T},Y,T;-\beta_0,\beta_1+\beta_0,\beta_2,\beta_3+\beta_0), \\
S_{1}: (X,Y,T;\beta_0,\beta_1,\beta_2,\beta_3) &\rightarrow (X,Y-\frac{\beta_1}{X},T;\beta_0+\beta_1,-\beta_1,\beta_2+\beta_1,\beta_3), \\
S_{2}: (X,Y,T;\beta_0,\beta_1,\beta_2,\beta_3) &\rightarrow (X+\frac{\beta_2}{Y},Y,T;\beta_0,\beta_1+\beta_2,-\beta_2,\beta_3+\beta_2), \\
S_{3}: (X,Y,T;\beta_0,\beta_1,\beta_2,\beta_3) &\rightarrow (X,Y-\frac{\beta_3}{X-1},T;\beta_0+\beta_3,\beta_1,\beta_2+\beta_3,-\beta_3), \\
S_4: (X,Y,T;\beta_0,\beta_1,\beta_2,\beta_3) &\rightarrow (X,Y+T,-T;\beta_2,\beta_1,\beta_0,\beta_3), \\
S_5: (X,Y,T;\beta_0,\beta_1,\beta_2,\beta_3) &\rightarrow (1-X,-Y,-T;\beta_0,\beta_3,\beta_2,\beta_1).
\end{align*}
The proof has thus been completed. \qed

\begin{Proposition}
The system \eqref{B3} admits the following transformation $\varphi$ as its B{\"a}cklund transformation{\rm:\rm}
\begin{align*}
\begin{split}
\varphi:&(x,y,t;\alpha_0,\alpha_1,\alpha_2,\alpha_3)\\
\ra &(-\frac{t}{x(xy+\alpha_2)},-\frac{x(xy+\alpha_2)(-xy+x^2y-\alpha_0-\alpha_2+\alpha_2x)}{t},t;\\
&\alpha_2+2\alpha_3,\alpha_2,\alpha_0,\frac{\alpha_1-\alpha_0}{2}).
\end{split}
\end{align*}
\end{Proposition}
We note that this transformation $\varphi$ is pulled back the diagram automorphism $\pi$ of the system \eqref{A3}
\begin{equation*}
\pi:(X,Y,T;\beta_0,\beta_1,\beta_2,\beta_3) \ra (-\frac{Y}{T},(X-1)T,T;\beta_1,\beta_2,\beta_3,\beta_0)
\end{equation*}
by transformations \eqref{PB3A3} and \eqref{VB3A3}.

\section{The system of type $G_2^{(1)}$}

In this section, we present a 2-parameter family of two-dimensional polynomial Hamiltonian systems given by
\begin{equation}\label{SG2}
  \left\{
  \begin{aligned}
   \frac{dx}{dt} &=4x^3y+2(\alpha_0+2\alpha_1)x^2+2tx+1,\\
   \frac{dy}{dt} &=-6x^2y^2-4(\alpha_0+2\alpha_1)xy-2ty-2\alpha_1(\alpha_0+\alpha_1)
   \end{aligned}
  \right. 
\end{equation}
with the polynomial Hamiltonian
\begin{align}\label{HG2}
\begin{split}
&H_{G_2^{(1)}}(x,y,t;\alpha_0,\alpha_1,\alpha_2)\\
&=2x^3y^2+2(\alpha_0+2\alpha_1)x^2y+2txy+2\alpha_1(\alpha_0+\alpha_1)x+y.
\end{split}
\end{align}
Here $x$ and $y$ denote unknown complex variables and $\alpha_0,\alpha_1,\alpha_2$ are complex parameters satisfying the relation:
\begin{equation}
\alpha_0+2\alpha_1+3\alpha_2=1.
\end{equation}

\begin{Theorem}\label{th:G2}
The system \eqref{SG2} admits extended affine Weyl group symmetry of type $G_2^{(1)}$ as the group of its B{\"a}cklund transformations {\rm (cf. \cite{N3}), \rm} whose generators are explicitly given as follows{\rm : \rm}
\begin{align*}
        s_{0}: (x,y,t;\alpha_0,\alpha_1,\alpha_2) \rightarrow &(x,y,t;-\alpha_0,\alpha_1+\alpha_0,\alpha_2), \\
        s_{1}: (x,y,t;\alpha_0,\alpha_1,\alpha_2) \rightarrow &(x+\frac{\alpha_1}{y},y,t;\alpha_0+\alpha_1,-\alpha_1,\alpha_2+\alpha_1), \\
        s_{2}: (x,y,t;\alpha_0,\alpha_1,\alpha_2) \rightarrow &(\sqrt{-1}x,-\sqrt{-1}(y-\frac{3\alpha_2}{x}+\frac{t}{x^2}+\frac{1}{2x^3}),-\sqrt{-1}t;\\
        &\alpha_0,\alpha_1+3\alpha_2,-\alpha_2).
        \end{align*}
\end{Theorem}

\begin{figure}
\unitlength 0.1in
\begin{picture}(30.00,6.70)(30.50,-15.90)
%
\special{pn 20}%
\special{ar 3356 1255 306 335  0.0000000 6.2831853}%
%
\special{pn 20}%
\special{ar 4530 1248 286 320  0.0000000 6.2831853}%
%
\special{pn 20}%
\special{ar 5763 1248 287 320  0.0000000 6.2831853}%
\put(30.6700,-13.6100){\makebox(0,0)[lb]{$x-\infty$}}%
\put(43.8100,-13.5400){\makebox(0,0)[lb]{$y$}}%
%
\special{pn 20}%
\special{pa 3668 1267}%
\special{pa 4224 1267}%
\special{fp}%
%
\special{pn 20}%
\special{pa 4791 1079}%
\special{pa 5470 1079}%
\special{fp}%
\special{sh 1}%
\special{pa 5470 1079}%
\special{pa 5403 1059}%
\special{pa 5417 1079}%
\special{pa 5403 1099}%
\special{pa 5470 1079}%
\special{fp}%
%
\special{pn 20}%
\special{pa 4814 1302}%
\special{pa 5448 1302}%
\special{fp}%
\special{sh 1}%
\special{pa 5448 1302}%
\special{pa 5381 1282}%
\special{pa 5395 1302}%
\special{pa 5381 1322}%
\special{pa 5448 1302}%
\special{fp}%
%
\special{pn 20}%
\special{pa 4747 1467}%
\special{pa 5526 1467}%
\special{fp}%
\special{sh 1}%
\special{pa 5526 1467}%
\special{pa 5459 1447}%
\special{pa 5473 1467}%
\special{pa 5459 1487}%
\special{pa 5526 1467}%
\special{fp}%
\end{picture}%
\label{fig:PS4}
\caption{Dynkin diagram of type $G_2^{(1)}$}
\end{figure}
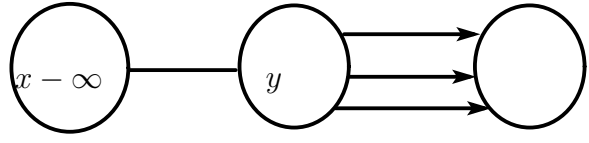

\begin{Theorem}\label{th:HoG2}
Let us consider a polynomial Hamiltonian system with Hamiltonian $H \in {\Bbb C}(t)[x,y]$. We assume that

$(A1)$ $deg(H)=5$ with respect to $x,y$.

$(A2)$ This system becomes again a polynomial Hamiltonian system in each coordinate $r_i \ (i=0,1,2)$ {\rm (cf. \cite{M}): \rm}
\begin{align*}
r_0&:x_0=\frac{1}{x}, \ y_0=-(yx+\alpha_0+\alpha_1)x,\\
r_1&:x_1=\frac{1}{x}, \ y_1=-(yx+\alpha_1)x, \\
r_2&:x_2=x, \ y_2=y-\frac{3\alpha_2}{x}+\frac{t}{x^2}+\frac{1}{2x^3}.
\end{align*}
Then such a system coincides with the system \eqref{SG2}.
\end{Theorem}

Theorems \ref{th:G2}, \ref{th:HoG2} can be cheched by a direct calculation, respectively.

\begin{Theorem}\label{th:SB3toSG2}
For the system \eqref{B3}, we make the change of parameters and variables
\begin{gather}
\begin{gathered}\label{PB3G2}
\alpha_0=-\frac{1}{2{\varepsilon}^2}, \quad \alpha_1=A_0, \quad \alpha_2=A_1, \quad \alpha_3=\frac{1+6A_2{\varepsilon}^2}{4{\varepsilon}^2},
\end{gathered}\\
\begin{gathered}\label{VB3G2}
t=\frac{-1-2{\varepsilon}T}{2e^2}, \quad x=\frac{\varepsilon-X}{\varepsilon}, \quad y=-\varepsilon Y
\end{gathered}
\end{gather}
from $\alpha_0,\alpha_1,\alpha_2,\alpha_3,t,x,y$ to $A_0,A_1,A_2,\varepsilon,T,X,Y$. Then the system \eqref{B3} can also be written in the new variables $T,X,Y$ and parameters $A_0,A_1,A_2,\varepsilon$ as a Hamiltonian system. This new system tends to the system \eqref{SG2} as $\varepsilon \ra 0$.
\end{Theorem}

By proving the following theorem, we see how the degeneration process given in Theorem \ref{th:SB3toSG2} works on the B{\"a}cklund transformation group $W(B_3^{(1)})$ (cf. \cite{T2}).
\begin{Theorem}\label{th:WB3toWG2}
For the degeneration process in Theorem \ref{th:SB3toSG2}, we can choose a subgroup $W_{B_3^{(1)} \rightarrow G_2^{(1)}}$ of the B{\"a}cklund transformation group $W(B_3^{(1)})$ so that $W_{B_3^{(1)} \rightarrow G_2^{(1)}}$ converges to the B{\"a}cklund transformation group $W(G_2^{(1)})$ of the system \eqref{SG2} as $\varepsilon \rightarrow 0$.
\end{Theorem}

{\bf Proof of Theorem \ref{th:WB3toWG2}.}
Notice that
$$
A_0+2A_1+3A_2=\alpha_0+\alpha_1+2\alpha_2+2\alpha_3=1
$$
and the change of variables from $(x,y)$ to $(X,Y)$ is symplectic, however the change of parameters \eqref{PB3G2} is not one to one differently from the case of $P_{VI} \rightarrow P_{V}$.

Choose $S_i \ (i=0,1,2)$ as
$$
S_0:=s_1, \ S_1:=s_2, \ S_2:=s_0s_3,
$$
and set $W_{B_3^{(1)} \rightarrow G_2^{(1)}}=<S_0,S_1,S_2>$. Then we immediately have
\begin{align*}
S_0(A_0,A_1,A_2)&=(-A_0,A_1+A_0,A_2),\\
S_1(A_0,A_1,A_2)&=(A_0+A_1,-A_1,A_2+A_1),\\
S_2(A_0,A_1,A_2)&=(A_0,A_1+3A_2,-A_2).
\end{align*}
However, we see that $S_i(\varepsilon)$ have ambiguities of signature. For example, since
\begin{equation*}
S_2(\varepsilon^2)=s_0 \circ s_3 \left(-\frac{1}{2\alpha_0}\right)=-\frac{1}{2}s_0\left(\frac{1}{\alpha_0}\right)=-\frac{1}{2}\left(-\frac{1}{\alpha_0}\right)=\frac{1}{2\alpha_0}=-{\varepsilon}^2,
\end{equation*}
we can choose any one of the two branches as $S_2(\varepsilon)$. Among such possibilities, we take a choice as
\begin{equation}\label{pr1:G2}
S_0(\varepsilon)=\varepsilon, \quad S_1(\varepsilon)=\varepsilon(1-2A_1 \varepsilon^2)^{-1/2}, \quad S_2(\varepsilon)=\sqrt{-1} \varepsilon,
\end{equation}
where $(1-2A_1 \varepsilon^2)^{-1/2}=1$ at $A_1 \varepsilon^2=0$, or considering in the category of formal power series, we make a convention that $(1-2A_1 \varepsilon^2)^{-1/2}$ is formal power series of $A_1 \varepsilon^2$ with constant term 1 according to
$$
(1+x)^c \sim 1+\sum_{n \geq 1} \binom{c}{n} x^n.
$$
We notice that the generators acting on parameters $A_0,A_1,A_2,\varepsilon$ satisfy the following relations\rm{:\rm}
$$
S_i^2=1, \quad (S_0S_2)^2=1, \quad (S_0S_1)^3=1, \quad (S_1S_2)^6=1.
$$

Now we observe the actions of $S_i, \ i=0,1,2$ on the variables $X,Y,T$. By means of \eqref{VB3G2},\eqref{pr1:G2} and
$$
S_0(t)=s_1(t)=t, \quad S_1(t)=s_2(t)=t, \quad S_2(t)=s_0 \circ s_3(t)=-t,
$$
we can easily check
\begin{equation*}
S_0(T)=T, \quad S_1(T)=(T+A_1 \varepsilon)(1-2A_1 \varepsilon^2)^{-1/2}, \quad S_2(T)=-\sqrt{-1}T.
\end{equation*}
By \eqref{PB3G2},\eqref{VB3G2},\eqref{pr1:G2} and the actions of $s_1,s_2$ on $x,y$, we can easily verify
\begin{align*}
S_0(X,Y)&=(X,Y),\\
S_1(X,Y)&=(X+\frac{A_1}{Y},Y).
\end{align*}
The form of the actions $S_2=s_0s_3$ on $X$ and $Y$ are complicated, but we can see that
\begin{equation*}
S_2(X,Y)=\left(\sqrt{-1}X,-\sqrt{-1}(Y-\frac{3A_2}{X}+\frac{T}{X^2}+\frac{1}{2X^3})\right).
\end{equation*}
The proof has thus been completed. \qed

By the following theorem, we will show that the system \eqref{SG2} coincides with the system of type $A_2^{(1)}$ (see \cite{N1,N2,N3}).
\begin{Theorem}\label{th:G2A2}
For the system \eqref{SG2}, we make the change of parameters and variables
\begin{gather}
\begin{gathered}\label{PG2A2}
\beta_0=\alpha_1+3\alpha_2. \quad \beta_1=\alpha_0, \quad \beta_2=\alpha_1,
\end{gathered}\\
\begin{gathered}\label{VG2A2}
X=\frac{1}{x}, \quad Y=-x(xy+\alpha_1), \quad T=t
\end{gathered}
\end{gather}
from $\alpha_0,\alpha_1,\alpha_2,t,x,y$ to $\beta_0,\beta_1,\beta_2,T,X,Y$. Then the system \eqref{SG2} can also be written in the new variables $T,X,Y$ and parameters $\beta_0,\beta_1,\beta_2$ as a Hamiltonian system. This new system tends to the system $A_2^{(1)}${\rm : \rm}
\begin{equation}\label{SA2}
  \left\{
  \begin{aligned}
   \frac{dX}{dT} &=-X^2+4XY-2TX-2\beta_1,\\
   \frac{dY}{dT} &=-2Y^2+2XY+2TY+\beta_2
   \end{aligned}
  \right.
\end{equation}
with the Hamiltonian
\begin{align}\label{HA2}
\begin{split}
&H_{A_2^{(1)}}(X,Y,T;\beta_0,\beta_1,\beta_2)\\
&=-2TXY-X^2Y+2XY^2-2\beta_1Y-\beta_2X.
\end{split}
\end{align}
\end{Theorem}

By putting $q=X$, we have the Painlev\'e IV equation:
\begin{equation*}
\frac{d^2q}{dT^2}=\frac{1}{2q} \left(\frac{dq}{dT}\right)^2+\frac{3}{2}q^3+4Tq^2+2(T^2-a)q+\frac{b}{q},
\end{equation*}
where
\begin{equation*}
a=1-\beta_1-2\beta_2, \quad b=-2\beta_1^2.
\end{equation*}

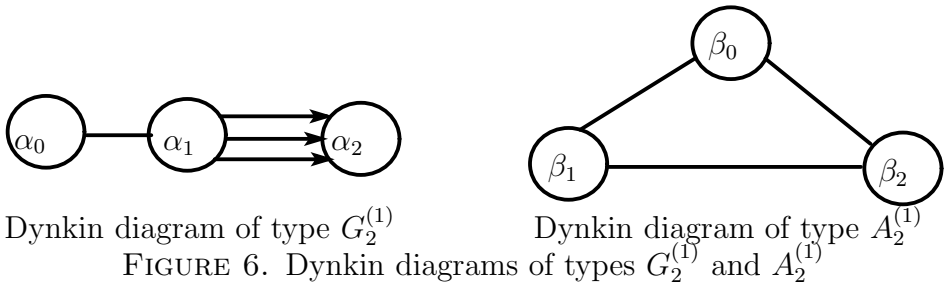
\begin{figure}[h]
\unitlength 0.1in
\begin{picture}(49.04,10.40)(11.70,-17.40)
%
\special{pn 20}%
\special{ar 3035 1389 202 189  0.0000000 6.2831853}%
%
\special{pn 20}%
\special{ar 1387 1352 201 188  0.0000000 6.2831853}%
%
\special{pn 20}%
\special{ar 2126 1368 202 188  0.0000000 6.2831853}%
%
\special{pn 20}%
\special{pa 1600 1370}%
\special{pa 1938 1370}%
\special{fp}%
%
\special{pn 20}%
\special{pa 2298 1271}%
\special{pa 2856 1271}%
\special{fp}%
\special{sh 1}%
\special{pa 2856 1271}%
\special{pa 2789 1251}%
\special{pa 2803 1271}%
\special{pa 2789 1291}%
\special{pa 2856 1271}%
\special{fp}%
%
\special{pn 20}%
\special{pa 2332 1389}%
\special{pa 2828 1389}%
\special{fp}%
\special{sh 1}%
\special{pa 2828 1389}%
\special{pa 2761 1369}%
\special{pa 2775 1389}%
\special{pa 2761 1409}%
\special{pa 2828 1389}%
\special{fp}%
%
\special{pn 20}%
\special{pa 2270 1496}%
\special{pa 2863 1496}%
\special{fp}%
\special{sh 1}%
\special{pa 2863 1496}%
\special{pa 2796 1476}%
\special{pa 2810 1496}%
\special{pa 2796 1516}%
\special{pa 2863 1496}%
\special{fp}%
\put(12.2000,-14.4000){\makebox(0,0)[lb]{$\alpha_0$}}%
\put(20.0000,-14.5000){\makebox(0,0)[lb]{$\alpha_1$}}%
\put(28.8000,-14.5000){\makebox(0,0)[lb]{$\alpha_2$}}%
%
\special{pn 20}%
\special{ar 4972 888 202 188  0.0000000 6.2831853}%
%
\special{pn 20}%
\special{ar 4122 1520 202 188  0.0000000 6.2831853}%
%
\special{pn 20}%
\special{ar 5872 1544 202 188  0.0000000 6.2831853}%
%
\special{pn 20}%
\special{pa 4782 969}%
\special{pa 4182 1342}%
\special{fp}%
%
\special{pn 20}%
\special{pa 4322 1536}%
\special{pa 5652 1536}%
\special{fp}%
%
\special{pn 20}%
\special{pa 5152 969}%
\special{pa 5702 1415}%
\special{fp}%
\put(48.5200,-9.4500){\makebox(0,0)[lb]{$\beta_0$}}%
\put(40.1200,-15.8500){\makebox(0,0)[lb]{$\beta_1$}}%
\put(57.4200,-16.0900){\makebox(0,0)[lb]{$\beta_2$}}%
\put(11.7000,-19.1000){\makebox(0,0)[lb]{Dynkin diagram of type $G_2^{(1)}$}}%
\put(39.4000,-19.1000){\makebox(0,0)[lb]{Dynkin diagram of type $A_2^{(1)}$}}%
\end{picture}%
\label{fig:CPVB512}
\caption{Dynkin diagrams of types $G_2^{(1)}$ and $A_2^{(1)}$}
\end{figure}

{\bf Proof of Theorem \ref{th:G2A2}.}
Notice that
\begin{equation*}
\alpha_0+2\alpha_1+3\alpha_2=\beta_0+\beta_1+\beta_2=1
\end{equation*}
and the change of variables from $(x,y,t)$ to $(X,Y,T)$ in Theorem \ref{th:G2A2} is symplectic. Choose $S_i$ as
\begin{equation*}
S_0:=s_2s_1s_2, \quad S_1:=s_0, \quad S_2:=s_1, \quad S_3:=s_2.
\end{equation*}
The transformations $S_0,S_1,S_2$ are reflections of
\begin{equation*}
\beta_0=\alpha_1+3\alpha_2,  \quad \beta_1=\alpha_0, \quad \beta_2=\alpha_1 \ {\rm respectively. \rm}
\end{equation*}
We can verify
\begin{align*}
        S_{0}: (X,Y,T;\beta_0,\beta_1,\beta_2) \rightarrow &(X+\frac{2\beta_0}{2Y-X-2T},Y+\frac{\beta_0}{2Y-X-2T},T;\\
        &-\beta_0,\beta_1+\beta_0,\beta_2+\beta_0), \\
        S_{1}: (X,Y,T;\beta_0,\beta_1,\beta_2) \rightarrow &(X,Y-\frac{\beta_1}{X},T;\beta_0+\beta_1,-\beta_1,\beta_2+\beta_1), \\
        S_{2}: (X,Y,T;\beta_0,\beta_1,\beta_2) \rightarrow &(X+\frac{\beta_2}{Y},Y,T;\beta_0+\beta_2,\beta_1+\beta_2,-\beta_2), \\
        S_{3}: (X,Y,T;\beta_0,\beta_1,\beta_2) \rightarrow &(-\sqrt{-1}X,-\frac{\sqrt{-1}(X-2Y+2T)}{2},-\sqrt{-1}T;\beta_2,\beta_1,\beta_0).
\end{align*}
The proof has thus been completed. \qed

\begin{Proposition}
The system \eqref{SG2} admits the following transformation $\varphi$ as its B{\"a}cklund transformation\rm{:\rm}
\begin{align*}
\begin{split}
\varphi: &(x,y,t;\alpha_0,\alpha_1,\alpha_2) \ra (\frac{1}{2x(xy+\alpha_1)},\\
&-2x(xy+\alpha_1)\{xy+2tx^2y+2x^4y^2+4\alpha_1 x^3y+2\alpha_1^2 x^2+2\alpha_1 tx+2\alpha_1+3\alpha_2\},t;\\
&\alpha_1,\alpha_1+3\alpha_2,\frac{\alpha_0-\alpha_1-3\alpha_2}{3}).
\end{split}
\end{align*}
\end{Proposition}
We note that this transformation $\varphi$ is pulled back the diagram automorphism $\pi$ of the system \eqref{SA2} by the symplectic transformation defined in Theorem \ref{th:G2A2}
\begin{equation*}
\pi:(X,Y,T;\beta_0,\beta_1,\beta_2) \ra (-2Y,\frac{2T+X-2Y}{2},T;\beta_1,\beta_2,\beta_0).
\end{equation*}

\section{The system of $D_3^{(2)}$}

In this section, we present a 2-parameter family of polynomial Hamiltonian systems given by
\begin{equation}\label{SD32}
  \left\{
  \begin{aligned}
   \frac{dx}{dt} &=\frac{2x^2y}{t}-x^2-\frac{2\alpha_0x}{t}+\frac{1}{t},\\
   \frac{dy}{dt} &=-\frac{2xy^2}{t}+2xy+\frac{2\alpha_0y}{t}+\alpha_1
   \end{aligned}
  \right. 
\end{equation}
with the polynomial Hamiltonian
\begin{align}\label{HD32}
H_{D_3^{(2)}}(x,y,t;\alpha_0,\alpha_1,\alpha_2)=\frac{x^2y(y-t)-(2\alpha_0y+t\alpha_1)x+y}{t}.
\end{align}
Here $x$ and $y$ denote unknown complex variables and $\alpha_0,\alpha_1,\alpha_2$ are complex parameters satisfying the relation:
\begin{equation}
\alpha_0+\alpha_1+\alpha_2=\frac{1}{2}.
\end{equation}

\begin{Theorem}\label{th:D32}
The system \eqref{SD32} admits extended affine Weyl group symmetry of type $D_3^{(2)}$ as the group of its B{\"a}cklund transformations {\rm (cf. \cite{N3}), \rm} whose generators are explicitly given as follows{\rm : \rm}
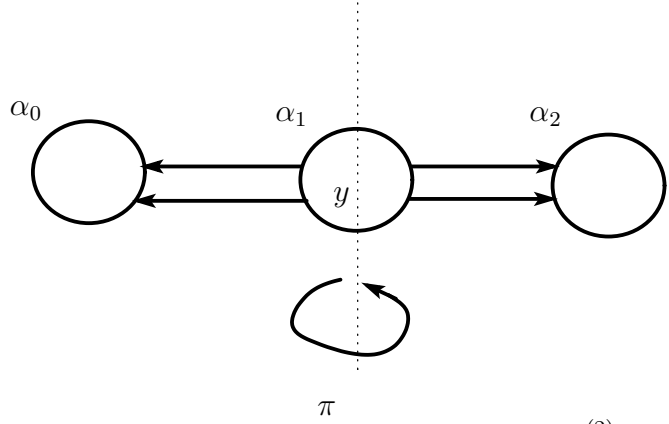
\begin{figure}[ht]
\unitlength 0.1in
\begin{picture}(34.26,19.90)(15.30,-24.20)
%
\special{pn 20}%
\special{ar 3343 1357 293 267  0.0000000 6.2831853}%
%
\special{pn 20}%
\special{ar 4663 1387 293 267  0.0000000 6.2831853}%
%
\special{pn 20}%
\special{pa 3633 1287}%
\special{pa 4373 1287}%
\special{fp}%
\special{sh 1}%
\special{pa 4373 1287}%
\special{pa 4306 1267}%
\special{pa 4320 1287}%
\special{pa 4306 1307}%
\special{pa 4373 1287}%
\special{fp}%
%
\special{pn 20}%
\special{pa 3623 1457}%
\special{pa 4363 1457}%
\special{fp}%
\special{sh 1}%
\special{pa 4363 1457}%
\special{pa 4296 1437}%
\special{pa 4310 1457}%
\special{pa 4296 1477}%
\special{pa 4363 1457}%
\special{fp}%
\put(32.2300,-14.6700){\makebox(0,0)[lb]{$y$}}%
\put(31.4000,-25.9000){\makebox(0,0)[lb]{$\pi$}}%
%
\special{pn 20}%
\special{ar 1943 1317 293 267  0.0000000 6.2831853}%
%
\special{pn 20}%
\special{pa 3053 1287}%
\special{pa 2243 1287}%
\special{fp}%
\special{sh 1}%
\special{pa 2243 1287}%
\special{pa 2310 1307}%
\special{pa 2296 1287}%
\special{pa 2310 1267}%
\special{pa 2243 1287}%
\special{fp}%
%
\special{pn 20}%
\special{pa 3083 1467}%
\special{pa 2203 1467}%
\special{fp}%
\special{sh 1}%
\special{pa 2203 1467}%
\special{pa 2270 1487}%
\special{pa 2256 1467}%
\special{pa 2270 1447}%
\special{pa 2203 1467}%
\special{fp}%
%
\special{pn 8}%
\special{pa 3350 430}%
\special{pa 3350 2350}%
\special{dt 0.045}%
\special{pa 3350 2350}%
\special{pa 3350 2349}%
\special{dt 0.045}%
%
\special{pn 20}%
\special{pa 3260 1880}%
\special{pa 3230 1887}%
\special{pa 3200 1896}%
\special{pa 3170 1907}%
\special{pa 3140 1921}%
\special{pa 3111 1939}%
\special{pa 3083 1962}%
\special{pa 3058 1990}%
\special{pa 3035 2019}%
\special{pa 3018 2051}%
\special{pa 3008 2081}%
\special{pa 3005 2110}%
\special{pa 3013 2136}%
\special{pa 3030 2158}%
\special{pa 3055 2176}%
\special{pa 3086 2192}%
\special{pa 3121 2206}%
\special{pa 3156 2220}%
\special{pa 3192 2233}%
\special{pa 3229 2244}%
\special{pa 3265 2255}%
\special{pa 3301 2263}%
\special{pa 3336 2269}%
\special{pa 3371 2273}%
\special{pa 3404 2274}%
\special{pa 3437 2272}%
\special{pa 3467 2266}%
\special{pa 3496 2256}%
\special{pa 3523 2242}%
\special{pa 3548 2224}%
\special{pa 3570 2200}%
\special{pa 3589 2172}%
\special{pa 3604 2141}%
\special{pa 3613 2108}%
\special{pa 3615 2076}%
\special{pa 3609 2046}%
\special{pa 3593 2021}%
\special{pa 3570 1999}%
\special{pa 3542 1980}%
\special{pa 3512 1961}%
\special{pa 3510 1960}%
\special{sp}%
%
\special{pn 20}%
\special{pa 3550 1980}%
\special{pa 3400 1910}%
\special{fp}%
\special{sh 1}%
\special{pa 3400 1910}%
\special{pa 3452 1956}%
\special{pa 3448 1933}%
\special{pa 3469 1920}%
\special{pa 3400 1910}%
\special{fp}%
\put(15.3000,-10.2000){\makebox(0,0)[lb]{$\alpha_0$}}%
\put(29.2000,-10.5000){\makebox(0,0)[lb]{$\alpha_1$}}%
\put(42.5000,-10.5000){\makebox(0,0)[lb]{$\alpha_2$}}%
\end{picture}%
\label{fig:CPIIID3.bak}
\caption{Dynkin diagram of type $D_3^{(2)}$}
\end{figure}
\begin{align*}
        s_{0}: (x,y,t;\alpha_0,\alpha_1,\alpha_2) &\rightarrow (-x,-y+\frac{2\alpha_0}{x}-\frac{1}{x^2},-t;-\alpha_0,\alpha_1+2\alpha_0,\alpha_2),\\
        s_{1}: (x,y,t;\alpha_0,\alpha_1,\alpha_2) &\rightarrow (x+\frac{\alpha_1}{y},y,t;\alpha_0+\alpha_1,-\alpha_1,\alpha_2+\alpha_1), \\
        s_{2}: (x,y,t;\alpha_0,\alpha_1,\alpha_2) &\rightarrow (x,y-t,-t;\alpha_0,\alpha_1+2\alpha_2,-\alpha_2), \\
        {\pi}: (x,y,t;\alpha_0,\alpha_1,\alpha_2) &\rightarrow (\frac{1}{tx},-tx(yx+\alpha_1),t;\alpha_2,\alpha_1,\alpha_0).
\end{align*}
\end{Theorem}

\begin{Theorem}\label{th:HoD32}
Let us consider a polynomial Hamiltonian system with Hamiltonian $H \in {\Bbb C}(t)[x,y]$. We assume that

$(A1)$ $deg(H)=5$ with respect to $x,y$.

$(A2)$ This system becomes again a polynomial Hamiltonian system in each coordinate $r_i \ (i=0,1,2)$ {\rm (cf. \cite{M}): \rm}
\begin{align*}
r_0&:x_0=x, \ y_0=y-\frac{2\alpha_0}{x}+\frac{1}{x^2},\\
r_1&:x_1=\frac{1}{x}, \ y_1=-(yx+\alpha_1)x, \\
r_2&:x_2=\frac{1}{x}, \ y_2=-((y-t)x+\alpha_1+2\alpha_2)x.
\end{align*}
Then such a system coincides with the system \eqref{SD32}.
\end{Theorem}

\begin{Theorem}\label{th:SB3SD32}
For the system \eqref{B3}, we make the change of parameters and variables\rm{:\rm}
\begin{gather}
\begin{gathered}\label{PB3D32}
x=\frac{\varepsilon TX}{1+\varepsilon TX}, \quad y=\frac{(1+\varepsilon TX)(\varepsilon TXY+Y+A_1\varepsilon T)}{\varepsilon T}, \quad t=\varepsilon T,
\end{gathered}\\
\begin{gathered}\label{VB3D32}
\alpha_0=-\frac{1-2A_2 \varepsilon}{\varepsilon}, \quad \alpha_1=\frac{1}{\varepsilon}, \quad \alpha_2=A_1, \quad \alpha_3=A_0
\end{gathered}
\end{gather}
from $x,y,t,\alpha_0,\alpha_1,\alpha_2,\alpha_3$ to $X,Y,T,A_0,A_1,A_2,\varepsilon$. Then the system \eqref{B3} can also be written in the new variables $T,X,Y$ and parameters $A_0,A_1,A_2,\varepsilon$ as a Hamiltonian system. This new system tends to the system \eqref{SD32} as $\varepsilon \ra 0$.
\end{Theorem}

By proving the following theorem, we see how the degeneration process given in Theorem \ref{th:SB3SD32} works on the B{\"a}cklund transformation group $W(B_3^{(1)})$ (cf. \cite{T2}).
\begin{Theorem}\label{th:WB3WD32}
For the degeneration process in Theorem \ref{th:SB3SD32}, we can choose a subgroup $W_{B_3^{(1)} \rightarrow D_3^{(2)}}$ of the B{\"a}cklund transformation group $W(B_3^{(1)})$ so that $W_{B_3^{(1)} \rightarrow D_3^{(2)}}$ converges to the B{\"a}cklund transformation group $W(D_3^{(2)})$ of the system \eqref{SD32} as $\varepsilon \rightarrow 0$.
\end{Theorem}

{\bf Proof of Therem \ref{th:WB3WD32}.}
Notice that
$$
2(A_0+A_1+A_2)=\alpha_0+\alpha_1+2\alpha_2+2\alpha_3=1
$$
and the change of variables from $(x,y)$ to $(X,Y)$ is symplectic. Choose $S_i \ (i=0,1,2)$ as
$$
S_0:=s_3, \ S_1:=s_2, \ S_2:=s_0s_1.
$$
Then the transformations $S_i$ are reflections of the parameters $A_0,A_1,A_2$. The transformation group $<S_0,S_1,S_2>$ coincides with the transformations given in Theorem \ref{th:D32}.
The proof has thus been completed. \qed

By the following theorem, we will show that the system \eqref{SD32} coincides with the system of type $C_2^{(1)}$ (see \cite{N3,Tsuda}).

\begin{Theorem}\label{th:D32C2}
For the system \eqref{SD32}, we make the change of parameters and variables\rm{:\rm}
\begin{gather}
\begin{gathered}\label{PD32C2}
\alpha_0=\frac{\beta_2-\beta_0}{2},  \quad \alpha_1=\beta_0, \quad \alpha_2=\beta_1,
\end{gathered}\\
\begin{gathered}\label{VD32C2}
x=\frac{1}{X}, \quad y=-X(YX+\beta_0), \quad T=t
\end{gathered}
\end{gather}
from $x,y,t,\alpha_0,\alpha_1,\alpha_2$ to $X,Y,T,\beta_0,\beta_1,\beta_2$. Then this new system coincides with the system of type $C_2^{(1)}$ \rm{(see \cite{N3}): \rm}
\begin{equation}\label{SC2}
  \left\{
  \begin{aligned}
   \frac{dX}{dT} &=\frac{2X^2Y}{T}-\frac{X^2}{T}+\frac{(\beta_0+\beta_2)X}{T}+1,\\
   \frac{dY}{dT} &=-\frac{2XY^2}{T}+\frac{2XY}{T}-\frac{(\beta_0+\beta_2)Y}{T}+\frac{\beta_0}{T}
   \end{aligned}
  \right. 
\end{equation}
with the Hamiltonian
\begin{align}
\begin{split}
H_{C_2^{(1)}}(X,Y,T;\beta_0,\beta_1,\beta_2)=&\frac{X^2Y(Y-1)+X[(\beta_0+\beta_2)Y-\beta_0]+TY}{T}\\
&(\beta_0+2\beta_1+\beta_2=1).
\end{split}
\end{align}
\end{Theorem}
By putting $q=\frac{x}{\tau}, \ T={\tau}^2$, we will see that this system \eqref{SC2} is equivalent to the third Painlev\'e equation:
\begin{equation*}
\frac{d^2q}{d{\tau}^2}=\frac{1}{q} \left(\frac{dq}{d{\tau}}\right)^2-\frac{1}{\tau} \frac{dq}{d{\tau}}+\frac{1}{\tau}(aq^2+b)+cq^3+\frac{d}{q},
\end{equation*}
where
\begin{equation*}
a=4(\beta_0-\beta_2), \quad b=-4(\beta_0+\beta_2-1), \quad c=4, \quad d=-4.
\end{equation*}

\begin{figure}[h]
\unitlength 0.1in
\begin{picture}(52.61,7.56)(9.30,-16.98)
%
\special{pn 20}%
\special{ar 2174 1388 201 217  0.0000000 6.2831853}%
%
\special{pn 20}%
\special{ar 3080 1412 201 218  0.0000000 6.2831853}%
%
\special{pn 20}%
\special{pa 2373 1330}%
\special{pa 2881 1330}%
\special{fp}%
\special{sh 1}%
\special{pa 2881 1330}%
\special{pa 2814 1310}%
\special{pa 2828 1330}%
\special{pa 2814 1350}%
\special{pa 2881 1330}%
\special{fp}%
%
\special{pn 20}%
\special{pa 2366 1470}%
\special{pa 2875 1470}%
\special{fp}%
\special{sh 1}%
\special{pa 2875 1470}%
\special{pa 2808 1450}%
\special{pa 2822 1470}%
\special{pa 2808 1490}%
\special{pa 2875 1470}%
\special{fp}%
\put(20.9200,-14.7800){\makebox(0,0)[lb]{$y$}}%
%
\special{pn 20}%
\special{ar 1214 1355 201 218  0.0000000 6.2831853}%
%
\special{pn 20}%
\special{pa 1975 1330}%
\special{pa 1419 1330}%
\special{fp}%
\special{sh 1}%
\special{pa 1419 1330}%
\special{pa 1486 1350}%
\special{pa 1472 1330}%
\special{pa 1486 1310}%
\special{pa 1419 1330}%
\special{fp}%
%
\special{pn 20}%
\special{pa 1996 1478}%
\special{pa 1392 1478}%
\special{fp}%
\special{sh 1}%
\special{pa 1392 1478}%
\special{pa 1459 1498}%
\special{pa 1445 1478}%
\special{pa 1459 1458}%
\special{pa 1392 1478}%
\special{fp}%
\put(9.3000,-11.1200){\makebox(0,0)[lb]{$\alpha_0$}}%
\put(18.8400,-11.3700){\makebox(0,0)[lb]{$\alpha_1$}}%
\put(27.9600,-11.3700){\makebox(0,0)[lb]{$\alpha_2$}}%
%
\special{pn 20}%
\special{ar 5084 1456 201 217  0.0000000 6.2831853}%
%
\special{pn 20}%
\special{ar 5990 1480 201 218  0.0000000 6.2831853}%
\put(40.4000,-15.3000){\makebox(0,0)[lb]{$Y$}}%
%
\special{pn 20}%
\special{ar 4124 1423 201 218  0.0000000 6.2831853}%
\put(38.4000,-11.8000){\makebox(0,0)[lb]{$\beta_0$}}%
\put(47.9400,-12.0500){\makebox(0,0)[lb]{$\beta_1$}}%
\put(57.0600,-12.0500){\makebox(0,0)[lb]{$\beta_2$}}%
\put(57.9000,-15.5000){\makebox(0,0)[lb]{$Y-1$}}%
%
\special{pn 20}%
\special{pa 4330 1360}%
\special{pa 4860 1360}%
\special{fp}%
\special{sh 1}%
\special{pa 4860 1360}%
\special{pa 4793 1340}%
\special{pa 4807 1360}%
\special{pa 4793 1380}%
\special{pa 4860 1360}%
\special{fp}%
%
\special{pn 20}%
\special{pa 4320 1510}%
\special{pa 4870 1510}%
\special{fp}%
\special{sh 1}%
\special{pa 4870 1510}%
\special{pa 4803 1490}%
\special{pa 4817 1510}%
\special{pa 4803 1530}%
\special{pa 4870 1510}%
\special{fp}%
%
\special{pn 20}%
\special{pa 5790 1370}%
\special{pa 5300 1370}%
\special{fp}%
\special{sh 1}%
\special{pa 5300 1370}%
\special{pa 5367 1390}%
\special{pa 5353 1370}%
\special{pa 5367 1350}%
\special{pa 5300 1370}%
\special{fp}%
%
\special{pn 20}%
\special{pa 5800 1540}%
\special{pa 5300 1540}%
\special{fp}%
\special{sh 1}%
\special{pa 5300 1540}%
\special{pa 5367 1560}%
\special{pa 5353 1540}%
\special{pa 5367 1520}%
\special{pa 5300 1540}%
\special{fp}%
\put(10.5000,-18.3000){\makebox(0,0)[lb]{Dynkin diagram of type $D_3^{(2)}$}}%
\put(39.7000,-18.3000){\makebox(0,0)[lb]{Dynkin diagram of type $C_2^{(1)}$}}%
\end{picture}%
\label{fig:CPIIID4}
\caption{Dynkin diagrams of types $D_3^{(2)}$ and $C_2^{(1)}$}
\end{figure}
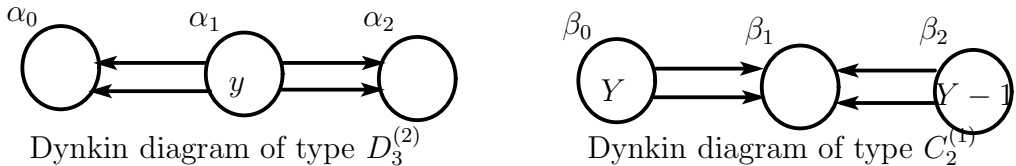

{\bf Proof of Theorem \ref{th:D32C2}.}
Notice that
\begin{equation*}
2(\alpha_0+\alpha_1+\alpha_2)=\beta_0+2\beta_1+\beta_2=1
\end{equation*}
and the change of variables from $(x,y,t)$ to $(X,Y,T)$ in Theorem \ref{th:D32C2} is symplectic. Choose $S_i, \ i=0,1,2,3$ as
\begin{equation*}
S_0:=s_1, \quad S_1:=s_2, \quad S_2:=s_0s_1s_0, \quad S_3:=s_0.
\end{equation*}
The transformations $S_0,S_1,S_2$ are reflections of
\begin{equation*}
\beta_0=\alpha_1,  \quad \beta_1=\alpha_2, \quad \beta_2=2\alpha_0+\alpha_1 \ {\rm respectively \rm}.
\end{equation*}
We can verify
\begin{align*}
        S_{0}: (X,Y,T;\beta_0,\beta_1,\beta_2) &\rightarrow (X+\frac{\beta_0}{Y},Y,T;-\beta_0,\beta_1+\beta_0,\beta_2),\\
        S_{1}: (X,Y,T;\beta_0,\beta_1,\beta_2) &\rightarrow (X,Y-\frac{2\beta_1}{X}+\frac{T}{X^2},-T;\beta_0+2\beta_1,-\beta_1,\beta_2+2\beta_1),\\
        S_{2}: (X,Y,T;\beta_0,\beta_1,\beta_2) &\rightarrow (X+\frac{\beta_2}{Y-1},Y,T;\beta_0,\beta_1+\beta_2,-\beta_2),\\
        S_{3}: (X,Y,T;\beta_0,\beta_1,\beta_2) &\rightarrow (-X,1-Y,-T;\beta_2,\beta_1,\beta_0).
\end{align*}
The proof has thus been completed. \qed

\section{The system of $A_2^{(2)}$}

In this section, we present a 1-parameter family of polynomial Hamiltonian systems given by
\begin{equation}\label{SA22}
  \left\{
  \begin{aligned}
   \frac{dx}{dt} &=x^4y+\alpha_0x^3+\frac{t}{2}x^2+1,\\
   \frac{dy}{dt} &=-2x^3y^2-3\alpha_0x^2y-txy-{\alpha_0}^2x-\frac{\alpha_0t}{2}
   \end{aligned}
  \right. 
\end{equation}
with the polynomial Hamiltonian
\begin{align}\label{HA22}
H_{A_2^{(2)}}(x,y,t;\alpha_0,\alpha_1)=\frac{x^4y^2}{2}+\alpha_0x^3y+\frac{1}{2}tx^2y+\frac{\alpha_0^2}{2}x^2+\frac{\alpha_0}{2}tx+y.
\end{align}
Here $x$ and $y$ denote unknown complex variables and $\alpha_0,\alpha_1$ are complex parameters satisfying the relation:
\begin{equation}
2\alpha_0+\alpha_1=1.
\end{equation}

By putting $q:=x$, we obtain the following equation:
\begin{equation*}
\frac{d^2q}{dt^2}=\left(\frac{2}{q}\right)\left(\frac{dq}{dt}\right)^2-\alpha_0 q^2-tq +\frac{q^2}{2}-\frac{2}{q}.
\end{equation*}

\begin{Theorem}\label{th:A22}
The system \eqref{SA22} admits affine Weyl group symmetry of type $A_2^{(2)}$ as the group of its B{\"a}cklund transformations {\rm (cf. \cite{N3}), \rm} whose generators are explicitly given as follows{\rm : \rm}
\begin{align*}
        s_{0}: &(x,y,t;\alpha_0,\alpha_1) \rightarrow (x+\frac{\alpha_0}{y},y,t;-\alpha_0,\alpha_1+4\alpha_0), \\
        s_{1}: &(x,y,t;\alpha_0,\alpha_1) \rightarrow (-x,-y+\frac{\alpha_1}{x}-\frac{t}{x^2}-\frac{2}{x^4},t;\alpha_0+\alpha_1,-\alpha_1).
\end{align*}
\end{Theorem}

\begin{figure}
\unitlength 0.1in
\begin{picture}(21.06,8.40)(23.40,-16.90)
%
\special{pn 20}%
\special{ar 4153 1357 293 267  0.0000000 6.2831853}%
\put(26.2000,-14.3000){\makebox(0,0)[lb]{$y$}}%
%
\special{pn 20}%
\special{ar 2753 1317 293 267  0.0000000 6.2831853}%
%
\special{pn 20}%
\special{pa 3863 1287}%
\special{pa 3053 1287}%
\special{fp}%
\special{sh 1}%
\special{pa 3053 1287}%
\special{pa 3120 1307}%
\special{pa 3106 1287}%
\special{pa 3120 1267}%
\special{pa 3053 1287}%
\special{fp}%
%
\special{pn 20}%
\special{pa 3893 1467}%
\special{pa 3013 1467}%
\special{fp}%
\special{sh 1}%
\special{pa 3013 1467}%
\special{pa 3080 1487}%
\special{pa 3066 1467}%
\special{pa 3080 1447}%
\special{pa 3013 1467}%
\special{fp}%
\put(23.4000,-10.2000){\makebox(0,0)[lb]{$\alpha_0$}}%
\put(37.3000,-10.5000){\makebox(0,0)[lb]{$\alpha_1$}}%
%
\special{pn 20}%
\special{pa 3830 1390}%
\special{pa 3020 1390}%
\special{fp}%
\special{sh 1}%
\special{pa 3020 1390}%
\special{pa 3087 1410}%
\special{pa 3073 1390}%
\special{pa 3087 1370}%
\special{pa 3020 1390}%
\special{fp}%
%
\special{pn 20}%
\special{pa 3890 1210}%
\special{pa 3010 1210}%
\special{fp}%
\special{sh 1}%
\special{pa 3010 1210}%
\special{pa 3077 1230}%
\special{pa 3063 1210}%
\special{pa 3077 1190}%
\special{pa 3010 1210}%
\special{fp}%
\put(24.3000,-18.6000){\makebox(0,0)[lb]{Dynkin diagram of type $A_2^{(2)}$}}%
\end{picture}%
\label{fig:A2}
\caption{Dynkin diagram of type $A_2^{(2)}$}
\end{figure}

\begin{Theorem}\label{th:HoA22}
Let us consider a polynomial Hamiltonian system with Hamiltonian $H \in {\Bbb C}(t)[x,y]$. We assume that

$(A1)$ $deg(H)=6$ with respect to $x,y$.

$(A2)$ This system becomes again a polynomial Hamiltonian system in each coordinate $r_i \ (i=0,1)${\rm (cf. \cite{M}): \rm}
\begin{align*}
r_0&:x_0=\frac{1}{x}, \ y_0=-(yx+\alpha_0)x,\\
r_1&:x_1=x, \ y_1=y-\frac{\alpha_1}{x}+\frac{t}{x^2}+\frac{2}{x^4}.
\end{align*}
Then such a system coincides with the system \eqref{SA22}.
\end{Theorem}

Theorems \ref{th:A22}, \ref{th:HoA22} can be cheched by a direct calculation, respectively.

\begin{Theorem}\label{th:SG2toSA22}
For the system \eqref{SG2}, we make the change of parameters and variables
\begin{gather}
\begin{gathered}\label{PG2A22}
\alpha_0=\frac{1}{4{\varepsilon}^6}, \quad \alpha_1=A_0, \quad \alpha_2=\frac{4A_1{\varepsilon}^6-1}{12{\varepsilon}^6}, \quad t=-\frac{1-{\varepsilon}^4 T}{\sqrt{2}{\varepsilon}^3},
\end{gathered}\\
\begin{gathered}\label{VG2A22}
x=\frac{\sqrt{2}{\varepsilon}^3 X}{2{\varepsilon}^2+X}, \quad y=\frac{(X+2{\varepsilon}^2)(XY+2{\varepsilon}^2Y+A_0)}{2\sqrt{2}{\varepsilon}^5}
\end{gathered}
\end{gather}
from $\alpha_0,\alpha_1,\alpha_2,t,x,y$ to $A_0,A_1,\varepsilon,T,X,Y$. Then the system \eqref{SG2} can also be written in the new variables $T,X,Y$ and parameters $A_0,A_1,\varepsilon$ as a Hamiltonian system. This new system tends to the system \eqref{SA22} as $\varepsilon \ra 0$.
\end{Theorem}

By proving the following theorem, we see how the degeneration process given in Theorem \ref{th:SG2toSA22} works on the B{\"a}cklund transformation group $W(G_2^{(1)})$ (cf. \cite{T2}).
\begin{Theorem}\label{th:WG2A22}
For the degeneration process in Theorem \ref{th:SG2toSA22}, we can choose a subgroup $W_{G_2^{(1)} \rightarrow A_2^{(2)}}$ of the B{\"a}cklund transformation group $W(G_2^{(1)})$ so that $W_{G_2^{(1)} \rightarrow A_2^{(2)}}$ converges to the B{\"a}cklund transformation group $W(A_2^{(2)})$ of the system \eqref{SA22} as $\varepsilon \rightarrow 0$.
\end{Theorem}

{\bf Proof of Theorem \ref{th:WG2A22}.}
Notice that
\begin{equation*}
2A_0+A_1=\alpha_0+2\alpha_1+3\alpha_2=1
\end{equation*}
and the change of variables from $(x,y)$ to $(X,Y)$ is symplectic. Since the change of parameters \eqref{PG2A22} is not one to one, we consider the degeneration process by introducing formal power series of the new parameters $A_0,A_1,\varepsilon$.

We choose $S_0,S_1$ as
\begin{equation*}
S_0:=s_1, \quad S_1:=s_0s_2
\end{equation*}
and put $W_{G_2^{(1)} \rightarrow A_2^{(2)}}=<S_0,S_1>$. Notice the $S_0,S_1$ are reflections of the parameters $A_0,A_1$, respectively.

Then we can obtain
\begin{align*}
S_0(A_0,A_1,\varepsilon)&=(-A_0,A_1+4A_0,\varepsilon(1+4A_0 \varepsilon^6)^{-1/6}),\\
S_1(A_0,A_1,\varepsilon)&=(A_0+A_1,-A_1,-\sqrt{-1} \varepsilon).
\end{align*}
Here, we make the same convention as in Section 2 that $(1+4A_0 \varepsilon^6)^{-1/6}$ means formal power series of $A_0 \varepsilon^6$ with 1 as constant term.

Then we can verify
\begin{align*}
S_0(X,Y,T)&=(X+\frac{A_0}{Y},Y,T),\\
S_1(X,Y,T)&=(-X,-Y+\frac{A_1}{X}-\frac{T}{X^2}-\frac{2}{X^4},T)
\end{align*}
as $\varepsilon \rightarrow 0$.

The proof has thus been completed. \qed

By the following theorem, we will show that the system \eqref{SA22} coincides with the system of type $A_1^{(1)}$.
\begin{Theorem}\label{A22toA1}
For the system \eqref{SA22}, we make the change of parameters and variables
\begin{gather}
\begin{gathered}\label{PA22A1}
\beta_0=\alpha_0+\alpha_1, \quad \beta_1=\alpha_0,
\end{gathered}\\
\begin{gathered}\label{VA22A1}
X=\frac{1}{x}, \quad Y=-x(xy+\alpha_0), \quad T=t
\end{gathered}
\end{gather}
from $x,y,t,\alpha_0,\alpha_1$ to $X,Y,t,\beta_0,\beta_1$. Then this new system coincides with the system of type $A_1^{(1)}$ \rm{(see \cite{N3}):\rm}
\begin{equation}\label{A1}
  \left\{
  \begin{aligned}
   \frac{dX}{dT} &=-X^2+Y-\frac{T}{2},\\
   \frac{dY}{dT} &=2XY+\beta_1
   \end{aligned}
  \right. 
\end{equation}
with the Hamiltonian
\begin{equation}
H_{A_1^{(1)}}(X,Y,T;\beta_0,\beta_1)=\frac{1}{2}Y^2-\left(X^2+\frac{T}{2}\right)Y-\beta_1X.
\end{equation}
\end{Theorem}
By putting $q=X$, we have the second Painlev\'e equation (see \cite{N3}):
\begin{equation*}
\frac{d^2q}{dT^2}=2q^3+Tq+\left(\beta_1-\frac{1}{2}\right).
\end{equation*}

{\bf Proof of Theorem \ref{A22toA1}.}
Notice that
\begin{equation}
2\alpha_0+\alpha_1=\beta_0+\beta_1=1
\end{equation}
and the change of variables from $(x,y,t)$ to $(X,Y,t)$ in Theorem \ref{A22toA1} is symplectic. Choose $S_1$ and $\pi$ as
\begin{equation*}
S_1:=s_0, \quad \pi:=s_1.
\end{equation*}
The transformations $S_1,\pi$ are reflections of
\begin{equation*}
\beta_0=\alpha_0+\alpha_1, \quad \beta_1=\alpha_0, \ \rm{respectively. \rm}
\end{equation*}
We can verify
\begin{align*}
        S_{1}: &(X,Y,t;\beta_0,\beta_1) \rightarrow (X+\frac{\beta_1}{Y},Y,t;\beta_0+2\beta_1,-\beta_1),\\
        \pi: &(X,Y,t;\beta_0,\beta_1) \rightarrow (-X,t+2X^2-Y,t;\beta_1,\beta_0).
\end{align*}
The proof has thus been completed. \qed

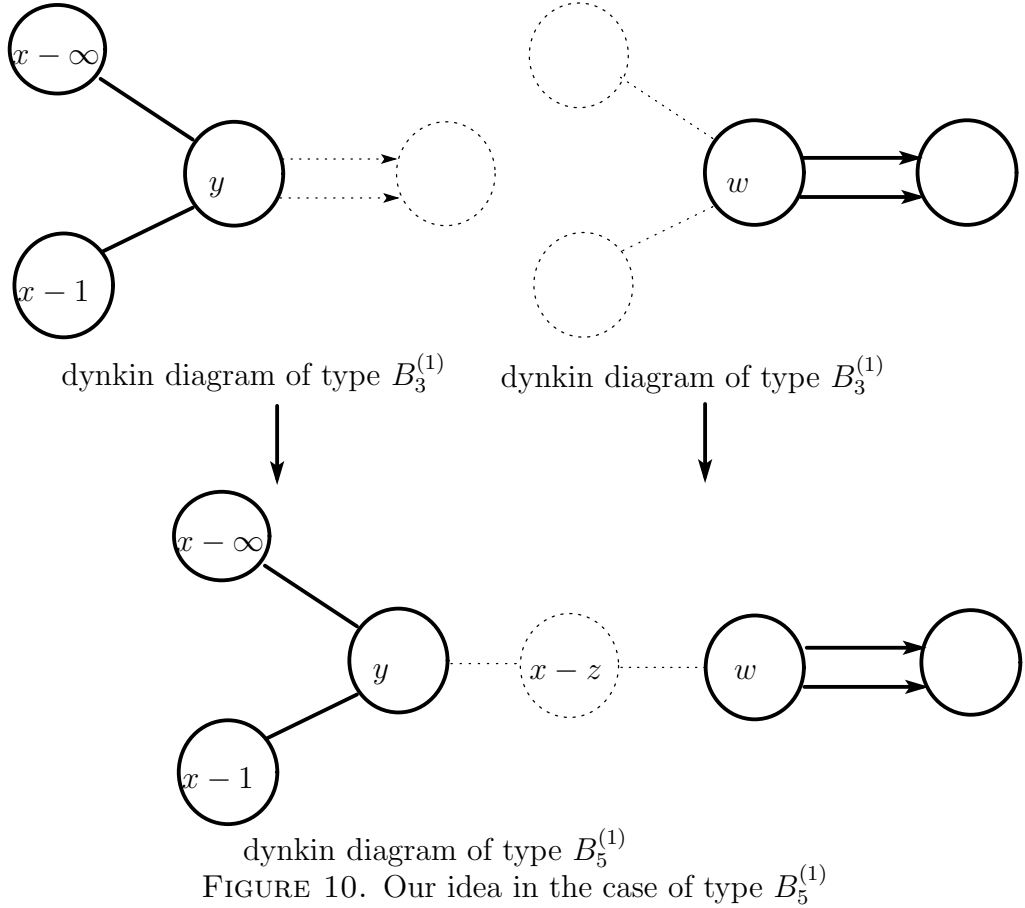
\begin{figure}[h]
\unitlength 0.1in
\begin{picture}(52.90,43.20)(20.40,-45.70)
%
\special{pn 20}%
\special{ar 2290 492 250 232  0.0000000 6.2831853}%
%
\special{pn 20}%
\special{ar 2323 1729 258 272  0.0000000 6.2831853}%
%
\special{pn 20}%
\special{ar 3215 1144 257 272  0.0000000 6.2831853}%
%
\special{pn 8}%
\special{ar 4323 1144 258 272  0.0000000 0.0452830}%
\special{ar 4323 1144 258 272  0.1811321 0.2264151}%
\special{ar 4323 1144 258 272  0.3622642 0.4075472}%
\special{ar 4323 1144 258 272  0.5433962 0.5886792}%
\special{ar 4323 1144 258 272  0.7245283 0.7698113}%
\special{ar 4323 1144 258 272  0.9056604 0.9509434}%
\special{ar 4323 1144 258 272  1.0867925 1.1320755}%
\special{ar 4323 1144 258 272  1.2679245 1.3132075}%
\special{ar 4323 1144 258 272  1.4490566 1.4943396}%
\special{ar 4323 1144 258 272  1.6301887 1.6754717}%
\special{ar 4323 1144 258 272  1.8113208 1.8566038}%
\special{ar 4323 1144 258 272  1.9924528 2.0377358}%
\special{ar 4323 1144 258 272  2.1735849 2.2188679}%
\special{ar 4323 1144 258 272  2.3547170 2.4000000}%
\special{ar 4323 1144 258 272  2.5358491 2.5811321}%
\special{ar 4323 1144 258 272  2.7169811 2.7622642}%
\special{ar 4323 1144 258 272  2.8981132 2.9433962}%
\special{ar 4323 1144 258 272  3.0792453 3.1245283}%
\special{ar 4323 1144 258 272  3.2603774 3.3056604}%
\special{ar 4323 1144 258 272  3.4415094 3.4867925}%
\special{ar 4323 1144 258 272  3.6226415 3.6679245}%
\special{ar 4323 1144 258 272  3.8037736 3.8490566}%
\special{ar 4323 1144 258 272  3.9849057 4.0301887}%
\special{ar 4323 1144 258 272  4.1660377 4.2113208}%
\special{ar 4323 1144 258 272  4.3471698 4.3924528}%
\special{ar 4323 1144 258 272  4.5283019 4.5735849}%
\special{ar 4323 1144 258 272  4.7094340 4.7547170}%
\special{ar 4323 1144 258 272  4.8905660 4.9358491}%
\special{ar 4323 1144 258 272  5.0716981 5.1169811}%
\special{ar 4323 1144 258 272  5.2528302 5.2981132}%
\special{ar 4323 1144 258 272  5.4339623 5.4792453}%
\special{ar 4323 1144 258 272  5.6150943 5.6603774}%
\special{ar 4323 1144 258 272  5.7962264 5.8415094}%
\special{ar 4323 1144 258 272  5.9773585 6.0226415}%
\special{ar 4323 1144 258 272  6.1584906 6.2037736}%
%
\special{pn 20}%
\special{pa 2516 647}%
\special{pa 2998 965}%
\special{fp}%
%
\special{pn 20}%
\special{pa 2528 1551}%
\special{pa 2998 1322}%
\special{fp}%
%
\special{pn 8}%
\special{pa 3468 1067}%
\special{pa 4058 1067}%
\special{dt 0.045}%
\special{sh 1}%
\special{pa 4058 1067}%
\special{pa 3991 1047}%
\special{pa 4005 1067}%
\special{pa 3991 1087}%
\special{pa 4058 1067}%
\special{fp}%
%
\special{pn 8}%
\special{pa 3456 1271}%
\special{pa 4070 1271}%
\special{dt 0.045}%
\special{sh 1}%
\special{pa 4070 1271}%
\special{pa 4003 1251}%
\special{pa 4017 1271}%
\special{pa 4003 1291}%
\special{pa 4070 1271}%
\special{fp}%
%
\special{pn 8}%
\special{ar 5019 524 259 274  0.0000000 0.0450281}%
\special{ar 5019 524 259 274  0.1801126 0.2251407}%
\special{ar 5019 524 259 274  0.3602251 0.4052533}%
\special{ar 5019 524 259 274  0.5403377 0.5853659}%
\special{ar 5019 524 259 274  0.7204503 0.7654784}%
\special{ar 5019 524 259 274  0.9005629 0.9455910}%
\special{ar 5019 524 259 274  1.0806754 1.1257036}%
\special{ar 5019 524 259 274  1.2607880 1.3058161}%
\special{ar 5019 524 259 274  1.4409006 1.4859287}%
\special{ar 5019 524 259 274  1.6210131 1.6660413}%
\special{ar 5019 524 259 274  1.8011257 1.8461538}%
\special{ar 5019 524 259 274  1.9812383 2.0262664}%
\special{ar 5019 524 259 274  2.1613508 2.2063790}%
\special{ar 5019 524 259 274  2.3414634 2.3864916}%
\special{ar 5019 524 259 274  2.5215760 2.5666041}%
\special{ar 5019 524 259 274  2.7016886 2.7467167}%
\special{ar 5019 524 259 274  2.8818011 2.9268293}%
\special{ar 5019 524 259 274  3.0619137 3.1069418}%
\special{ar 5019 524 259 274  3.2420263 3.2870544}%
\special{ar 5019 524 259 274  3.4221388 3.4671670}%
\special{ar 5019 524 259 274  3.6022514 3.6472795}%
\special{ar 5019 524 259 274  3.7823640 3.8273921}%
\special{ar 5019 524 259 274  3.9624765 4.0075047}%
\special{ar 5019 524 259 274  4.1425891 4.1876173}%
\special{ar 5019 524 259 274  4.3227017 4.3677298}%
\special{ar 5019 524 259 274  4.5028143 4.5478424}%
\special{ar 5019 524 259 274  4.6829268 4.7279550}%
\special{ar 5019 524 259 274  4.8630394 4.9080675}%
\special{ar 5019 524 259 274  5.0431520 5.0881801}%
\special{ar 5019 524 259 274  5.2232645 5.2682927}%
\special{ar 5019 524 259 274  5.4033771 5.4484053}%
\special{ar 5019 524 259 274  5.5834897 5.6285178}%
\special{ar 5019 524 259 274  5.7636023 5.8086304}%
\special{ar 5019 524 259 274  5.9437148 5.9887430}%
\special{ar 5019 524 259 274  6.1238274 6.1688555}%
%
\special{pn 8}%
\special{ar 5043 1726 259 274  0.0000000 0.0450281}%
\special{ar 5043 1726 259 274  0.1801126 0.2251407}%
\special{ar 5043 1726 259 274  0.3602251 0.4052533}%
\special{ar 5043 1726 259 274  0.5403377 0.5853659}%
\special{ar 5043 1726 259 274  0.7204503 0.7654784}%
\special{ar 5043 1726 259 274  0.9005629 0.9455910}%
\special{ar 5043 1726 259 274  1.0806754 1.1257036}%
\special{ar 5043 1726 259 274  1.2607880 1.3058161}%
\special{ar 5043 1726 259 274  1.4409006 1.4859287}%
\special{ar 5043 1726 259 274  1.6210131 1.6660413}%
\special{ar 5043 1726 259 274  1.8011257 1.8461538}%
\special{ar 5043 1726 259 274  1.9812383 2.0262664}%
\special{ar 5043 1726 259 274  2.1613508 2.2063790}%
\special{ar 5043 1726 259 274  2.3414634 2.3864916}%
\special{ar 5043 1726 259 274  2.5215760 2.5666041}%
\special{ar 5043 1726 259 274  2.7016886 2.7467167}%
\special{ar 5043 1726 259 274  2.8818011 2.9268293}%
\special{ar 5043 1726 259 274  3.0619137 3.1069418}%
\special{ar 5043 1726 259 274  3.2420263 3.2870544}%
\special{ar 5043 1726 259 274  3.4221388 3.4671670}%
\special{ar 5043 1726 259 274  3.6022514 3.6472795}%
\special{ar 5043 1726 259 274  3.7823640 3.8273921}%
\special{ar 5043 1726 259 274  3.9624765 4.0075047}%
\special{ar 5043 1726 259 274  4.1425891 4.1876173}%
\special{ar 5043 1726 259 274  4.3227017 4.3677298}%
\special{ar 5043 1726 259 274  4.5028143 4.5478424}%
\special{ar 5043 1726 259 274  4.6829268 4.7279550}%
\special{ar 5043 1726 259 274  4.8630394 4.9080675}%
\special{ar 5043 1726 259 274  5.0431520 5.0881801}%
\special{ar 5043 1726 259 274  5.2232645 5.2682927}%
\special{ar 5043 1726 259 274  5.4033771 5.4484053}%
\special{ar 5043 1726 259 274  5.5834897 5.6285178}%
\special{ar 5043 1726 259 274  5.7636023 5.8086304}%
\special{ar 5043 1726 259 274  5.9437148 5.9887430}%
\special{ar 5043 1726 259 274  6.1238274 6.1688555}%
%
\special{pn 20}%
\special{ar 5938 1138 259 274  0.0000000 6.2831853}%
%
\special{pn 20}%
\special{ar 7051 1138 259 274  0.0000000 6.2831853}%
%
\special{pn 8}%
\special{pa 5237 639}%
\special{pa 5720 959}%
\special{dt 0.045}%
\special{pa 5720 959}%
\special{pa 5719 959}%
\special{dt 0.045}%
%
\special{pn 8}%
\special{pa 5249 1547}%
\special{pa 5720 1317}%
\special{dt 0.045}%
\special{pa 5720 1317}%
\special{pa 5719 1317}%
\special{dt 0.045}%
%
\special{pn 20}%
\special{pa 6192 1061}%
\special{pa 6785 1061}%
\special{fp}%
\special{sh 1}%
\special{pa 6785 1061}%
\special{pa 6718 1041}%
\special{pa 6732 1061}%
\special{pa 6718 1081}%
\special{pa 6785 1061}%
\special{fp}%
%
\special{pn 20}%
\special{pa 6180 1266}%
\special{pa 6797 1266}%
\special{fp}%
\special{sh 1}%
\special{pa 6797 1266}%
\special{pa 6730 1246}%
\special{pa 6744 1266}%
\special{pa 6730 1286}%
\special{pa 6797 1266}%
\special{fp}%
\put(23.1000,-22.5000){\makebox(0,0)[lb]{dynkin diagram of type $B_3^{(1)}$}}%
\put(46.1000,-22.6000){\makebox(0,0)[lb]{dynkin diagram of type $B_3^{(1)}$}}%
%
\special{pn 20}%
\special{pa 3440 2360}%
\special{pa 3440 2740}%
\special{fp}%
\special{sh 1}%
\special{pa 3440 2740}%
\special{pa 3460 2673}%
\special{pa 3440 2687}%
\special{pa 3420 2673}%
\special{pa 3440 2740}%
\special{fp}%
%
\special{pn 20}%
\special{pa 5680 2350}%
\special{pa 5680 2730}%
\special{fp}%
\special{sh 1}%
\special{pa 5680 2730}%
\special{pa 5700 2663}%
\special{pa 5680 2677}%
\special{pa 5660 2663}%
\special{pa 5680 2730}%
\special{fp}%
\put(20.8100,-18.1400){\makebox(0,0)[lb]{$x-1$}}%
\put(20.5100,-5.7000){\makebox(0,0)[lb]{$x-\infty$}}%
\put(30.8100,-12.3400){\makebox(0,0)[lb]{$y$}}%
\put(57.9000,-12.4000){\makebox(0,0)[lb]{$w$}}%
%
\special{pn 20}%
\special{ar 3150 3042 250 232  0.0000000 6.2831853}%
%
\special{pn 20}%
\special{ar 3183 4279 258 272  0.0000000 6.2831853}%
%
\special{pn 20}%
\special{ar 4075 3694 257 272  0.0000000 6.2831853}%
%
\special{pn 20}%
\special{pa 3376 3197}%
\special{pa 3858 3515}%
\special{fp}%
%
\special{pn 20}%
\special{pa 3388 4101}%
\special{pa 3858 3872}%
\special{fp}%
\put(29.4100,-43.6400){\makebox(0,0)[lb]{$x-1$}}%
\put(29.1100,-31.2000){\makebox(0,0)[lb]{$x-\infty$}}%
\put(39.4100,-37.8400){\makebox(0,0)[lb]{$y$}}%
%
\special{pn 8}%
\special{ar 4970 3720 257 272  0.0000000 0.0453686}%
\special{ar 4970 3720 257 272  0.1814745 0.2268431}%
\special{ar 4970 3720 257 272  0.3629490 0.4083176}%
\special{ar 4970 3720 257 272  0.5444234 0.5897921}%
\special{ar 4970 3720 257 272  0.7258979 0.7712665}%
\special{ar 4970 3720 257 272  0.9073724 0.9527410}%
\special{ar 4970 3720 257 272  1.0888469 1.1342155}%
\special{ar 4970 3720 257 272  1.2703214 1.3156900}%
\special{ar 4970 3720 257 272  1.4517958 1.4971645}%
\special{ar 4970 3720 257 272  1.6332703 1.6786389}%
\special{ar 4970 3720 257 272  1.8147448 1.8601134}%
\special{ar 4970 3720 257 272  1.9962193 2.0415879}%
\special{ar 4970 3720 257 272  2.1776938 2.2230624}%
\special{ar 4970 3720 257 272  2.3591682 2.4045369}%
\special{ar 4970 3720 257 272  2.5406427 2.5860113}%
\special{ar 4970 3720 257 272  2.7221172 2.7674858}%
\special{ar 4970 3720 257 272  2.9035917 2.9489603}%
\special{ar 4970 3720 257 272  3.0850662 3.1304348}%
\special{ar 4970 3720 257 272  3.2665406 3.3119093}%
\special{ar 4970 3720 257 272  3.4480151 3.4933837}%
\special{ar 4970 3720 257 272  3.6294896 3.6748582}%
\special{ar 4970 3720 257 272  3.8109641 3.8563327}%
\special{ar 4970 3720 257 272  3.9924386 4.0378072}%
\special{ar 4970 3720 257 272  4.1739130 4.2192817}%
\special{ar 4970 3720 257 272  4.3553875 4.4007561}%
\special{ar 4970 3720 257 272  4.5368620 4.5822306}%
\special{ar 4970 3720 257 272  4.7183365 4.7637051}%
\special{ar 4970 3720 257 272  4.8998110 4.9451796}%
\special{ar 4970 3720 257 272  5.0812854 5.1266541}%
\special{ar 4970 3720 257 272  5.2627599 5.3081285}%
\special{ar 4970 3720 257 272  5.4442344 5.4896030}%
\special{ar 4970 3720 257 272  5.6257089 5.6710775}%
\special{ar 4970 3720 257 272  5.8071834 5.8525520}%
\special{ar 4970 3720 257 272  5.9886578 6.0340265}%
\special{ar 4970 3720 257 272  6.1701323 6.2155009}%
%
\special{pn 20}%
\special{ar 5940 3730 257 272  0.0000000 6.2831853}%
\put(58.3000,-37.9000){\makebox(0,0)[lb]{$w$}}%
%
\special{pn 20}%
\special{ar 7071 3704 259 274  0.0000000 6.2831853}%
%
\special{pn 20}%
\special{pa 6212 3627}%
\special{pa 6805 3627}%
\special{fp}%
\special{sh 1}%
\special{pa 6805 3627}%
\special{pa 6738 3607}%
\special{pa 6752 3627}%
\special{pa 6738 3647}%
\special{pa 6805 3627}%
\special{fp}%
%
\special{pn 20}%
\special{pa 6200 3832}%
\special{pa 6817 3832}%
\special{fp}%
\special{sh 1}%
\special{pa 6817 3832}%
\special{pa 6750 3812}%
\special{pa 6764 3832}%
\special{pa 6750 3852}%
\special{pa 6817 3832}%
\special{fp}%
%
\special{pn 8}%
\special{pa 4350 3710}%
\special{pa 4690 3710}%
\special{dt 0.045}%
\special{pa 4690 3710}%
\special{pa 4689 3710}%
\special{dt 0.045}%
%
\special{pn 8}%
\special{pa 5240 3730}%
\special{pa 5680 3730}%
\special{dt 0.045}%
\special{pa 5680 3730}%
\special{pa 5679 3730}%
\special{dt 0.045}%
\put(47.6000,-37.9000){\makebox(0,0)[lb]{$x-z$}}%
\put(32.6000,-47.4000){\makebox(0,0)[lb]{dynkin diagram of type $B_5^{(1)}$}}%
\end{picture}%
\label{fig:PS1}
\caption{Our idea in the case of type $B_5^{(1)}$}
\end{figure}

\section{The case of $B_5^{(1)}$}
In this section and next section, we present polynomial Hamiltonian systems in dimension four with affine Weyl group symmetry of types $B_5^{(1)}$ and $D_5^{(2)}$. Our idea is the following way:
\begin{enumerate}
\item We make a dynkin diagram given by connecting two dynkin diagrams of type $B_3^{(1)}$ (resp. $D_3^{(2)}$) by adding the term with invariant divisor $x-z$.

\item We make the symmetry associated with the dynkin diagram given by 1.

\item We make the holomorphy conditions $r_i$ associated with the symmetry given by 2.

\item We look for a polynomial Hamiltonian system with the holomorphy conditions given by 3, that is,
\begin{Problem}\label{Problem1}
Let us consider a polynomial Hamiltonian system with Hamiltonian $H \in {\Bbb C}(t)[x,y,z,w]$. We assume that

$(A1)$ $deg(H)=5$ with respect to $x,y,z,w$.

$(A2)$ This system becomes again a polynomial Hamiltonian system in each coordinate $r_i$ given in $4$.
\end{Problem}
\end{enumerate}
To solve Problem \ref{Problem1}, for the Hamiltonian satisfying the assumption $(A1)$ we only have to determine unknown coefficients so that they cancel the poles of Hamiltonian transformed by each $r_i$.

In this section, we present a 5-parameter family of polynomial Hamiltonian systems that can be considered as four-dimensional coupled Painlev\'e V systems, which is explicitly given by
\begin{equation}\label{SB5}
  \left\{
  \begin{aligned}
   \frac{dx}{dt} &=-\frac{2x^3y-2x^2y+(\alpha_1+2\alpha_2)x^2+(t-1+2\alpha_3+2\alpha_4+2\alpha_5)x-t}{t}\\
                 &-\frac{2(x-1)z(zw+\alpha_4)}{t},\\
   \frac{dy}{dt} &=\frac{3x^2y^2-2xy^2+2(\alpha_1+2\alpha_2)xy+(t-1+2\alpha_3+2\alpha_4+2\alpha_5)y+\alpha_2(\alpha_1+\alpha_2)}{t}\\
                 &+\frac{2yz(zw+\alpha_4)}{t},\\
   \frac{dz}{dt} &=-\frac{2z^3w-2z^2w+(\alpha_1+2\alpha_2+\alpha_3+2\alpha_4)z^2+(t-1+2\alpha_5)z-t}{t}\\
                 &-\frac{2(x-1)yz^2}{t},\\
   \frac{dw}{dt} &=\frac{3z^2w^2-2zw^2+2(\alpha_1+2\alpha_2+\alpha_3+2\alpha_4)zw+(t-1+2\alpha_5)w}{t}\\
                 &+\frac{\alpha_4(\alpha_1+2\alpha_2+\alpha_3+\alpha_4)+2(x-1)y(2zw+\alpha_4)}{t}
   \end{aligned}
  \right. 
\end{equation}
with the Hamiltonian
\begin{align}\label{HB5}
\begin{split}
&H_{{B_5}^{(1)}}(x,y,z,w,t;\alpha_0,\alpha_1, \dots ,\alpha_5)\\
&=H_{{B_3}^{(1)}}(x,y,t;\alpha_1,\alpha_2,\alpha_3+\alpha_4+\alpha_5)+H_{{B_3}^{(1)}}(z,w,t;\alpha_1+2\alpha_2+\alpha_3,\alpha_4,\alpha_5)\\
&-\frac{2(x-1)yz(zw+\alpha_4)}{t}.
\end{split}
\end{align}
Here $x,y,z$ and $w$ denote unknown complex variables, and $\alpha_0,\alpha_1,..,\alpha_5$ are complex parameters satisfying the relation:
\begin{equation}
\alpha_0+\alpha_1+2\alpha_2+2\alpha_3+2\alpha_4+2\alpha_5=1.
\end{equation}

\begin{Theorem}
The system \eqref{SB5} admits extended affine Weyl group symmetry of type $B_5^{(1)}$ as the group of its B{\"a}cklund transformations {\rm (cf. \cite{N3}), \rm} whose generators are explicitly given as follows{\rm : \rm}with the notation
$$
(*):=(x,y,z,w,t;\alpha_0,\alpha_1,\alpha_2,\alpha_3,\alpha_4,\alpha_5),
$$
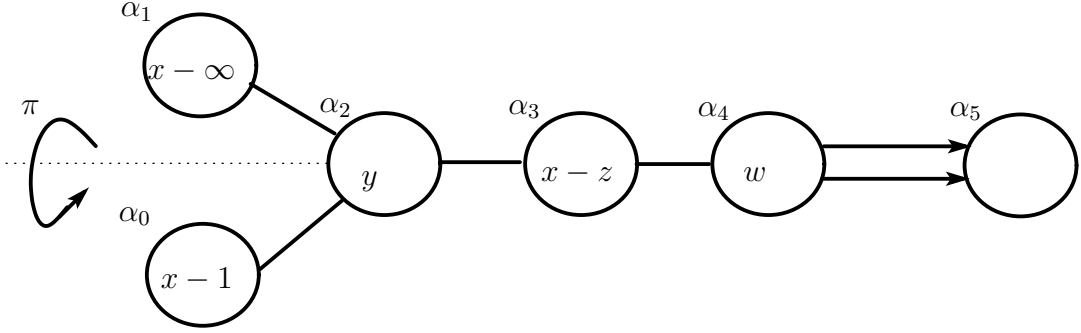
\begin{figure}
\unitlength 0.1in
\begin{picture}(56.03,17.91)(4.10,-23.01)
%
\special{pn 20}%
\special{ar 1427 937 293 267  0.0000000 6.2831853}%
%
\special{pn 20}%
\special{ar 1440 2034 293 267  0.0000000 6.2831853}%
%
\special{pn 20}%
\special{ar 2390 1454 293 267  0.0000000 6.2831853}%
%
\special{pn 20}%
\special{ar 5720 1461 293 267  0.0000000 6.2831853}%
%
\special{pn 20}%
\special{pa 1690 1034}%
\special{pa 2130 1294}%
\special{fp}%
%
\special{pn 20}%
\special{pa 1740 2004}%
\special{pa 2170 1644}%
\special{fp}%
%
\special{pn 20}%
\special{pa 4690 1361}%
\special{pa 5430 1361}%
\special{fp}%
\special{sh 1}%
\special{pa 5430 1361}%
\special{pa 5363 1341}%
\special{pa 5377 1361}%
\special{pa 5363 1381}%
\special{pa 5430 1361}%
\special{fp}%
%
\special{pn 20}%
\special{pa 4680 1531}%
\special{pa 5420 1531}%
\special{fp}%
\special{sh 1}%
\special{pa 5420 1531}%
\special{pa 5353 1511}%
\special{pa 5367 1531}%
\special{pa 5353 1551}%
\special{pa 5420 1531}%
\special{fp}%
\put(11.5000,-10.1000){\makebox(0,0)[lb]{$x-\infty$}}%
\put(12.2000,-21.0400){\makebox(0,0)[lb]{$x-1$}}%
\put(22.7000,-15.6400){\makebox(0,0)[lb]{$y$}}%
%
\special{pn 20}%
\special{ar 3420 1454 293 267  0.0000000 6.2831853}%
%
\special{pn 20}%
\special{ar 4400 1454 293 267  0.0000000 6.2831853}%
%
\special{pn 20}%
\special{pa 2680 1444}%
\special{pa 3100 1444}%
\special{fp}%
%
\special{pn 20}%
\special{pa 3710 1454}%
\special{pa 4090 1454}%
\special{fp}%
\put(32.1000,-15.4400){\makebox(0,0)[lb]{$x-z$}}%
\put(42.7000,-15.4400){\makebox(0,0)[lb]{$w$}}%
%
\special{pn 8}%
\special{pa 410 1450}%
\special{pa 2090 1450}%
\special{dt 0.045}%
\special{pa 2090 1450}%
\special{pa 2089 1450}%
\special{dt 0.045}%
%
\special{pn 20}%
\special{pa 880 1360}%
\special{pa 854 1332}%
\special{pa 828 1304}%
\special{pa 802 1279}%
\special{pa 776 1257}%
\special{pa 751 1239}%
\special{pa 727 1226}%
\special{pa 703 1220}%
\special{pa 680 1221}%
\special{pa 658 1229}%
\special{pa 638 1244}%
\special{pa 618 1264}%
\special{pa 601 1289}%
\special{pa 585 1319}%
\special{pa 571 1353}%
\special{pa 560 1390}%
\special{pa 551 1430}%
\special{pa 544 1472}%
\special{pa 540 1516}%
\special{pa 539 1560}%
\special{pa 542 1605}%
\special{pa 547 1649}%
\special{pa 556 1689}%
\special{pa 567 1725}%
\special{pa 581 1753}%
\special{pa 598 1772}%
\special{pa 617 1780}%
\special{pa 638 1775}%
\special{pa 660 1759}%
\special{pa 684 1735}%
\special{pa 709 1706}%
\special{pa 730 1680}%
\special{sp}%
%
\special{pn 20}%
\special{pa 690 1730}%
\special{pa 820 1600}%
\special{fp}%
\special{sh 1}%
\special{pa 820 1600}%
\special{pa 759 1633}%
\special{pa 782 1638}%
\special{pa 787 1661}%
\special{pa 820 1600}%
\special{fp}%
\put(4.9000,-11.9000){\makebox(0,0)[lb]{$\pi$}}%
\put(10.0000,-17.5000){\makebox(0,0)[lb]{$\alpha_0$}}%
\put(10.1000,-6.8000){\makebox(0,0)[lb]{$\alpha_1$}}%
\put(20.5000,-11.8000){\makebox(0,0)[lb]{$\alpha_2$}}%
\put(30.4000,-11.9000){\makebox(0,0)[lb]{$\alpha_3$}}%
\put(40.3000,-12.0000){\makebox(0,0)[lb]{$\alpha_4$}}%
\put(53.5000,-12.0000){\makebox(0,0)[lb]{$\alpha_5$}}%
\end{picture}%
\label{fig:CPVB57}
\caption{Dynkin diagram of type $B_5^{(1)}$}
\end{figure}
\begin{align*}
        s_{0}: (*) &\rightarrow (x,y-\frac{\alpha_0}{x-1},z,w,t;-\alpha_0,\alpha_1,\alpha_2+\alpha_0,\alpha_3,\alpha_4,\alpha_5), \\
        s_{1}: (*) &\rightarrow (x,y,z,w,t;\alpha_0,-\alpha_1,\alpha_2+\alpha_1,\alpha_3,\alpha_4,\alpha_5), \\
        s_{2}: (*) &\rightarrow (x+\frac{\alpha_2}{y},y,z,w,t;\alpha_0+\alpha_2,\alpha_1+\alpha_2,-\alpha_2,\alpha_3+\alpha_2,\alpha_4,\alpha_5), \\
        s_{3}: (*) &\rightarrow (x,y-\frac{\alpha_3}{x-z},z,w+\frac{\alpha_3}{x-z},t;\alpha_0,\alpha_1,\alpha_2+\alpha_3,-\alpha_3,\alpha_4+\alpha_3,\alpha_5), \\
        s_{4}: (*) &\rightarrow (x,y,z+\frac{\alpha_4}{w},w,t;\alpha_0,\alpha_1,\alpha_2,\alpha_3+\alpha_4,-\alpha_4,\alpha_5+\alpha_4), \\
        s_{5}: (*) &\rightarrow (x,y,z,w-\frac{2\alpha_5}{z}+\frac{t}{z^2},-t;\alpha_0,\alpha_1,\alpha_2,\alpha_3,\alpha_4+2\alpha_5,-\alpha_5), \\
        {\pi}: (*) &\rightarrow  (\frac{x}{x-1},-(x-1)((x-1)y+\alpha_2),\frac{z}{z-1},-(z-1)((z-1)w+\alpha_4),-t;\\
&\alpha_1,\alpha_0,\alpha_2,\alpha_3,\alpha_4,\alpha_5).
\end{align*}
\end{Theorem}

\begin{Theorem}\label{th:HoB5}
Let us consider a polynomial Hamiltonian system with Hamiltonian $H \in {\Bbb C}(t)[x,y,z,w]$. We assume that

$(A1)$ $deg(H)=5$ with respect to $x,y,z,w$.

$(A2)$ This system becomes again a polynomial Hamiltonian system in each coordinate $r_i \ (i=0,1, \dots ,5)${\rm : \rm}
\begin{align*}
r_0&:x_0=-((x-1)y-\alpha_0)y, \quad y_0=\frac{1}{y}, \quad z_0=z, \quad w_0=w,\\
r_1&:x_1=\frac{1}{x}, \quad y_1=-(yx+\alpha_1+\alpha_2)x, \quad z_1=z, \quad w_1=w, \\
r_2&:x_2=\frac{1}{x}, \quad y_2=-(yx+\alpha_2)x, \quad z_2=z, \quad w_2=w, \\
r_3&:x_3=-((x-z)y-\alpha_3)y, \quad y_3=\frac{1}{y}, \quad z_3=z, \quad w_3=w+y, \\
r_4&:x_4=x, \quad y_4=y, \quad z_4=\frac{1}{z}, \quad w_4=-(wz+\alpha_4)z, \\
r_5&:x_5=x, \quad y_5=y, \quad z_5=z, \quad w_5=w-\frac{2\alpha_5}{z}+\frac{t}{z^2}.
\end{align*}
Then such a system coincides with the system \eqref{SB5}.
\end{Theorem}

By the following theorem, we will show that the system \eqref{SB5} coincides with a 5-parameter family of four-dimensional coupled Painlev\'e $V$ systems with the affine Weyl group symmetry of type $D_5^{(1)}$ (see \cite{Sasa6}).
\begin{Theorem}\label{th:B5D5}
For the system \eqref{SB5}, we make the change of parameters and variables
\begin{gather}
\begin{gathered}\label{PB5D5}
\beta_0=\alpha_4+2\alpha_5,  \quad \beta_1=\alpha_4, \quad \beta_2=\alpha_3, \quad \beta_3=\alpha_2, \quad \beta_4=\alpha_0, \quad \beta_5=\alpha_1,
\end{gathered}\\
\begin{gathered}\label{VB5D5}
X=\frac{1}{z}, \quad Y=-z(zw+\alpha_4), \quad Z=\frac{1}{x}, \quad {W}=-x(xy+\alpha_2),  \quad T=-t
\end{gathered}
\end{gather}
from $\alpha_0,\alpha_1,\alpha_2,\alpha_3,\alpha_4,\alpha_5,t,x,y,z,w$ to $\beta_0,\beta_1,\beta_2,\beta_3,\beta_4,\beta_5,T,X,Y,Z,W$. Then the system \eqref{SB5} can also be written in the new variables $T,X,Y,Z,W$ and parameters $\beta_0,\beta_1,\beta_2,\beta_3,\beta_4,\beta_5$ as a Hamiltonian system. This new system coincides with the system of type $D_5^{(1)}${\rm : \rm}
\begin{equation}\label{SD5}
  \left\{
  \begin{aligned}
   \frac{dX}{dT} &=\frac{2X^2Y}{T}+X^2-\frac{2XY}{T}-(1+\frac{2\beta_2+2\beta_3+\beta_5+\beta_4}{T})X+\frac{\beta_2+\beta_5}{T}\\
                 &+\frac{2Z((Z-1)W+\beta_3)}{T},\\
   \frac{dY}{dT} &=-\frac{2XY^2}{T}+\frac{Y^2}{T}-2XY+(1+\frac{2\beta_2+2\beta_3+\beta_5+\beta_4}{T})Y-\beta_1,\\
   \frac{dZ}{dT} &=\frac{2Z^2W}{T}+Z^2-\frac{2ZW}{T}-(1+\frac{\beta_5+\beta_4}{T})Z+\frac{\beta_5}{T}+\frac{2YZ(Z-1)}{T},\\
   \frac{dW}{dT} &=-\frac{2ZW^2}{T}+\frac{W^2}{T}-2ZW+(1+\frac{\beta_5+\beta_4}{T})W-\beta_3-\frac{2Y(-W+2ZW+\beta_3)}{T}
   \end{aligned}
  \right. 
\end{equation}
with the Hamiltonian
\begin{align}\label{HD5}
\begin{split}
H &=H_{A_3^{(1)}}(X,Y,T;\beta_2+\beta_5,\beta_1,\beta_2+2\beta_3+\beta_4)+H_{A_3^{(1)}}(Z,W,T;\beta_5,\beta_3,\beta_4)\\
&+\frac{2YZ\{(Z-1)W+\beta_3\}}{T}.
\end{split}
\end{align}
\end{Theorem}

{\bf Proof of Theorem \ref{th:B5D5}.}
Notice that
\begin{align*}
\begin{split}
&\alpha_0+\alpha_1+2\alpha_2+2\alpha_3+2\alpha_4+2\alpha_5\\
&=\beta_0+\beta_1+2\beta_2+2\beta_3+\beta_4+\beta_5=1
\end{split}
\end{align*}
and the change of variables from $(x,y,z,w,t)$ to $(X,Y,Z,W,T)$ in Theorem \ref{th:B5D5} is symplectic. Choose $S_i, \ i=0,1, \dots ,8$ as
\begin{align*}
\begin{split}
&S_0:=s_5s_4s_5, \quad S_1:=s_4, \quad S_2:=s_3, \quad S_3:=s_2, \quad S_4:=s_0,\\
&S_5:=s_1, \quad S_6:=s_5, \quad S_7:=\pi, \quad S_8:=\pi{s_5}.
\end{split}
\end{align*}
The transformations $S_0,S_1, \dots ,S_5$ are reflections of
\begin{align*}
\begin{split}
\beta_0=\alpha_4+2\alpha_5, \quad \beta_1=\alpha_4, \quad \beta_2=\alpha_3, \quad \beta_3=\alpha_2, \quad \beta_4=\alpha_0, \quad \beta_5=\alpha_1 \ {\rm respectively. \rm}
\end{split}
\end{align*}
By using the notation
$$
(*):=(X,Y,Z,W,T;\beta_0,\beta_1,\beta_2,\beta_3,\beta_4,\beta_5),
$$
we can verify
\begin{figure}
\unitlength 0.1in
\begin{picture}(33.91,19.42)(20.40,-25.10)
%
\special{pn 20}%
\special{ar 2399 1049 251 283  0.0000000 6.2831853}%
%
\special{pn 20}%
\special{ar 3218 1598 249 283  0.0000000 6.2831853}%
%
\special{pn 20}%
\special{pa 2623 1151}%
\special{pa 2997 1427}%
\special{fp}%
%
\special{pn 20}%
\special{pa 2665 2180}%
\special{pa 3030 1799}%
\special{fp}%
%
\special{pn 20}%
\special{ar 4093 1598 250 283  0.0000000 6.2831853}%
%
\special{pn 20}%
\special{pa 3465 1587}%
\special{pa 3821 1587}%
\special{fp}%
%
\special{pn 20}%
\special{ar 5115 1033 316 271  0.0000000 6.2831853}%
%
\special{pn 20}%
\special{ar 5168 2225 250 284  0.0000000 6.2831853}%
%
\special{pn 20}%
\special{pa 4287 1462}%
\special{pa 4873 1177}%
\special{fp}%
%
\special{pn 20}%
\special{pa 4308 1721}%
\special{pa 4934 2134}%
\special{fp}%
\put(21.3500,-7.3800){\makebox(0,0)[lb]{$\beta_1$}}%
\put(21.3500,-19.0200){\makebox(0,0)[lb]{$\beta_0$}}%
\put(30.1600,-12.9400){\makebox(0,0)[lb]{$\beta_2$}}%
\put(39.1500,-13.0700){\makebox(0,0)[lb]{$\beta_3$}}%
\put(48.9300,-7.3800){\makebox(0,0)[lb]{$\beta_4$}}%
\put(49.3400,-19.1500){\makebox(0,0)[lb]{$\beta_5$}}%
\put(21.1200,-23.2900){\makebox(0,0)[lb]{$Y+T$}}%
\put(22.6000,-11.6000){\makebox(0,0)[lb]{$Y$}}%
\put(29.6400,-17.2600){\makebox(0,0)[lb]{$X-Z$}}%
\put(39.1500,-17.3400){\makebox(0,0)[lb]{$W$}}%
\put(48.5800,-11.3600){\makebox(0,0)[lb]{$Z-1$}}%
\put(50.2100,-23.4200){\makebox(0,0)[lb]{$Z$}}%
%
\special{pn 20}%
\special{ar 2355 2239 315 271  0.0000000 6.2831853}%
\end{picture}%
\label{fig:CPVB513}
\caption{Dynkin diagram of type $B_5^{(1)}$}
\end{figure}
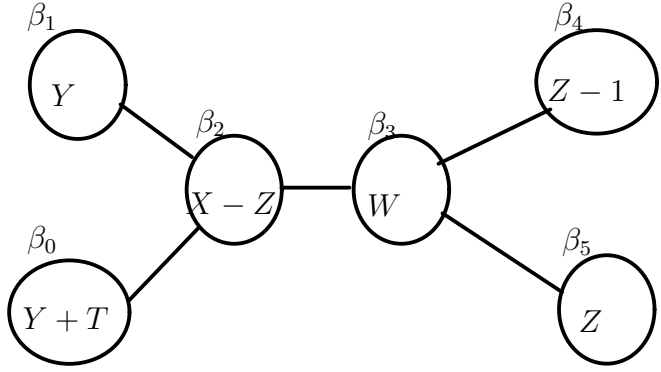
\begin{align*}
        S_0: (*) &\rightarrow (X+\frac{\beta_0}{Y+T},Y,Z,W,T;-\beta_0,\beta_1,\beta_2+\beta_0,\beta_3,\beta_4,\beta_5), \\
        S_1: (*) &\rightarrow (X+\frac{\beta_1}{Y},Y,Z,W,T;\beta_0,-\beta_1,\beta_2+\beta_1,\beta_3,\beta_4,\beta_5), \\
        S_2: (*) &\rightarrow (X,Y-\frac{\beta_2}{X-Z},Z,W+\frac{\beta_2}{X-Z},T;\beta_0+\beta_2,\beta_1+\beta_2,-\beta_2,\beta_3+\beta_2,\beta_4,\beta_5), \\
        S_3: (*) &\rightarrow (X,Y,Z+\frac{\beta_3}{W},W,T;\beta_0,\beta_1,\beta_2+\beta_3,-\beta_3,\beta_4+\beta_3,\beta_5+\beta_3), \\
        S_4: (*) &\rightarrow (X,Y,Z,W-\frac{\beta_4}{Z-1},T;\beta_0,\beta_1,\beta_2,\beta_3+\beta_4,-\beta_4,\beta_5), \\
        S_5: (*) &\rightarrow (X,Y,Z,W-\frac{\beta_5}{Z},T;\beta_0,\beta_1,\beta_2,\beta_3+\beta_5,\beta_4,-\beta_5), \\
        S_6: (*) &\rightarrow (X,Y+T,Z,W,-T;\beta_1,\beta_0,\beta_2,\beta_3,\beta_4,\beta_5), \\
        S_7: (*) &\rightarrow (1-X,-Y,1-Z,-W,-T;\beta_0,\beta_1,\beta_2,\beta_3,\beta_5,\beta_4), \\
        S_8: (*) &\rightarrow (1-X,-Y-T,1-Z,-W,T;\beta_1,\beta_0,\beta_2,\beta_3,\beta_5,\beta_4).
\end{align*}
The transformations $S_i, \ i=0,1,\dots ,8,$ define a represention of the affine Weyl group of type $D_5^{(1)}$, that is, they satisfy the following relations\rm{:\rm}
\begin{align*}
&S_0^2=S_1^2=\dots =S_5^2=S_6^2=S_7^2=S_8^2=1,\\
&(S_0S_1)^2=(S_0S_3)^2=(S_0S_4)^2=(S_0S_5)^2=(S_1S_3)^2\\
&=(S_1S_4)^2=(S_1S_5)^2=(S_2S_4)^2=(S_2S_5)^2=(S_4S_5)^2=1,\\
&(S_0S_2)^3=(S_1S_2)^3=(S_2S_3)^3=(S_3S_4)^3=(S_3S_5)^3=1,\\
&S_6(S_0,S_1,S_2,S_3,S_4,S_5)=(S_1,S_0,S_2,S_3,S_4,S_5)S_6,\\
&S_7(S_0,S_1,S_2,S_3,S_4,S_5)=(S_0,S_1,S_2,S_3,S_5,S_4)S_7,\\
&S_8(S_0,S_1,S_2,S_3,S_4,S_5)=(S_1,S_0,S_2,S_3,S_5,S_4)S_8.
\end{align*}
The proof has thus been completed. \qed

\begin{Proposition}
The system \eqref{SB5} admits the following transformation as its B{\"a}-\\
cklund transformation\rm{:\rm}
\begin{align*}
\begin{split}
&\varphi: (x,y,z,w,t;\alpha_0,\alpha_1,\alpha_2,\alpha_3,\alpha_4,\alpha_5) \ra\\
&(\frac{t}{(t+z^2w+\alpha_4z)},\frac{(t+z^2w+\alpha_4z)(-tx+tz-xz^2w+z^3w-(\alpha_3+\alpha_4)xz+\alpha_4z^2)}{txz},\\
&\frac{t}{(t+x^2y+z^2w+\alpha_2x+\alpha_4z)},\\
&\frac{(t+x^2y+z^2w+\alpha_2x+\alpha_4z)\{(x-t)(t+x^2y+z^2w+\alpha_4z)-(\alpha_0+\alpha_2)x+\alpha_2x^2\}}{tx},\\
&-t;\alpha_4,\alpha_4+2\alpha_5,\alpha_3,\alpha_2,\alpha_0,\frac{\alpha_1-\alpha_0}{2}).
\end{split}
\end{align*}
\end{Proposition}
We note that this transformation $\varphi$ is pulled back the following diagram automorphism $\pi$ of the system \eqref{SD5} by the symplectic transformation defined in Theorem \ref{th:B5D5}
\begin{align*}
\begin{split}
\pi:&(X,Y,Z,W,T;\beta_0,\beta_1,\beta_2,\beta_3,\beta_4,\beta_5) \ra\\
&((Y+W+T)/T,-T(Z-1),(Y+T)/T,-T(X-Z),-T;\beta_5,\beta_4,\beta_3,\beta_2,\beta_1,\beta_0).
\end{split}
\end{align*}

\section{The system of $D_5^{(2)}$}

In this section, we present a 4-parameter family of polynomial Hamiltonian systems that can be considered as four-dimensional coupled Painlev\'e III systems, which is given as follows
\begin{equation}\label{SD52}
  \left\{
  \begin{aligned}
   \frac{dx}{dt} &=\frac{2x^2y}{t}-x^2-\frac{2(\alpha_2+\alpha_3+\alpha_4)x}{t}+\frac{1}{t}+\frac{2z(zw+\alpha_3)}{t},\\
   \frac{dy}{dt} &=-\frac{2xy^2}{t}+2xy+\frac{2(\alpha_2+\alpha_3+\alpha_4)y}{t}+\alpha_1,\\
   \frac{dz}{dt} &=\frac{2z^2w}{t}-z^2-\frac{2\alpha_4z}{t}+\frac{1}{t}+\frac{2yz^2}{t},\\
   \frac{dw}{dt} &=-\frac{2zw^2}{t}+2zw+\frac{2\alpha_4w}{t}+\alpha_3-\frac{2y(2zw+\alpha_3)}{t}
   \end{aligned}
  \right. 
\end{equation}
with the Hamiltonian
\begin{align}
\begin{split}
&H_{D_5^{(2)}}(x,y,z,w,t;\alpha_0,\alpha_1, \dots ,\alpha_4)\\
&=H_{D_3^{(2)}}(x,y,t;\alpha_2+\alpha_3+\alpha_4,\alpha_1,-\frac{1}{2}+\alpha_0)+H_{D_3^{(2)}}(z,w,t;\alpha_4,\alpha_3,\frac{1}{2}-\alpha_3-\alpha_4)\\
&+\frac{2yz(zw+\alpha_3)}{t}.
\end{split}
\end{align}
Here $x,y,z$ and $w$ denote unknown complex variables, and $\alpha_0,\alpha_1, \dots,\alpha_4$ are complex parameters satisfying the relation:
\begin{equation}
\alpha_0+\alpha_1+\alpha_2+\alpha_3+\alpha_4=\frac{1}{2}.
\end{equation}

\begin{Theorem}\label{th:D52}
The system \eqref{SD52} admits extended affine Weyl group symmetry of type $D_5^{(2)}$ as the group of its B{\"a}cklund transformations {\rm (cf. \cite{N3}), \rm} whose generators are explicitly given as follows{\rm : \rm}with the notation
$$
(*):=(x,y,z,w,t;\alpha_0,\alpha_1,\alpha_2,\alpha_3,\alpha_4),
$$
\begin{align*}
s_{0}: (*) &\rightarrow (x,y-t,z,w,-t;-\alpha_0,\alpha_1+2\alpha_0,\alpha_2,\alpha_3,\alpha_4), \\
s_{1}: (*) &\rightarrow (x+\frac{\alpha_1}{y},y,z,w,t;\alpha_0+\alpha_1,-\alpha_1,\alpha_2+\alpha_1,\alpha_3,\alpha_4), \\
s_{2}: (*) &\rightarrow  (x,y-\frac{\alpha_2}{x-z},z,w+\frac{\alpha_2}{x-z},t;\alpha_0,\alpha_1+\alpha_2,-\alpha_2,\alpha_3+\alpha_2,\alpha_4), \\
s_{3}: (*) &\rightarrow (x,y,z+\frac{\alpha_3}{w},w,t;\alpha_0,\alpha_1,\alpha_2+\alpha_3,-\alpha_3,\alpha_4+\alpha_3), \\
s_{4}: (*) &\rightarrow (-x,-y,-z,-w+\frac{2\alpha_4}{z}-\frac{1}{z^2},-t;\alpha_0,\alpha_1,\alpha_2,\alpha_3+2\alpha_4,-\alpha_4), \\
{\pi}: (*) &\rightarrow (\frac{1}{tz},-t(zw+\alpha_3)z,\frac{1}{tx},-t(xy+\alpha_1)x,t;\alpha_4,\alpha_3,\alpha_2,\alpha_1,\alpha_0).
\end{align*}
\end{Theorem}

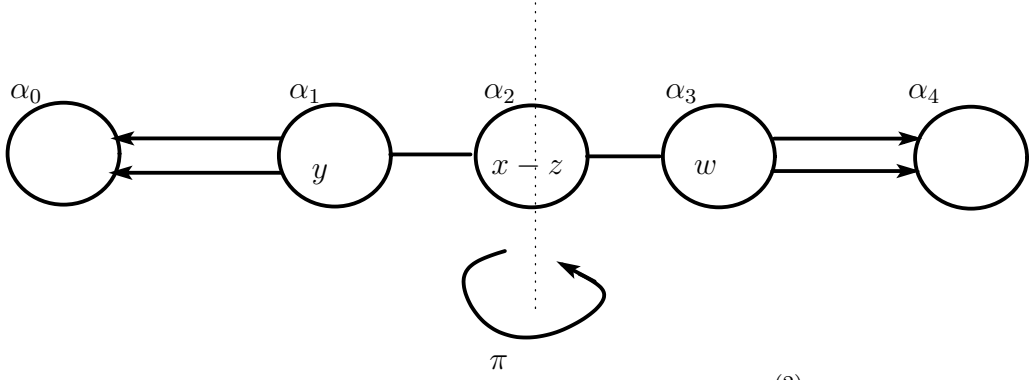
\begin{figure}
\unitlength 0.1in
\begin{picture}(53.36,17.54)(6.77,-24.04)
%
\special{pn 20}%
\special{ar 2390 1450 293 267  0.0000000 6.2831853}%
%
\special{pn 20}%
\special{ar 5720 1461 293 267  0.0000000 6.2831853}%
%
\special{pn 20}%
\special{pa 4690 1361}%
\special{pa 5430 1361}%
\special{fp}%
\special{sh 1}%
\special{pa 5430 1361}%
\special{pa 5363 1341}%
\special{pa 5377 1361}%
\special{pa 5363 1381}%
\special{pa 5430 1361}%
\special{fp}%
%
\special{pn 20}%
\special{pa 4680 1531}%
\special{pa 5420 1531}%
\special{fp}%
\special{sh 1}%
\special{pa 5420 1531}%
\special{pa 5353 1511}%
\special{pa 5367 1531}%
\special{pa 5353 1551}%
\special{pa 5420 1531}%
\special{fp}%
\put(22.7000,-15.6400){\makebox(0,0)[lb]{$y$}}%
%
\special{pn 20}%
\special{ar 3420 1454 293 267  0.0000000 6.2831853}%
%
\special{pn 20}%
\special{ar 4400 1454 293 267  0.0000000 6.2831853}%
%
\special{pn 20}%
\special{pa 2680 1444}%
\special{pa 3100 1444}%
\special{fp}%
%
\special{pn 20}%
\special{pa 3710 1454}%
\special{pa 4090 1454}%
\special{fp}%
\put(32.1000,-15.4400){\makebox(0,0)[lb]{$x-z$}}%
\put(42.7000,-15.4400){\makebox(0,0)[lb]{$w$}}%
%
\special{pn 20}%
\special{ar 970 1440 293 267  0.0000000 6.2831853}%
%
\special{pn 20}%
\special{pa 2110 1360}%
\special{pa 1250 1360}%
\special{fp}%
\special{sh 1}%
\special{pa 1250 1360}%
\special{pa 1317 1380}%
\special{pa 1303 1360}%
\special{pa 1317 1340}%
\special{pa 1250 1360}%
\special{fp}%
%
\special{pn 20}%
\special{pa 2110 1540}%
\special{pa 1240 1540}%
\special{fp}%
\special{sh 1}%
\special{pa 1240 1540}%
\special{pa 1307 1560}%
\special{pa 1293 1540}%
\special{pa 1307 1520}%
\special{pa 1240 1540}%
\special{fp}%
%
\special{pn 8}%
\special{pa 3440 650}%
\special{pa 3440 2250}%
\special{dt 0.045}%
\special{pa 3440 2250}%
\special{pa 3440 2249}%
\special{dt 0.045}%
%
\special{pn 20}%
\special{pa 3280 1950}%
\special{pa 3243 1961}%
\special{pa 3208 1973}%
\special{pa 3174 1986}%
\special{pa 3143 2000}%
\special{pa 3116 2016}%
\special{pa 3093 2035}%
\special{pa 3077 2056}%
\special{pa 3066 2082}%
\special{pa 3063 2110}%
\special{pa 3066 2141}%
\special{pa 3074 2173}%
\special{pa 3087 2206}%
\special{pa 3103 2238}%
\special{pa 3123 2268}%
\special{pa 3145 2296}%
\special{pa 3169 2321}%
\special{pa 3195 2342}%
\special{pa 3222 2360}%
\special{pa 3251 2375}%
\special{pa 3281 2386}%
\special{pa 3312 2395}%
\special{pa 3344 2401}%
\special{pa 3377 2403}%
\special{pa 3410 2403}%
\special{pa 3445 2400}%
\special{pa 3480 2394}%
\special{pa 3515 2386}%
\special{pa 3551 2375}%
\special{pa 3587 2361}%
\special{pa 3623 2344}%
\special{pa 3659 2325}%
\special{pa 3694 2304}%
\special{pa 3726 2281}%
\special{pa 3755 2257}%
\special{pa 3778 2232}%
\special{pa 3794 2208}%
\special{pa 3801 2184}%
\special{pa 3797 2161}%
\special{pa 3782 2140}%
\special{pa 3759 2121}%
\special{pa 3730 2102}%
\special{pa 3696 2083}%
\special{pa 3670 2070}%
\special{sp}%
%
\special{pn 20}%
\special{pa 3750 2110}%
\special{pa 3580 2020}%
\special{fp}%
\special{sh 1}%
\special{pa 3580 2020}%
\special{pa 3630 2069}%
\special{pa 3627 2045}%
\special{pa 3648 2034}%
\special{pa 3580 2020}%
\special{fp}%
\put(32.0000,-25.7000){\makebox(0,0)[lb]{$\pi$}}%
\put(6.9000,-11.5000){\makebox(0,0)[lb]{$\alpha_0$}}%
\put(21.5000,-11.5000){\makebox(0,0)[lb]{$\alpha_1$}}%
\put(31.7000,-11.5000){\makebox(0,0)[lb]{$\alpha_2$}}%
\put(41.2000,-11.5000){\makebox(0,0)[lb]{$\alpha_3$}}%
\put(53.9000,-11.5000){\makebox(0,0)[lb]{$\alpha_4$}}%
\end{picture}%
\label{fig:CPIIID32}
\caption{Dynkin diagram of type $D_5^{(2)}$}
\end{figure}

\begin{Theorem}\label{th:HoD52}
Let us consider a polynomial Hamiltonian system with Hamiltonian $H \in {\Bbb C}(t)[x,y,z,w]$. We assume that

$(A1)$ $deg(H)=5$ with respect to $x,y,z,w$.

$(A2)$ This system becomes again a polynomial Hamiltonian system in each coordinate $r_i \ (i=0,1, \dots ,4)${\rm : \rm}
\begin{align*}
r_0&:x_0=\frac{1}{x}, \quad y_0=-((y-t)x+2\alpha_0+\alpha_1)x, \quad z_0=z, \quad w_0=w,\\
r_1&:x_1=\frac{1}{x}, \quad y_1=-(yx+\alpha_1)x, \quad z_1=z, \quad w_1=w, \\
r_2&:x_2=-((x-z)y-\alpha_2)y, \quad y_2=\frac{1}{y}, \quad z_2=z, \quad w_2=w+y, \\
r_3&:x_3=x, \quad y_3=y, \quad z_3=\frac{1}{z}, \quad w_3=-(wz+\alpha_3)z, \\
r_4&:x_4=x, \quad y_4=y, \quad z_4=z, \quad w_4=w-\frac{2\alpha_4}{z}+\frac{1}{z^2}.
\end{align*}
Then such a system coincides with the system \eqref{SD52}.
\end{Theorem}

\begin{Theorem}\label{th:SB5SD52}
For the system \eqref{SB5}, we make the change of parameters and variables\rm{:\rm}
\begin{align}\label{PB5D52}
\begin{split}
x=\frac{\varepsilon TX}{1+\varepsilon TX}, \quad y=\frac{(1+\varepsilon TX)(\varepsilon TXY+Y+A_1\varepsilon T)}{\varepsilon T}, \quad z=\frac{\varepsilon TZ}{1+\varepsilon TZ},\\
w=\frac{(1+\varepsilon TZ)(\varepsilon TZW+W+A_3\varepsilon T)}{\varepsilon T}, \quad t=\varepsilon T,
\end{split}
\end{align}\\
\begin{equation}\label{VB5D52}
\alpha_0=-\frac{1}{\varepsilon}+2A_0, \quad \alpha_1=\frac{1}{\varepsilon}, \quad \alpha_2=A_1,\quad \alpha_3=A_2, \quad \alpha_4=A_3, \quad \alpha_5=A_4
\end{equation}
from $x,y,z,w,t,\alpha_0,\alpha_1,\dots,\alpha_5$ to $X,Y,Z,W,T,A_0,A_1,\dots,A_4,\varepsilon$. Then the system \eqref{SB5} can also be written in the new variables $T,X,Y,Z,W$ and parameters $A_0,A_1,\\
\dots,A_4,\varepsilon$ as a Hamiltonian system. This new system tends to the system \eqref{SD52} as $\varepsilon \ra 0$.
\end{Theorem}

By proving the following theorem, we see how the degeneration process given in Theorem \ref{th:SB5SD52} works on the B{\"a}cklund transformation group $W(B_5^{(1)})$ (cf. \cite{T2}).
\begin{Theorem}\label{th:WB5WD52}
For the degeneration process in Theorem \ref{th:SB5SD52}, we can choose a subgroup $W_{B_5^{(1)} \rightarrow D_5^{(2)}}$ of the B{\"a}cklund transformation group $W(B_5^{(1)})$ so that $W_{B_5^{(1)} \rightarrow D_5^{(2)}}$ converges to the B{\"a}cklund transformation group $W(D_5^{(2)})$ of the system \eqref{SD52} as $\varepsilon \rightarrow 0$.
\end{Theorem}

{\bf Proof of Theorem \ref{th:WB5WD52}.}
Note that
\begin{equation*}
2(A_0+A_1+A_2+A_3+A_4)
=\alpha_0+\alpha_1+2\alpha_2+2\alpha_3+2\alpha_4+2\alpha_5=1
\end{equation*}
and the change of variables from $(x,y,z,w)$ to $(X,Y,Z,W)$ is symplectic. Choose $S_i \ (i=0,1,\dots ,4)$ as
\begin{equation*}
S_0:=s_0s_1, \quad S_1:=s_2, \quad S_2:=s_3, \quad S_3=s_4, \quad S_4=s_5.
\end{equation*}
Then the transformations $S_i$ are reflections of the parameters $A_0,A_1,\dots,A_4$. The transformation group $<S_0,S_1,\dots .S_4>$ coincides with the transformations given in Theorem \ref{th:D52}.
The proof has thus been completed. \qed

By the following theorem, we will show that the system \eqref{SD52} coincides with a 4-parameter family of four-dimensional coupled Painlev\'e III systems with the affine Weyl group symmetry of type $B_4^{(1)}$ {\rm (see \cite{Sasa6}) \rm}.
\begin{Theorem}\label{th:D52B4}
For the system \eqref{SD52}, we first make the change of parameters and variables
\begin{gather}
\begin{gathered}\label{PD52B4}
\beta_0=\alpha_3+2\alpha_4, \quad \beta_1=\alpha_3, \quad \beta_2=\alpha_2, \quad \beta_3=\alpha_1, \quad \beta_4=\alpha_0,
\end{gathered}\\
\begin{gathered}\label{VD52B4}
X=\frac{1}{z}, \quad Y=-z(zw+\alpha_3), \quad Z=\frac{1}{x}, \quad W=-x(xy+\alpha_1), \quad T=t
\end{gathered}
\end{gather}
from $\alpha_0,\alpha_1, \dots,\alpha_4,t,x,y,z,w$ to $\beta_0,\beta_1, \dots ,\beta_4,T,X,Y,Z,W$. Then the system \eqref{SD52} can also be written in the new variables $T,X,Y,Z,W$ and parameters $\beta_0,\beta_1, \dots,\beta_4$ as a Hamiltonian system. This new system coincides with the system of type $B_4^{(1)}$\rm{:\rm}
\begin{equation}\label{SB4}
  \left\{
  \begin{aligned}
   \frac{dX}{dT} &=\frac{2X^2Y-{X^2}+(\beta_0+\beta_1)X+2\beta_3Z+2Z^2W}{T}+1,\\
   \frac{dY}{dT} &=\frac{-2XY^2+2XY-(\beta_0+\beta_1)Y+{\beta_1}}{T},\\
   \frac{dZ}{dT} &=\frac{2Z^2W-{Z^2}+(1-2\beta_4)Z+2YZ^2}{T}+1,\\
   \frac{dW}{dT} &=\frac{-2ZW^2+2{ZW}-(1-2\beta_4)W-2\beta_3Y-4YZW+\beta_3}{T}
   \end{aligned}
  \right. 
\end{equation}
with the Hamiltonian
\begin{align}
\begin{split}
H &=H_{C_2^{(1)}}(X,Y,T;\beta_0,\frac{1-\beta_0-\beta_1}{2},\beta_1)+H_{C_2^{(1)}}(Z,W,T;\beta_0+\beta_1+2\beta_2+\beta_3,\beta_4,\beta_3)\\
&+\frac{2YZ(ZW+\beta_3)}{T}\\
&=\frac{X^2Y(Y-1)+X\{(\beta_0+\beta_1)Y-{\beta_1}\}+TY}{T}\\
&+\frac{Z^2W(W-1)+Z\{(1-2\beta_4)W-{\beta_3}\}+TW}{T}+\frac{2YZ(ZW+\beta_3)}{T}.
\end{split}
\end{align}
\end{Theorem}

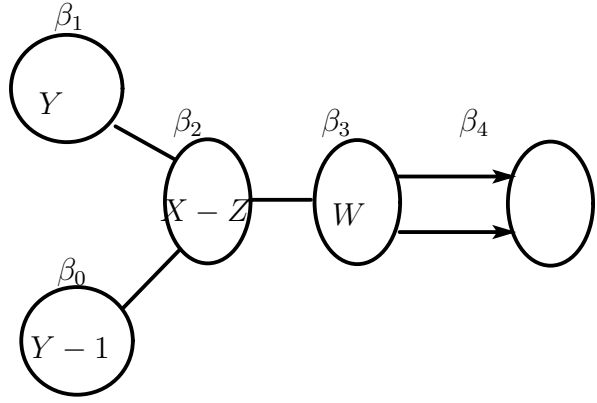
\begin{figure}
\unitlength 0.1in
\begin{picture}(30.47,21.05)(21.60,-27.90)
%
\special{pn 20}%
\special{ar 3190 1777 224 325  0.0000000 6.2831853}%
%
\special{pn 20}%
\special{ar 4984 1789 223 325  0.0000000 6.2831853}%
%
\special{pn 20}%
\special{ar 3973 1781 224 325  0.0000000 6.2831853}%
%
\special{pn 20}%
\special{pa 3411 1769}%
\special{pa 3730 1769}%
\special{fp}%
\put(30.0800,-14.1000){\makebox(0,0)[lb]{$\beta_2$}}%
\put(37.8400,-14.1000){\makebox(0,0)[lb]{$\beta_3$}}%
\put(45.0500,-14.1000){\makebox(0,0)[lb]{$\beta_4$}}%
%
\special{pn 20}%
\special{pa 4180 1647}%
\special{pa 4770 1647}%
\special{fp}%
\special{sh 1}%
\special{pa 4770 1647}%
\special{pa 4703 1627}%
\special{pa 4717 1647}%
\special{pa 4703 1667}%
\special{pa 4770 1647}%
\special{fp}%
%
\special{pn 20}%
\special{pa 4194 1938}%
\special{pa 4770 1938}%
\special{fp}%
\special{sh 1}%
\special{pa 4770 1938}%
\special{pa 4703 1918}%
\special{pa 4717 1938}%
\special{pa 4703 1958}%
\special{pa 4770 1938}%
\special{fp}%
%
\special{pn 20}%
\special{pa 2707 1384}%
\special{pa 3015 1556}%
\special{fp}%
%
\special{pn 20}%
\special{pa 2747 2336}%
\special{pa 3042 2032}%
\special{fp}%
\put(23.8500,-8.5500){\makebox(0,0)[lb]{$\beta_1$}}%
\put(24.0000,-21.9000){\makebox(0,0)[lb]{$\beta_0$}}%
\put(22.6000,-26.0000){\makebox(0,0)[lb]{$Y-1$}}%
\put(23.0500,-13.1800){\makebox(0,0)[lb]{$Y$}}%
\put(29.4000,-18.8000){\makebox(0,0)[lb]{$X-Z$}}%
\put(38.4000,-19.0000){\makebox(0,0)[lb]{$W$}}%
%
\special{pn 20}%
\special{ar 2506 2508 292 282  0.0000000 6.2831853}%
%
\special{pn 20}%
\special{ar 2452 1186 292 281  0.0000000 6.2831853}%
\end{picture}%
\label{fig:CPIIID34}
\caption{Dynkin diagram of type $B_4^{(1)}$}
\end{figure}

{\bf Proof of Theorem \ref{th:D52B4}.}
Notice that
\begin{equation*}
2(\alpha_0+\alpha_1+\alpha_2+\alpha_3+\alpha_4)=\beta_0+\beta_1+2\beta_2+2\beta_3+2\beta_4=1
\end{equation*}
and the change of variables from $(x,y,z,w,t)$ to $(X,Y,Z,W,T)$ in Theorem \ref{th:D52B4} is symplectic. Choose $S_i, \ i=0,1, \dots ,6$ as
\begin{align*}
\begin{split}
&S_0:=s_4s_3s_4, \quad S_1:=s_3, \quad S_2:=s_2,\\
&S_3:=s_1, \quad S_4:=s_0, \quad S_5:=s_4, \quad S_6:=\pi.
\end{split}
\end{align*}
The transformations $S_0,S_1, \dots ,S_4$ are reflections of
\begin{align*}
\beta_0=\alpha_3+2\alpha_4,  \quad \beta_1=\alpha_3, \quad \beta_2=\alpha_2, \quad \beta_3=\alpha_1, \quad \beta_4=\alpha_0 \ {\rm respectively. \rm}
\end{align*}
By using the notation
$$
(*):=(X,Y,Z,W,T;\beta_0,\beta_1,\beta_2,\beta_3,\beta_4),
$$
we can verify
\begin{align*}
        S_0: (*) &\rightarrow (X+\frac{\beta_0}{Y-1},Y,Z,W,T;-\beta_0,\beta_1,\beta_2+\beta_0,\beta_3,\beta_4),\\
        S_1: (*) &\rightarrow (X+\frac{\beta_1}{Y},Y,Z,W,T;\beta_0,-\beta_1,\beta_2+\beta_1,\beta_3,\beta_4),\\
        S_2: (*) &\rightarrow (X,Y-\frac{\beta_2}{X-Z},Z,W+\frac{\beta_2}{X-Z},T;\beta_0+\beta_2,\beta_1+\beta_2,-\beta_2,\beta_3+\beta_2,\beta_4),\\
        S_3: (*) &\rightarrow (X,Y,Z+\frac{\beta_3}{W},W,T;\beta_0,\beta_1,\beta_2+\beta_3,-\beta_3,\beta_4+\beta_3),\\
        S_4: (*) &\rightarrow (X,Y,Z,W-\frac{2\beta_4}{Z}+\frac{T}{Z^2},-T;\beta_0,\beta_1,\beta_2,\beta_3+2\beta_4,-\beta_4),\\
        S_5: (*) &\rightarrow (-X,1-Y,-Z,-W,-T;\beta_1,\beta_0,\beta_2,\beta_3,\beta_4),\\
        S_6: (*) &\rightarrow (\frac{T}{Z},-\frac{(ZW+\beta_3)Z}{T},\frac{T}{X},-\frac{(XY+\beta_1)X}{T},T;\beta_3+2\beta_4,\beta_3,\beta_2,\beta_1,\frac{\beta_0-\beta_1}{2}).
\end{align*}
The transformations $S_i, \ i=0,1,\dots ,6,$ define a represention of the extended affine Weyl group of type $B_4^{(1)}$, that is, they satisfy the following relations\rm{:\rm}
\begin{align*}
&S_0^2=S_1^2=\dots =S_4^2=S_5^2=S_6^2=1,\\
&(S_0S_1)^2=(S_0S_3)^2=(S_0S_4)^2=(S_1S_3)^2=(S_1S_4)^2=(S_2S_4)^2=1,\\
&(S_0S_2)^3=(S_1S_2)^3=(S_2S_3)^3=1, \ (S_3S_4)^4=1,\\
&S_5(S_0,S_1,S_2,S_3,S_4)=(S_1,S_0,S_2,S_3,S_4)S_5.
\end{align*}
The proof has thus been completed. \qed

\section{Further problems}

For the system of type $G_2^{(1)}$, let us consider a generalization of this system with the same way in Sections 5 and 6.

By using the coupling term $x-z$, we make the following representation:
\begin{align}\label{symmetryG4}
\begin{split}
s_0(*) \rightarrow &(x,y,z,w,t;-\alpha_0,\alpha_1+\alpha_0,\alpha_2,\alpha_3,\alpha_4),\\
s_1(*) \rightarrow &(x+\frac{\alpha_1}{y},y,z,w,t;\alpha_0+\alpha_1,-\alpha_1,\alpha_2+\alpha_1,\alpha_3,\alpha_4),\\
s_2(*) \rightarrow &(x,y-\frac{\alpha_2}{x-z},z,w+\frac{\alpha_2}{x-z},t;\alpha_0,\alpha_1+\alpha_2,-\alpha_2,\alpha_3+\alpha_2,\alpha_4),\\
s_3(*) \rightarrow &(x,y,z+\frac{\alpha_3}{w},w,t;\alpha_0,\alpha_1,\alpha_2+\alpha_3,-\alpha_3,\alpha_4+\alpha_3),\\
s_4(*) \rightarrow &(\sqrt{-1}x,-\sqrt{-1}y,\sqrt{-1}z,-\sqrt{-1}(w-\frac{3\alpha_4}{z}+\frac{t}{z^2}+\frac{1}{2z^3}),-\sqrt{-1}t;\\
&\alpha_0,\alpha_1,\alpha_2,\alpha_3+3\alpha_4,-\alpha_4).
\end{split}
\end{align}

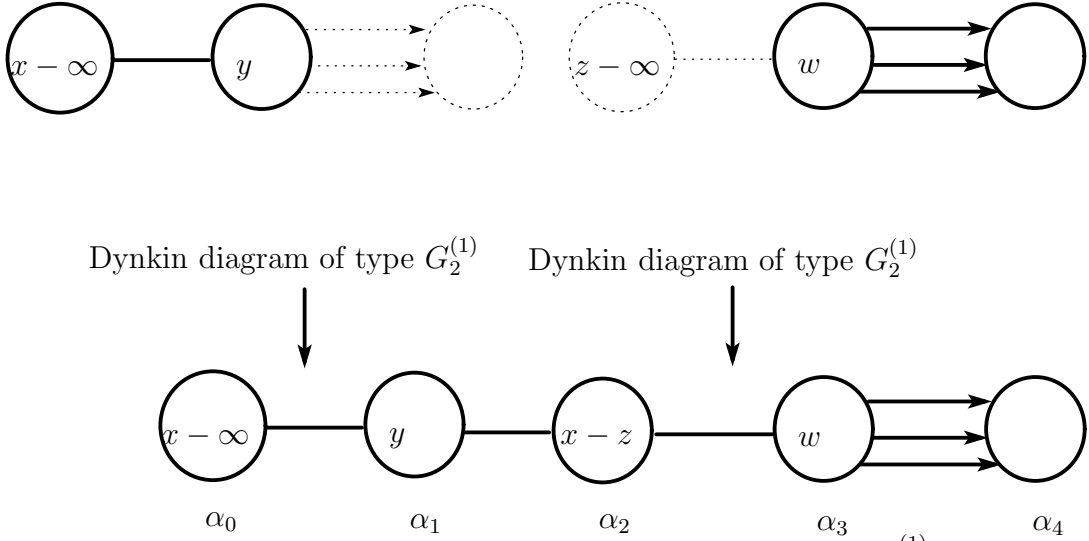
\begin{figure}[ht]
\unitlength 0.1in
\begin{picture}(56.38,25.90)(18.85,-34.50)
%
\special{pn 20}%
\special{ar 2160 1150 275 285  0.0000000 6.2831853}%
%
\special{pn 20}%
\special{ar 3215 1144 257 272  0.0000000 6.2831853}%
%
\special{pn 8}%
\special{ar 4323 1144 258 272  0.0000000 0.0452830}%
\special{ar 4323 1144 258 272  0.1811321 0.2264151}%
\special{ar 4323 1144 258 272  0.3622642 0.4075472}%
\special{ar 4323 1144 258 272  0.5433962 0.5886792}%
\special{ar 4323 1144 258 272  0.7245283 0.7698113}%
\special{ar 4323 1144 258 272  0.9056604 0.9509434}%
\special{ar 4323 1144 258 272  1.0867925 1.1320755}%
\special{ar 4323 1144 258 272  1.2679245 1.3132075}%
\special{ar 4323 1144 258 272  1.4490566 1.4943396}%
\special{ar 4323 1144 258 272  1.6301887 1.6754717}%
\special{ar 4323 1144 258 272  1.8113208 1.8566038}%
\special{ar 4323 1144 258 272  1.9924528 2.0377358}%
\special{ar 4323 1144 258 272  2.1735849 2.2188679}%
\special{ar 4323 1144 258 272  2.3547170 2.4000000}%
\special{ar 4323 1144 258 272  2.5358491 2.5811321}%
\special{ar 4323 1144 258 272  2.7169811 2.7622642}%
\special{ar 4323 1144 258 272  2.8981132 2.9433962}%
\special{ar 4323 1144 258 272  3.0792453 3.1245283}%
\special{ar 4323 1144 258 272  3.2603774 3.3056604}%
\special{ar 4323 1144 258 272  3.4415094 3.4867925}%
\special{ar 4323 1144 258 272  3.6226415 3.6679245}%
\special{ar 4323 1144 258 272  3.8037736 3.8490566}%
\special{ar 4323 1144 258 272  3.9849057 4.0301887}%
\special{ar 4323 1144 258 272  4.1660377 4.2113208}%
\special{ar 4323 1144 258 272  4.3471698 4.3924528}%
\special{ar 4323 1144 258 272  4.5283019 4.5735849}%
\special{ar 4323 1144 258 272  4.7094340 4.7547170}%
\special{ar 4323 1144 258 272  4.8905660 4.9358491}%
\special{ar 4323 1144 258 272  5.0716981 5.1169811}%
\special{ar 4323 1144 258 272  5.2528302 5.2981132}%
\special{ar 4323 1144 258 272  5.4339623 5.4792453}%
\special{ar 4323 1144 258 272  5.6150943 5.6603774}%
\special{ar 4323 1144 258 272  5.7962264 5.8415094}%
\special{ar 4323 1144 258 272  5.9773585 6.0226415}%
\special{ar 4323 1144 258 272  6.1584906 6.2037736}%
\put(23.1000,-22.5000){\makebox(0,0)[lb]{Dynkin diagram of type $G_2^{(1)}$}}%
\put(46.1000,-22.6000){\makebox(0,0)[lb]{Dynkin diagram of type $G_2^{(1)}$}}%
%
\special{pn 20}%
\special{pa 3440 2360}%
\special{pa 3440 2740}%
\special{fp}%
\special{sh 1}%
\special{pa 3440 2740}%
\special{pa 3460 2673}%
\special{pa 3440 2687}%
\special{pa 3420 2673}%
\special{pa 3440 2740}%
\special{fp}%
%
\special{pn 20}%
\special{pa 5680 2350}%
\special{pa 5680 2730}%
\special{fp}%
\special{sh 1}%
\special{pa 5680 2730}%
\special{pa 5700 2663}%
\special{pa 5680 2677}%
\special{pa 5660 2663}%
\special{pa 5680 2730}%
\special{fp}%
\put(19.0000,-12.4000){\makebox(0,0)[lb]{$x-\infty$}}%
\put(30.8100,-12.3400){\makebox(0,0)[lb]{$y$}}%
%
\special{pn 20}%
\special{pa 2440 1160}%
\special{pa 2940 1160}%
\special{fp}%
%
\special{pn 8}%
\special{pa 3450 1000}%
\special{pa 4060 1000}%
\special{dt 0.045}%
\special{sh 1}%
\special{pa 4060 1000}%
\special{pa 3993 980}%
\special{pa 4007 1000}%
\special{pa 3993 1020}%
\special{pa 4060 1000}%
\special{fp}%
%
\special{pn 8}%
\special{pa 3470 1190}%
\special{pa 4040 1190}%
\special{dt 0.045}%
\special{sh 1}%
\special{pa 4040 1190}%
\special{pa 3973 1170}%
\special{pa 3987 1190}%
\special{pa 3973 1210}%
\special{pa 4040 1190}%
\special{fp}%
%
\special{pn 8}%
\special{pa 3410 1330}%
\special{pa 4110 1330}%
\special{dt 0.045}%
\special{sh 1}%
\special{pa 4110 1330}%
\special{pa 4043 1310}%
\special{pa 4057 1330}%
\special{pa 4043 1350}%
\special{pa 4110 1330}%
\special{fp}%
%
\special{pn 8}%
\special{ar 5100 1145 275 285  0.0000000 0.0428571}%
\special{ar 5100 1145 275 285  0.1714286 0.2142857}%
\special{ar 5100 1145 275 285  0.3428571 0.3857143}%
\special{ar 5100 1145 275 285  0.5142857 0.5571429}%
\special{ar 5100 1145 275 285  0.6857143 0.7285714}%
\special{ar 5100 1145 275 285  0.8571429 0.9000000}%
\special{ar 5100 1145 275 285  1.0285714 1.0714286}%
\special{ar 5100 1145 275 285  1.2000000 1.2428571}%
\special{ar 5100 1145 275 285  1.3714286 1.4142857}%
\special{ar 5100 1145 275 285  1.5428571 1.5857143}%
\special{ar 5100 1145 275 285  1.7142857 1.7571429}%
\special{ar 5100 1145 275 285  1.8857143 1.9285714}%
\special{ar 5100 1145 275 285  2.0571429 2.1000000}%
\special{ar 5100 1145 275 285  2.2285714 2.2714286}%
\special{ar 5100 1145 275 285  2.4000000 2.4428571}%
\special{ar 5100 1145 275 285  2.5714286 2.6142857}%
\special{ar 5100 1145 275 285  2.7428571 2.7857143}%
\special{ar 5100 1145 275 285  2.9142857 2.9571429}%
\special{ar 5100 1145 275 285  3.0857143 3.1285714}%
\special{ar 5100 1145 275 285  3.2571429 3.3000000}%
\special{ar 5100 1145 275 285  3.4285714 3.4714286}%
\special{ar 5100 1145 275 285  3.6000000 3.6428571}%
\special{ar 5100 1145 275 285  3.7714286 3.8142857}%
\special{ar 5100 1145 275 285  3.9428571 3.9857143}%
\special{ar 5100 1145 275 285  4.1142857 4.1571429}%
\special{ar 5100 1145 275 285  4.2857143 4.3285714}%
\special{ar 5100 1145 275 285  4.4571429 4.5000000}%
\special{ar 5100 1145 275 285  4.6285714 4.6714286}%
\special{ar 5100 1145 275 285  4.8000000 4.8428571}%
\special{ar 5100 1145 275 285  4.9714286 5.0142857}%
\special{ar 5100 1145 275 285  5.1428571 5.1857143}%
\special{ar 5100 1145 275 285  5.3142857 5.3571429}%
\special{ar 5100 1145 275 285  5.4857143 5.5285714}%
\special{ar 5100 1145 275 285  5.6571429 5.7000000}%
\special{ar 5100 1145 275 285  5.8285714 5.8714286}%
\special{ar 5100 1145 275 285  6.0000000 6.0428571}%
\special{ar 5100 1145 275 285  6.1714286 6.2142857}%
%
\special{pn 20}%
\special{ar 6155 1139 257 272  0.0000000 6.2831853}%
%
\special{pn 20}%
\special{ar 7263 1139 258 272  0.0000000 6.2831853}%
\put(48.5000,-12.4000){\makebox(0,0)[lb]{$z-\infty$}}%
\put(60.2100,-12.2900){\makebox(0,0)[lb]{$w$}}%
%
\special{pn 8}%
\special{pa 5380 1155}%
\special{pa 5880 1155}%
\special{dt 0.045}%
\special{pa 5880 1155}%
\special{pa 5879 1155}%
\special{dt 0.045}%
%
\special{pn 20}%
\special{pa 6390 995}%
\special{pa 7000 995}%
\special{fp}%
\special{sh 1}%
\special{pa 7000 995}%
\special{pa 6933 975}%
\special{pa 6947 995}%
\special{pa 6933 1015}%
\special{pa 7000 995}%
\special{fp}%
%
\special{pn 20}%
\special{pa 6410 1185}%
\special{pa 6980 1185}%
\special{fp}%
\special{sh 1}%
\special{pa 6980 1185}%
\special{pa 6913 1165}%
\special{pa 6927 1185}%
\special{pa 6913 1205}%
\special{pa 6980 1185}%
\special{fp}%
%
\special{pn 20}%
\special{pa 6350 1325}%
\special{pa 7050 1325}%
\special{fp}%
\special{sh 1}%
\special{pa 7050 1325}%
\special{pa 6983 1305}%
\special{pa 6997 1325}%
\special{pa 6983 1345}%
\special{pa 7050 1325}%
\special{fp}%
%
\special{pn 20}%
\special{ar 2960 3075 275 285  0.0000000 6.2831853}%
%
\special{pn 20}%
\special{ar 4015 3069 257 272  0.0000000 6.2831853}%
\put(26.9000,-31.7000){\makebox(0,0)[lb]{$x-\infty$}}%
\put(38.8100,-31.5900){\makebox(0,0)[lb]{$y$}}%
%
\special{pn 20}%
\special{pa 3240 3085}%
\special{pa 3740 3085}%
\special{fp}%
%
\special{pn 20}%
\special{ar 5000 3100 257 272  0.0000000 6.2831853}%
%
\special{pn 20}%
\special{ar 6157 3092 257 272  0.0000000 6.2831853}%
%
\special{pn 20}%
\special{ar 7265 3092 258 272  0.0000000 6.2831853}%
\put(60.2300,-31.8200){\makebox(0,0)[lb]{$w$}}%
%
\special{pn 20}%
\special{pa 6392 2948}%
\special{pa 7002 2948}%
\special{fp}%
\special{sh 1}%
\special{pa 7002 2948}%
\special{pa 6935 2928}%
\special{pa 6949 2948}%
\special{pa 6935 2968}%
\special{pa 7002 2948}%
\special{fp}%
%
\special{pn 20}%
\special{pa 6412 3138}%
\special{pa 6982 3138}%
\special{fp}%
\special{sh 1}%
\special{pa 6982 3138}%
\special{pa 6915 3118}%
\special{pa 6929 3138}%
\special{pa 6915 3158}%
\special{pa 6982 3138}%
\special{fp}%
%
\special{pn 20}%
\special{pa 6352 3278}%
\special{pa 7052 3278}%
\special{fp}%
\special{sh 1}%
\special{pa 7052 3278}%
\special{pa 6985 3258}%
\special{pa 6999 3278}%
\special{pa 6985 3298}%
\special{pa 7052 3278}%
\special{fp}%
%
\special{pn 20}%
\special{pa 4270 3110}%
\special{pa 4720 3110}%
\special{fp}%
%
\special{pn 20}%
\special{pa 5280 3120}%
\special{pa 5890 3120}%
\special{fp}%
\put(47.8000,-31.7000){\makebox(0,0)[lb]{$x-z$}}%
\put(29.2000,-36.0000){\makebox(0,0)[lb]{$\alpha_0$}}%
\put(39.9000,-36.1000){\makebox(0,0)[lb]{$\alpha_1$}}%
\put(49.8000,-36.1000){\makebox(0,0)[lb]{$\alpha_2$}}%
\put(61.2000,-36.2000){\makebox(0,0)[lb]{$\alpha_3$}}%
\put(72.5000,-36.2000){\makebox(0,0)[lb]{$\alpha_4$}}%
\end{picture}%
\label{fig:PS2}
\caption{A generalization of the system of type $G_2^{(1)}$}
\end{figure}

\begin{Lemma}
These transformations $s_i \ (i=0,1,\dots,4)$ satisfy the following relations{\rm : \rm}
\begin{align}
\begin{split}
&s_i^2=1, \quad (s_0s_2)^2=(s_0s_3)^2=(s_0s_4)^2=(s_1s_3)^2=(s_1s_4)^2=(s_2s_4)^2=1,\\
&(s_0s_1)^3=(s_1s_2)^3=(s_2s_3)^3=1, \quad (s_3s_4)^6=1.
\end{split}
\end{align}
\end{Lemma}

Let us make the holomorphy conditions associated with \eqref{symmetryG4}.
\begin{align}
\begin{split}
r_0:&x_0=\frac{1}{x}, \quad y_0=-(yx+\alpha_0+\alpha_1)x, \quad z_0=z, \quad w_0=w,\\
r_1:&x_1=\frac{1}{x}, \quad y_1=-(yx+\alpha_1)x, \quad z_1=z, \quad w_1=w,\\
r_2:&x_2=-((x-z)y-\alpha_2)y, \quad y_2=\frac{1}{y}, \quad z_2=z, \quad w_2=w+y,\\
r_3:&x_3=x, \quad y_3=y, \quad z_3=\frac{1}{z}, \quad w_3=-(wz+\alpha_3)z,\\
r_4:&x_4=x, \quad y_4=y, \quad z_4=z, \quad w_4=w-\frac{3\alpha_4}{z}+\frac{t}{z^2}+\frac{1}{2z^3}.
\end{split}
\end{align}

\begin{Problem}
Let us consider a polynomial Hamiltonian system with Hamiltonian $H \in {\Bbb C}(t)[x,y,z,w]$. We assume that

$(A1)$ $deg(H)=5$ with respect to $x,y,z,w$.

$(A2)$ This system becomes again a polynomial Hamiltonian system in each coordinate $r_i \ (i=0,1, \dots ,4)$.
\end{Problem}
It is a pity that we can not find a polynomial Hamiltonian system satisfying the assumptions $(A1)$ and $(A2)$.

It is still an open question whether we can find a generalization of the system  of type $G_2^{(1)}$.

For the system of type $A_2^{(2)}$, let us consider a generalization of this system with the same way in Sections 5 and 6.

By using the coupling term $x-z$, we make the following representation:
\begin{align}\label{symmetryA3}
\begin{split}
s_0(*) \rightarrow &(x+\frac{\alpha_0}{y},y,z,w,t;-\alpha_0,\alpha_1+\alpha_0,\alpha_2,\alpha_3),\\
s_1(*) \rightarrow &(x,y-\frac{\alpha_1}{x-z},z,w+\frac{\alpha_1}{x-z},t;\alpha_0+\alpha_1,-\alpha_1,\alpha_2+\alpha_1,\alpha_3),\\
s_2(*) \rightarrow &(x,y,z+\frac{\alpha_2}{w},w,t;\alpha_0,\alpha_1+\alpha_2,-\alpha_2,\alpha_3+4\alpha_2),\\
s_3(*) \rightarrow &(-x,-y,-z,-w+\frac{\alpha_3}{z}-\frac{t}{z^2}-\frac{2}{z^4},t;\alpha_0,\alpha_1,\alpha_2+\alpha_3,-\alpha_3).
\end{split}
\end{align}

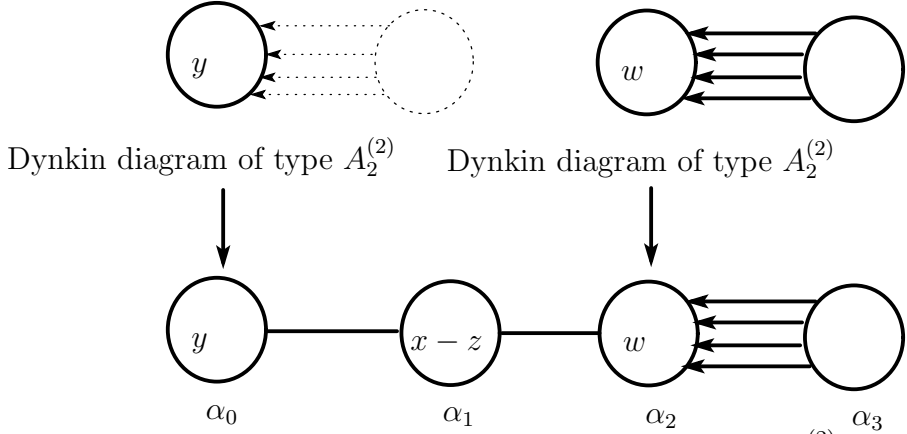
\begin{figure}[h]
\unitlength 0.1in
\begin{picture}(46.89,20.40)(23.10,-34.20)
%
\special{pn 20}%
\special{ar 3407 1652 257 272  0.0000000 6.2831853}%
\put(23.1000,-22.5000){\makebox(0,0)[lb]{Dynkin diagram of type $A_2^{(2)}$}}%
\put(46.1000,-22.6000){\makebox(0,0)[lb]{Dynkin diagram of type $A_2^{(2)}$}}%
%
\special{pn 20}%
\special{pa 3440 2360}%
\special{pa 3440 2740}%
\special{fp}%
\special{sh 1}%
\special{pa 3440 2740}%
\special{pa 3460 2673}%
\special{pa 3440 2687}%
\special{pa 3420 2673}%
\special{pa 3440 2740}%
\special{fp}%
%
\special{pn 20}%
\special{pa 5680 2350}%
\special{pa 5680 2730}%
\special{fp}%
\special{sh 1}%
\special{pa 5680 2730}%
\special{pa 5700 2663}%
\special{pa 5680 2677}%
\special{pa 5660 2663}%
\special{pa 5680 2730}%
\special{fp}%
\put(32.7300,-17.4200){\makebox(0,0)[lb]{$y$}}%
%
\special{pn 8}%
\special{ar 4492 1688 257 272  0.0000000 0.0453686}%
\special{ar 4492 1688 257 272  0.1814745 0.2268431}%
\special{ar 4492 1688 257 272  0.3629490 0.4083176}%
\special{ar 4492 1688 257 272  0.5444234 0.5897921}%
\special{ar 4492 1688 257 272  0.7258979 0.7712665}%
\special{ar 4492 1688 257 272  0.9073724 0.9527410}%
\special{ar 4492 1688 257 272  1.0888469 1.1342155}%
\special{ar 4492 1688 257 272  1.2703214 1.3156900}%
\special{ar 4492 1688 257 272  1.4517958 1.4971645}%
\special{ar 4492 1688 257 272  1.6332703 1.6786389}%
\special{ar 4492 1688 257 272  1.8147448 1.8601134}%
\special{ar 4492 1688 257 272  1.9962193 2.0415879}%
\special{ar 4492 1688 257 272  2.1776938 2.2230624}%
\special{ar 4492 1688 257 272  2.3591682 2.4045369}%
\special{ar 4492 1688 257 272  2.5406427 2.5860113}%
\special{ar 4492 1688 257 272  2.7221172 2.7674858}%
\special{ar 4492 1688 257 272  2.9035917 2.9489603}%
\special{ar 4492 1688 257 272  3.0850662 3.1304348}%
\special{ar 4492 1688 257 272  3.2665406 3.3119093}%
\special{ar 4492 1688 257 272  3.4480151 3.4933837}%
\special{ar 4492 1688 257 272  3.6294896 3.6748582}%
\special{ar 4492 1688 257 272  3.8109641 3.8563327}%
\special{ar 4492 1688 257 272  3.9924386 4.0378072}%
\special{ar 4492 1688 257 272  4.1739130 4.2192817}%
\special{ar 4492 1688 257 272  4.3553875 4.4007561}%
\special{ar 4492 1688 257 272  4.5368620 4.5822306}%
\special{ar 4492 1688 257 272  4.7183365 4.7637051}%
\special{ar 4492 1688 257 272  4.8998110 4.9451796}%
\special{ar 4492 1688 257 272  5.0812854 5.1266541}%
\special{ar 4492 1688 257 272  5.2627599 5.3081285}%
\special{ar 4492 1688 257 272  5.4442344 5.4896030}%
\special{ar 4492 1688 257 272  5.6257089 5.6710775}%
\special{ar 4492 1688 257 272  5.8071834 5.8525520}%
\special{ar 4492 1688 257 272  5.9886578 6.0340265}%
\special{ar 4492 1688 257 272  6.1701323 6.2155009}%
%
\special{pn 20}%
\special{ar 5657 1692 257 272  0.0000000 6.2831853}%
\put(55.2300,-17.8200){\makebox(0,0)[lb]{$w$}}%
%
\special{pn 20}%
\special{ar 6742 1728 257 272  0.0000000 6.2831853}%
%
\special{pn 20}%
\special{ar 3407 3092 257 272  0.0000000 6.2831853}%
\put(32.7300,-31.8200){\makebox(0,0)[lb]{$y$}}%
%
\special{pn 20}%
\special{ar 4630 3110 257 272  0.0000000 6.2831853}%
%
\special{pn 20}%
\special{ar 5667 3112 257 272  0.0000000 6.2831853}%
\put(55.3300,-32.0200){\makebox(0,0)[lb]{$w$}}%
%
\special{pn 20}%
\special{pa 3670 3100}%
\special{pa 4350 3100}%
\special{fp}%
%
\special{pn 20}%
\special{pa 4900 3110}%
\special{pa 5410 3110}%
\special{fp}%
\put(44.2000,-32.0000){\makebox(0,0)[lb]{$x-z$}}%
\put(33.5000,-35.6000){\makebox(0,0)[lb]{$\alpha_0$}}%
\put(45.9000,-35.7000){\makebox(0,0)[lb]{$\alpha_1$}}%
\put(56.5000,-35.8000){\makebox(0,0)[lb]{$\alpha_2$}}%
\put(67.3000,-35.9000){\makebox(0,0)[lb]{$\alpha_3$}}%
%
\special{pn 8}%
\special{pa 4310 1500}%
\special{pa 3640 1500}%
\special{dt 0.045}%
\special{sh 1}%
\special{pa 3640 1500}%
\special{pa 3707 1520}%
\special{pa 3693 1500}%
\special{pa 3707 1480}%
\special{pa 3640 1500}%
\special{fp}%
%
\special{pn 8}%
\special{pa 4220 1650}%
\special{pa 3670 1650}%
\special{dt 0.045}%
\special{sh 1}%
\special{pa 3670 1650}%
\special{pa 3737 1670}%
\special{pa 3723 1650}%
\special{pa 3737 1630}%
\special{pa 3670 1650}%
\special{fp}%
%
\special{pn 8}%
\special{pa 4240 1770}%
\special{pa 3640 1770}%
\special{dt 0.045}%
\special{sh 1}%
\special{pa 3640 1770}%
\special{pa 3707 1790}%
\special{pa 3693 1770}%
\special{pa 3707 1750}%
\special{pa 3640 1770}%
\special{fp}%
%
\special{pn 8}%
\special{pa 4280 1860}%
\special{pa 3590 1860}%
\special{dt 0.045}%
\special{sh 1}%
\special{pa 3590 1860}%
\special{pa 3657 1880}%
\special{pa 3643 1860}%
\special{pa 3657 1840}%
\special{pa 3590 1860}%
\special{fp}%
%
\special{pn 20}%
\special{pa 6540 1540}%
\special{pa 5890 1540}%
\special{fp}%
\special{sh 1}%
\special{pa 5890 1540}%
\special{pa 5957 1560}%
\special{pa 5943 1540}%
\special{pa 5957 1520}%
\special{pa 5890 1540}%
\special{fp}%
%
\special{pn 20}%
\special{pa 6470 1650}%
\special{pa 5930 1650}%
\special{fp}%
\special{sh 1}%
\special{pa 5930 1650}%
\special{pa 5997 1670}%
\special{pa 5983 1650}%
\special{pa 5997 1630}%
\special{pa 5930 1650}%
\special{fp}%
%
\special{pn 20}%
\special{pa 6460 1770}%
\special{pa 5920 1770}%
\special{fp}%
\special{sh 1}%
\special{pa 5920 1770}%
\special{pa 5987 1790}%
\special{pa 5973 1770}%
\special{pa 5987 1750}%
\special{pa 5920 1770}%
\special{fp}%
%
\special{pn 20}%
\special{pa 6510 1880}%
\special{pa 5870 1880}%
\special{fp}%
\special{sh 1}%
\special{pa 5870 1880}%
\special{pa 5937 1900}%
\special{pa 5923 1880}%
\special{pa 5937 1860}%
\special{pa 5870 1880}%
\special{fp}%
%
\special{pn 20}%
\special{ar 6742 3132 257 272  0.0000000 6.2831853}%
%
\special{pn 20}%
\special{pa 6540 2944}%
\special{pa 5890 2944}%
\special{fp}%
\special{sh 1}%
\special{pa 5890 2944}%
\special{pa 5957 2964}%
\special{pa 5943 2944}%
\special{pa 5957 2924}%
\special{pa 5890 2944}%
\special{fp}%
%
\special{pn 20}%
\special{pa 6470 3054}%
\special{pa 5930 3054}%
\special{fp}%
\special{sh 1}%
\special{pa 5930 3054}%
\special{pa 5997 3074}%
\special{pa 5983 3054}%
\special{pa 5997 3034}%
\special{pa 5930 3054}%
\special{fp}%
%
\special{pn 20}%
\special{pa 6460 3174}%
\special{pa 5920 3174}%
\special{fp}%
\special{sh 1}%
\special{pa 5920 3174}%
\special{pa 5987 3194}%
\special{pa 5973 3174}%
\special{pa 5987 3154}%
\special{pa 5920 3174}%
\special{fp}%
%
\special{pn 20}%
\special{pa 6510 3284}%
\special{pa 5870 3284}%
\special{fp}%
\special{sh 1}%
\special{pa 5870 3284}%
\special{pa 5937 3304}%
\special{pa 5923 3284}%
\special{pa 5937 3264}%
\special{pa 5870 3284}%
\special{fp}%
\end{picture}%
\label{fig:PS3}
\caption{A generalization of the system of type $A_2^{(2)}$}
\end{figure}

\begin{Lemma}
These transformations $s_i \ (i=0,1,2,3)$ satisfy the following relations{\rm : \rm}
\begin{align}
s_i^2=1, \quad (s_0s_2)^2=(s_0s_3)^2=(s_1s_3)^2=1, \quad (s_0s_1)^3=(s_1s_2)^3=1.
\end{align}
\end{Lemma}

Let us make the holomorphy conditions associated with \eqref{symmetryA3}.
\begin{align}
\begin{split}
r_0:&x_0=\frac{1}{x}, \quad y_0=-(yx+\alpha_0)x, \quad z_0=z, \quad w_0=w,\\
r_1:&x_1=-((x-z)y-\alpha_1)y, \quad y_1=\frac{1}{y}, \quad z_1=z, \quad w_1=w+y,\\
r_2:&x_2=x, \quad y_2=y, \quad z_2=\frac{1}{z}, \quad w_2=-(wz+\alpha_2)z,\\
r_3:&x_3=x, \quad y_3=y, \quad z_3=z, \quad w_3=w-\frac{\alpha_3}{z}+\frac{t}{z^2}+\frac{2}{z^4}.
\end{split}
\end{align}

\begin{Problem}
Let us consider a polynomial Hamiltonian system with Hamiltonian $H \in {\Bbb C}(t)[x,y,z,w]$. We assume that

$(A1)$ $deg(H)=6$ with respect to $x,y,z,w$.

$(A2)$ This system becomes again a polynomial Hamiltonian system in each coordinate $r_i \ (i=0,1,2,3)$.
\end{Problem}
It is a pity that we can not find a polynomial Hamiltonian system satisfying the assumptions $(A1)$ and $(A2)$.

It is still an open question whether we can find a generalization of the system  of type $A_2^{(2)}$.

{\it Acknowledgements.} The author would like to thank K. Fuji, H. Kawamuko, T. Masuda, Y Ohta, W. Rossman, T. Suzuki, N. Tahara, K. Takano, T. Tsuda and Y. Yamada  for useful suggestions and discussions.

\begin{flushleft}{\it
Graduate School of Mathematical Sciences,\\
University of Tokyo,\\
3-8-1 Komaba,Megro-ku,153-8914 Tokyo, Japan\\
\it}
{\rm
E-mail address: sasano@ms.u-tokyo.ac.jp
\rm}
\end{flushleft}

\end{document}